\input amstex
\input amsppt.sty
\NoBlackBoxes
\nologo
\topmatter
\title A note on the factorization theorem \linebreak
of toric birational maps\linebreak
after Morelli \linebreak
and \linebreak
its toroidal extension
\endtitle
\author
Dan Abramovich, \footnote{The first author is partially supported by NSF grant
DMS-9700520 and by an Alfred P. Sloan research fellowship.\hfill\hfill}
Kenji
Matsuki
\footnote{The second author is partially supported by NSA grant
MDA904-96-1-0008.\hfill\hfill}
and Suliman Rashid \footnote{The third author is partially supported by the
Purdue Research Foundation.\hfill\hfill}
\endauthor
\rightheadtext{A note on the factorization theorem}
\leftheadtext{D. Abramovich, K. Matsuki, S. Rashid}
\rightheadtext{A note on the factorization theorem}
\leftheadtext{D. Abramovich, K. Matsuki, S. Rashid}

\abstract Building upon the work of [Morelli1,2], we give a
coherent presentation of Morelli's algorithm for the weak and
strong factorization of toric birational maps.  We also discuss
its toroidal extension, which plays a crucial role in the recent
solutions [W{\l}odarczyk2,3, Abramovich-Karu-Matsuki-W{\l}odarczyk] of the weak
factorization conjecture of general birational maps.
\endabstract

\endtopmatter

\document

\footnote""{1991 {\it Mathematics Subject Classification}. Primary 14M25;
Secondary 14E05}

$$\bold{Table\ of\ Contents}$$

\vskip.1in

\S 0. Introduction

\vskip.1in

\S 1. Basic Ideas
 
\vskip.1in

\S 2. Cobordism

\vskip.1in

\S 3. Circuits and Bistellar Operations

\vskip.1in

\S 4. Collapsibility

\vskip.1in

\S 5. $\pi$-Desingularization

\vskip.1in

\S 6. The Weak Factorization Theorem

\vskip.1in

\S 7. The Strong Factorization Theorem

\vskip.1in

\S 8. The Toroidal Case

\newpage

\S 0. Introduction

\vskip.1in

This paper is a result of series of seminars held by the authors during
the summer of 1997 and continued from then on, toward a thorough understanding
of the following weak and strong factorization theorem of toric birational
maps by [Morelli1] (cf.[W{\l}odarczyk1]).  

\proclaim{Theorem 0.1 (Factorization Theorem for Toric Birational Maps)} Every
proper and equivariant birational map
$f:X_{\Delta}
\dashrightarrow X_{\Delta'}$ (``proper" in the sense of [Iitaka]) between
two nonsingular toric varieties can be factored into a sequence of blowups and
blowdowns with smooth centers which are the closures of orbits.  
\endproclaim

If we allow the sequence
to consist of blowups and blowdowns in any order, then the factorization is
called $\bold{weak}$.

If we insist on the sequence to consist only of blowups immediately followed by
blowdowns, then the factorization is called $\bold{strong}$.  

\vskip.1in

Our purpose is two-fold.  

The first is to
give a coherent presentation of the proof in [Morelli1] both for the weak
factorization and the strong factorization, modifying some discrepancies found
by King [King2] and by the authors in the due course of the seminars checking
the original arguments.  Most of these discrepancies are minor, except for a
couple of essential points in the process of
$\pi$-desingularization and in the process of showing that the weak
factorization implies the strong factorization.  \footnote{As of Jan.1998 we
learned from Prof. Fulton that Morelli himself offers correction in his homepage
[Morelli2] to the discrepancies in the process of $\pi$-desingularization found
by King.  We still need some clarification, as is presented in this paper, to
understand the correction.  We thank Prof. Morelli for guiding us toward a
better understanding through private communication.}  It is a mere attempt to
see in a transparent way the beautiful and brilliant original ideas of
[Morelli1,2] by sweeping dust off the surface.  

The second is the generalization to the
toroidal case, whose details are worked out as a part of the Ph.D. thesis of
the third author.  Though it may be said that the toroidal generalization is
straightforward and even implicit in the original papers [Morelli1,2]
(cf.[W{\l}odarczyk1]), we would like to emphasize its importance in a more
far-reaching problem formulated as below, with a view toward its application to
the factorization problem of general birational maps.\footnote{Recently two
independent proofs have appeared for the weak factorization conjecture of
general birational maps, one by [W{\l}odarczyk3] another by
[Abramovich-Karu-Matsuki-W{\l}odarczyk].  (Both proofs are based upon the theory
of birational cobordism of [W{\l}odarczyk2], which is inspired by the
combinatorial cobordism of [Morelli1] discussed in \S 2 of this paper.)  The
former uses the algorithm for
$\pi$-desingularization, while the latter uses the strong factorization of
toroidal birational maps directly in their proofs.  Thus the importance of the
toroidal extension has only increased, as well as the need for a clear coherent
presentation for the $\pi$-desingularization process.  The toroidalization
conjecture and the strong factorization conjecture remain open.}

In a most naive way the
``far-reaching" problem can be stated as follows: Let $f:X \rightarrow Y$ be a
morphism (one may put the condition ``with connected fibers" if one wishes)
between nonsingular complete (or projective) varieties.  By replacing
$X$ and $Y$ with their modifications $X'$ and $Y'$, how ``NICE" can one make the
morphism $f':X' \rightarrow Y'$ ?  
$$\CD
X @<<< X' \\
@VVfV @VVf'V \\
Y @<<< Y' \\
\endCD$$
Depending upon how we interpret the word
``NICE" mathematically and what restrictions we put on the modifications, we get
the corresponding interesting questions such as semistable reduction (when the
morphism $f'$ is ``NICE" if every fiber is reduced with only simple normal
crossings, $\dim Y = 1$ and the modifications for $Y$ are restricted to finite
morphisms while the modifications for $X$ are restricted to smooth blowups
after base change), resolution of hypersurface singularities (when the morphism
$f'$ is ``NICE" if every fiber has only simple normal crossings, this time not
necessarily reduced, $\dim Y = 1$ and no modification for $Y$ and only smooth
blowups are allowed for
$X$).  When
$f$ is birational and we require $f'$ to be an isomorphism in order for it to be
``NICE", restricting the modifications of $X$ and $Y$ to be smooth blowups, we
obtain the long standing (and perhaps notorious) strong factorization problem
for general birational morphisms (cf.[Hironaka]).

Our interpretation is that we put ``toroidal" for the word ``NICE" and
restrict the modifications of $X$ and $Y$ to be only smooth blowups.  

\proclaim{Conjecture 0.2 (Toroidalization Conjecture)} Let $f:X \rightarrow Y$
be a morphism between nonsingular complete varieties.  Then there exist
sequences of blowups with smooth centers for $X$ and $Y$ so that the induced
morphism
$f':X'
\rightarrow Y'$ is toroidal:

$$\CD
X @<\text{smooth\ blowups}<< X' \\
@VVfV @VVf'V \hskip.2in \text{toroidal}\\
Y @<\text{smooth\ blowups}<< Y'. \\
\endCD$$

\endproclaim

The conjecture is closely related to the recent work of [Abramovich-Karu], which
introduces the notion of ``toroidal" morphisms explicitly for the first time,
though implicitly it can be recognized in [Kempf-Knudsen-Mumford-SaintDonat]. 
By only requiring ``NICE" morphisms to be toroidal instead of being
isomorphisms, we can start dealing not only with birational morphisms but also
with fibering morphisms between varieties of different dimensions.  This seems
to give us more freedom to seek some inductional structure.  Actually we expect
that the powerful inductive method of [Bierstone-Milman] for the canonical
resolution of singularities, proceeding from the hypersurface case with only one
defining equation to the general case with several defining equations
through the ingeneous use of invariants, should be modified to be applied to
our toroidalization problem, proceeding similarly from the case $\dim Y = 1$ to
the general case $\dim Y > 1$. 

\vskip.1in

This interpretation not only generalizes the statement of the classical
factorization problem but also gives the following approach to it:

\vskip.1in

\proclaim{Expectation 0.3 (A Conjectural Approach to the Strong Factorization
Problem via Toroidalization)} Given a birational morphism $f:X \rightarrow Y$
between nonsingular complete varieties,

\vskip.1in

(I) make it ``toroidal" $f':X' \rightarrow Y'$ modifying $X$ and $Y$ into
$X'$ and $Y'$ by blowing up along smooth centers via some Bierstone-Milman type
argument, 

and then

(II) factor the toroidal birational morphism $f':X' \rightarrow Y'$ into
(equivariant) smooth blowups and blowdowns by applying the toroidal version of
the method of [Morelli1,2] (or [W{\l}odarczyk1]). 

\endproclaim

The toroidalization conjecture and the strong factorization of toroidal
birational morphisms would imply the strong factorization of general
birational maps between nonsingular complete varieties.

\vskip.1in

This line of ideas came up in our conversation as
a day-dreaming inspired by [DeJong], only to find out later that an almost
identical approach was already presented in [King1] and has been pursued by him
in reality, who has (privately) announced the affirmative solution to the
toroidalization conjecture in the case $\dim X = 3$.  Actually our formulation
above follows his presentation in [King1].  He has also read [Morelli1]
carefully and his correspondence with Morelli himself was kindly communicated
to us by Bierstone.  We thank both professors for their generosity sharing
their ideas with us and our indebtedness to them is both explicitly and
implicitly clear as well as to the original papers [Morelli1,2] and
[W{\l}odarczyk1].  Another big inspiration for the factorization problem comes
from the recent result of [Cutkosky1], which affirmatively solves the local
factorization problem in dimension 3 using valuation theory.  We thank
Prof. Cutkosky for kindly teaching us his method using valuation theory via
preprints and private conversations.  In response, we communicated to him our
idea above for the global factorization, which turns out to be very similar
to the idea of [Christensen] toward the local factorization problem:

\vskip.1in

(I) First ``monomialize" the given local birational morphism via valuation
theory ([Cutkosky2] uses the word ``monomialization", which is nothing but
``toroidalization" in the local case.), then

(II) factor the local monomial birational morphism (which is a
toroidal birational morphism). 

\vskip.1in

[Cutkosky2,3] achieves the local factorization in arbitrary dimension along
this line of ideas, extending his method using valuation theory.
\footnote{After monomializing a birational map in (I), which is the most
subtle and difficult part, [Cutkosky2] refers to the results of Morelli in
(II).  The first version of [Cutkosky3] factors the monomialized map in his
own algorithm in (II) avoiding the use of results of Morelli, which
were found to contain discrepancies at the time.  The second version of
[Cutkosky3], upon our communication, uses the strong factorization theorem of
this paper by Morelli in (II) and hence provides the strong factorization
theorem in the local case.} 

\vskip.1in

We remark that [Reid3] gives factorization of toric birational maps into
extremal divisorial contractions and flips by establishing the Minimal Model
Program for toric varieties in arbitrary dimension.   The Minimal Model
Program in general, also known as the Mori Program, is only established in
dimension 3 (cf.[Mori1,2, Kawamata1,2,3, Koll\'ar, Reid1,2, Shokurov]).  We
also remark that recently a new algorithm called the Sarkisov Program has
emerged (cf.[Sarkisov, Reid4]) to factor birational maps among uniruled
varieties.  Though it is only established in dimension 3 in general
(cf.[Corti]), the toric case is rather straightforward in arbitrary dimension
(cf.[Matsuki]).  We do not know of a way to solve the classical factorization
problem into smooth blowups and blowdowns using such factorizations as above.

\vskip.1in

Our organization, as being a note to [Morelli1,2], follows exactly the
structure of the original paper [Morelli1,2] with one last section on the
toroidal case added.  The content of each section is outlined at the end of
\S 1, where we explain the main ideas of Morelli.

\vskip.1in

Our hearty thanks go to Prof. Oda for giving us invaluable suggestions at many
critical points of the paper.  We thank the referee for a very careful reading
of the first draft of the paper and for providing us with meticulous and
constructive comments.

\newpage

\S 1. Basic Ideas

\vskip.1in

The purpose of this section is to present the basic ideas of the brilliant
solution of [Morelli1,2] (see also [W{\l}odarczyk1]) to the following
conjecture of Miyake and Oda (cf. [Oda1]).  We follow the usual notation and
terminology  concerning the toric varieties
$X_{\Delta}$ and their corresponding fans $\Delta$, as presented in
[Danilov, Fulton, Oda2].

\proclaim{Conjecture 1.1 (Weak and Strong Factorization of Toric Birational
Maps by Miyake and Oda)} Every proper and equivariant birational map
$f:X_{\Delta}
\dashrightarrow X_{\Delta'}$ (``proper" in the sense of [Iitaka]) between
two nonsingular toric varieties can be factored into a sequence of blowups and
blowdowns with smooth centers which are the closures of orbits.  

If we allow the sequence
to consist of blowups and blowdowns in any order, then the factorization is
called $\bold{weak}$.

If we insist on the sequence to consist only of blowups immediately followed by
blowdowns, then the factorization is called $\bold{strong}$.  

In short, a toric birational map admits not only a weak factorization but also
a strong factorization.
\endproclaim
  
\vskip.1in

As the toric varieties $X_{\Delta}$ correspond to the fans $\Delta$ in $N_{\Bbb
Q} = N \otimes {\Bbb Q}$, where $N$ is the lattice of one-parameter subgroups
of the torus, and blowups to the smooth star subdivisions of $\Delta$, we can
reformulate the above conjecture in the following purely combinatorial language:

\proclaim{Conjecture 1.2 (= Conjecture 1.1 in terms of Fans)} Let
$\Delta$ and
$\Delta'$ be two nonsingular fans in $N_{\Bbb Q}$ with the same support.  Then
there is a sequence of smooth star subdivisions and inverse operations called
smooth star assemblings starting from
$\Delta$ and ending with $\Delta'$. 

If we allow the sequence
to consist of smooth star subdivisions and smooth star assemblings in any order,
then the factorization is called weak.

If we insist on the sequence to consist only of smooth star subdivisions
immediately followed by smooth star assemblings, then the factorization is
called strong.  
\endproclaim

In order to understand Morelli's strategy toward the solution of Conjecture
1.1, we look at the following simple example.

\proclaim{Example 1.3}\endproclaim

We take two fans $\Delta$ and $\Delta'$ to consist of the maximal cones in
$N_{\Bbb Q}
\cong {\Bbb Z}^3 \otimes {\Bbb Q}$
$$\align
\Delta &= \{\gamma_{123},\gamma_{124}\} \\
\Delta' &= \{\gamma_{134},\gamma_{234}\} \\
\text{where\ }& \gamma_{ijk} = \langle v_i,v_j,v_k \rangle \text{\ with\ }\\
 v_1 &= (1,0,0), v_2 = (0,1,0), v_3 = (0,0,1), v_4 = (1,1,-1).\\
\endalign$$

Then we observe that by taking the common refinement ${\tilde \Delta}$ of
$\Delta$ and $\Delta'$, subdivided by the vector
$$v_1 + v_2 = v_3 + v_4,$$
the transformation from $\Delta$ to $\Delta'$ can be factored into a
smooth star subdivision immediately followed by a smooth star assembling
$$\Delta \leftarrow {\tilde \Delta} \rightarrow \Delta'$$
as asserted by Conjecture 1.2.

\vskip.1in

Morelli's great idea is to incorporate all the information of this
factorization into a ``cobordism" $\Sigma$, a fan in the vector space
$N_{\Bbb Q}^+ = N_{\Bbb Q} \oplus {\Bbb Q}$ of one dimensional higher, with
its lower face
$\partial_-\Sigma$ being $\Delta$ and its upper face
$\partial_+\Sigma$ being $\Delta'$.  Namely, we take

$$\align
\Sigma &= \{\sigma \text{\ and\ its\ proper\ faces}\} \subset N_{\Bbb Q}^+ =
N_{\Bbb Q}
\oplus {\Bbb Q} \\
\text{where\ }&\sigma = \langle \rho_1,\rho_2,\rho_3,\rho_4 \rangle \text{\
with\ }\\
\rho_1 &= (v_1,0), \rho_2 = (v_2,0), \rho_3 = (v_3,1), \rho_4 =
(v_4,1)\\
\endalign$$
and where the projection is denoted by
$$\pi:N_{\Bbb Q}^+ = N_{\Bbb Q} \oplus {\Bbb Q} \rightarrow N_{\Bbb Q}.$$
The lower face
$$\partial_-\Sigma = \{\langle
\rho_1,\rho_2,\rho_3 \rangle,\langle \rho_1,\rho_2,\rho_4 \rangle\}$$ maps
isomorphically onto
$\Delta$ by the projection $\pi$ and so does the upper face
$$\partial_+\Sigma = \{\langle \rho_1,\rho_3,\rho_4
\rangle,\langle \rho_2,\rho_3,\rho_4 \rangle\}$$ isomorphically onto $\Delta'$.

Moreover, since $\sigma$ does not map isomorphically onto its image by
$\pi$, i.e., since $\sigma$ is $\pi$-dependent, we have the linear relation,
unique up to scalar multiple, among the the primitive vectors
$v_i$ of the projections of the generators
$\rho_i$ of $\sigma$
$$v_1 + v_2 - v_3 - v_4 = 0.$$
From this linear relation, we can read off the point 
$$v_1 + v_2 = v_3 + v_4$$
by which we have to
subdivide $\Delta$ and $\Delta'$ to reach the common refinement ${\tilde
\Delta}$.

In short, we can realize the factorization from constructing the cobordism.

\vskip.1in

We can summarize Morelli's idea, demonstrated by the above example, in the
following.

\proclaim{Basic Idea 1.4 (Morelli's Idea for Factorization)} Let $\Delta$ and
$\Delta'$ be two nonsingular fans in $N_{\Bbb Q}$ with the same support.  Then
we can realize the (weak) factorization by constructing a cobordism $\Sigma$, a
simplicial fan consisting of $\pi$-strongly convex cones (See \S 3 for the
precise definition.) in
$N_{\Bbb Q}^+ = N_{\Bbb Q}
\oplus {\Bbb Q}$ such that

$(1.4.1)$ the lower face $\partial_-\Sigma$ and upper face
$\partial_+\Sigma$ of
$\Sigma$ map isomorphically onto
$\Delta$ and $\Delta'$ by the projection $\pi$
$$\pi:\partial_-\Sigma \overset{\sim}\to{\rightarrow}
\Delta, \hskip.1in \pi:\partial_+\Sigma \overset{\sim}\to{\rightarrow}
\Delta',$$

$(1.4.2)$ $\Sigma$ is $\pi$-nonsingular (See \S 3 for the precise
definition.),

$(1.4.3)$ $\Sigma$ is collapsible (See \S 4 for the precise definition.).
\endproclaim

In fact, let $\sigma$ be a minimal simplex in $\Sigma$ which is
$\pi$-dependent.  (We call such simplex $\sigma$ a circuit.)  If $\sigma$ is
generated by the extremal rays $\rho_i$
$$\sigma = \langle \rho_1, \rho_2, \cdot\cdot\cdot, \rho_k \rangle,$$
then we have the linear relation among the primitive vectors $v_i$ of the
projections of the generators $\rho_i$
$$\Sigma_{i = 1}^k r_iv_i = 0.$$
Now by the $\pi$-nonsingularity of $\Sigma$ and minimality of $\sigma$, it
follows that we may assume that all the coefficients $r_i$ are either $+1$ or
$-1$ after rescaling (cf. Theorem 3.2).  Thus after renumbering the
$v_i$, we may assume that the linear relation is given by
$$v_1 + v_2 + \cdot\cdot\cdot + v_l - v_{l+1} - \cdot\cdot\cdot - v_k = 0.$$ 
We then observe that by taking the common refinement subdivided by the vector
$$v_1 + v_2 + \cdot\cdot\cdot + v_l = v_{l+1} + \cdot\cdot\cdot + v_k$$
the transformation from $\partial_-\sigma$ to $\partial_+\sigma$ can be
factored into a smooth star subdivision of $\partial_-\sigma$ immediately
followed by a smooth star assembling into $\partial_+\sigma$.

Or more generally, we obtain the factorization
between the lower face $\partial_-\overline{\roman{Star}(\sigma)}$ and upper
face
$\partial_+\overline{\roman{Star}(\sigma)}$ of the closed star
$\overline{\roman{Star}(\sigma)}$ of $\sigma$, where
$$\overline{\roman{Star}(\sigma)} = \{\zeta \in \Sigma;\zeta \subset \eta
\supset
\sigma \text{\ for\ some\ cone\ }\eta \in \Sigma\}.$$
The $\pi$-nonsingularity also guarantees that the $\pi$-projections of all the
lower and upper faces and the common refinement obtained through the star
subdivisions are nonsingular and the star subdivisions are smooth.

This achieves the (weak) factorization for $\overline{\roman{Star}(\sigma)}$ for
one circuit $\sigma$ of $\Sigma$.  In order to achieve the (weak) factorization
for the entire $\Sigma = \cup \overline{\roman{Star}(\sigma)}$, where the union
is taken over all the circuits $\sigma$ in $\Sigma$, we have to coordinate the
way we take the (weak) factorizatons for all the circuits.  This is
done by requiring the collapsibility of the cobordism
$\Sigma$.

\vskip.1in

In \S 2, we construct a cobordism $\Sigma$ between two simplicial fans
$\Delta$ and $\Delta'$ with the same support.  The simplicial cobordism
constructed in this section only satisfies the condition (1.4.1) above of
Morelli's idea.  The construction is done via a slick use of Sumihiro's
equivariant completion theorem [Sumihiro1,2]. 

\vskip.1in

In \S 3, we discuss the (weak) factorization between the lower face
$\partial_-\overline{\roman{Star}(\sigma)}$ and upper face
$\partial_+\overline{\roman{Star}(\sigma)}$, which we call the bistellar
operation, more in detail assuming the $\pi$-nonsingularity.

\vskip.1in

In \S 4, we achieve the condition (1.4.3), the collapsibility for the
simplicial cobordism
$\Sigma$.  By star subdividing $\Sigma$ further to obtain ${\tilde \Sigma}$, we
can make ${\tilde \Sigma}$ projective via the use of toric version of
Moishezon's theorem.  Projectivity implies collapsibility, achieving a
collapsible and simplicial cobordism ${\tilde \Sigma}$ between
$\partial_-{\tilde
\Sigma}$ and
$\partial_+{\tilde
\Sigma}$.  We can explicitly construct a collapsible and simplicial cobordism
$\Sigma_1$ (resp.
$\Sigma_2$) between
$\Delta$ and $\partial_-{\tilde
\Sigma}$ (resp. between $\partial_+{\tilde \Sigma}$ and $\Delta'$), as the
latter is obtained through star subdivisions (resp. star assemblings) from the
former.  Now we only have to take the composite $\Sigma_1 \circ {\tilde
\Sigma} \circ \Sigma_2$ to be the one providing a new collapsible and simplicial
cobordism between
$\Delta$ and $\Delta'$. 
 
\vskip.1in 

\S 5 is the most sutble and difficult part of the proof, achieving the
condition (1.4.2), i.e., the
$\pi$-nonsingularity of the cobordism $\Sigma$.  We introduce the invariant
``$\pi$-multiplicity profile" of a simplicial cobordism, which measures how far
$\Sigma$ is from being $\pi$-nonsingular, and observe that it strictly drops
after some appropriate star subdivisions.  By the descending chain condition on
the set of the $\pi$-multiplicity profiles, we acquire the $\pi$-nonsingularity
after finitely many star subdivisions.

\vskip.1in

Th arguments in \S 2 $\sim$ \S 5 put together provide the weak factorization,
solving the weak form of Conjecture 1.1 affirmatively.  The results are
summarized in \S 6.

\vskip.1in

We should emphasize that the weak form of Conjecture 1.1 is also
solved by [W{\l}odarczyk1] along a similar line of ideas but in a more
combinatorial language.

\vskip.1in

In \S 7, we finally show the strong factorization, based upon the weak
factorization achieved in the previous sections.  We obtain ${\tilde \Sigma}$ by
further star subdividing the cobordism $\Sigma$ corresponding to the weak
factorization between $\Delta$ and $\Delta'$, without affecting the lower face
of
$\Sigma$ but possibly smooth star subdividing the upper face of $\Sigma$, so
that the bistellar operations of the circuits in ${\tilde
\Sigma}$ only provide blowups starting from the lower face. 
We achieve the strong factorization
$$\Delta \cong \partial_-\Sigma = \partial_-{\tilde \Sigma} \leftarrow
\partial_+{\tilde \Sigma} \rightarrow \partial_+\Sigma \cong \Delta',$$
the first left arrow representing a sequence of smooth star subdivisions and the
second right arrow representing a sequence of smooth star assemblings
immediately after.

\vskip.1in

\S 7 discusses the generalization to the toroidal case.  All the arguments
above for the toric case can be lifted immediately to the toroidal case,
except for the existence of a cobordism and $\pi$-collapsibility, where we
used the global results like Sumihiro's and Moishezon's theorems only valid
in the toric case.  We circumvent these difficulties by a trick embedding a
toroidal conical complex into a usual toric fan after barycentric
star subdivisions.

\newpage

\S 2. Cobordism

\vskip.1in

We follow the usual notation and terminology concerning the toric varieties
$X_{\Delta}$ and their corresponding fans $\Delta$, as presented in
[Danilov,Fulton,Oda2].

We recall the notion of star subdivisions of a fan $\Delta$, the key operation
repeatedly used in this note.

\vskip.1in

$\bold{Definition\ 2.1.}$ Let $\tau \in \Delta$ be a
cone in a fan
$\Delta$.  Let $\rho$ be a ray passing in the relative interior of $\tau$. 
(Note that such $\tau \in \Delta$ containing $\rho$ in its relative interior
is uniquely determined once the ray $\rho$ is fixed.)  Then we define the star
subdivision $\rho \cdot \Delta$ of
$\Delta$ with respect to $\rho$ to be

$$\rho \cdot \Delta = (\Delta - \roman{Star}(\tau)) \cup \{\rho + \tau'
+ \nu;\tau' \text{\ a\ proper\ face\ of\ }\tau, \nu \in
\roman{link}_{\Delta}(\tau)\}$$
where
$$\align
\roman{Star}(\tau) &= \{\zeta \in \Delta;\zeta \supset \tau\}\\
\overline{\roman{Star}(\tau)} &= \{\zeta \in \Delta;\zeta \subset \eta \text{\
for\ some\ }\eta \in \roman{Star}(\tau)\}\\
\roman{link}_{\Delta}(\tau) &= \{\zeta \in \overline{\roman{Star}(\tau)};\zeta
\cap
\tau = \emptyset\}\\
\endalign$$

We call the inverse of a star subdivision a star assembling.

When $\tau = \langle\rho_1, \cdot\cdot\cdot, \rho_l\rangle$ is generated by
extremal rays $\rho_i$ with the primitive vectors $v_i = n(\rho_i)$ and the ray
$\rho$ is generated by the vector $v_1 + \cdot\cdot\cdot + v_l$, the star
subdivision is called the barycentric star subdivision with respect to
$\tau$.

When $\Delta$ is nonsingular, the barycentric star subdivision with respect to a
face $\tau$ is called a smooth star subdivision and its inverse a smooth star
assembling.

\vskip.1in

The notion of a cobordism as defined below sits in the center of Morelli's
idea.

\vskip.1in

$\bold{Definition\ 2.2.}$ Let $\Delta$ and $\Delta'$ be two fans in
$N_{\Bbb Q} = N \otimes {\Bbb Q}$ with the same support, where $N$ is the
lattice of one-parameter subgroups of the torus.  A cobordism $\Sigma$ is a fan
in
$N_{\Bbb Q}^+ = (N \oplus {\Bbb Z}) \otimes {\Bbb Q} = N_{\Bbb Q} \oplus {\Bbb
Q}$ equipped with the natural projection
$$\pi:N_{\Bbb Q}^+ = N_{\Bbb Q} \oplus {\Bbb Q} \rightarrow N_{\Bbb Q}$$
such that

\ \ (2.2.1) any cone $\tau \in \Sigma$ is $\pi$-strongly convex, i.e.,
$$x,y \in \tau, \pi(x) = - \pi(y) \Longrightarrow x = y = 0,$$

\ \ (2.2.2) the projection $\pi$ gives an isomorphism between $\partial_-\Sigma$
and
$\Delta$ (resp. $\partial_+\Sigma$ and $\Delta'$) as linear complexes, i.e.,
there is a one-to-one correspondence between the cones $\sigma_-$ of
$\partial_-\Sigma$ (resp.
$\sigma'_+$ of $\partial_+\Sigma$) and the cones
$\sigma$ of $\Delta$ (resp. $\sigma'$ of $\Delta'$) such that $\pi:\sigma_-
\rightarrow
\sigma$ (resp. $\pi:\sigma'_+
\rightarrow
\sigma'$) is a linear isomorphism for each $\sigma_-$ (resp. $\sigma'_+$) and
its corresponding
$\sigma$ (resp. $\sigma'$). (Note that we do NOT require the map of lattices
$\pi:(N \oplus {\Bbb Z}) \cap
\sigma_- \rightarrow N \cap \sigma$ (resp. $\pi:(N \oplus {\Bbb Z}) \cap
\sigma'_+ \rightarrow N \cap \sigma'$) to be an isomorphism.)  We denote this
isomorphism by
$$\pi:\partial_-\Sigma \overset{\sim}\to{\rightarrow} \Delta
\hskip.1in (\text{resp.\ }
\pi:\partial_+\Sigma \overset{\sim}\to{\rightarrow} \Delta')$$
where
$$\align
\partial_-\Sigma = &\{\tau \in \Sigma;(x, y - \epsilon) \not\in
\roman{Supp}(\Sigma) \\
&\text{\ for\ any\ }(x, y) \in \tau \text{\ with\ }x \in N_{\Bbb Q}, y \in
{\Bbb Q} \text{\ and\ any\ sufficiently\ small\ }\epsilon > 0\}\\
(\text{resp.\ } \partial_+\Sigma = &\{\tau \in \Sigma;(x, y + \epsilon)
\not\in
\roman{Supp}(\Sigma) \\
&\text{\ for\ any\ }(x, y) \in \tau \text{\ with\ }x \in N_{\Bbb Q}, y
\in {\Bbb Q} \text{\ and\ any\ sufficiently\ small\ }\epsilon > 0\})\\
\endalign$$

\ \ (2.2.3) the support $\roman{Supp}(\Sigma)$ of $\Sigma$ lies between the
lower face
$\partial_-\Sigma$ and the upper face $\partial_+\Sigma$, i.e.,
$$\align
\roman{Supp}(\Sigma) &= \{(x, y) \in N_{\Bbb Q}^+;x \in \roman{Supp}(\Delta)
= \roman{Supp}(\Delta') \text{\ and\ } y^x_- \leq y \leq y^x_+ \\
&\text{\ where\ }(x, y^x_-) \in
\roman{Supp}(\partial_-\Sigma) \text{\ and\ }(x, y^x_+) \in
\roman{Supp}(\partial_+\Sigma).\}\\
\endalign$$

\vskip.1in

We remark that actually we only need the condition (2.2.2) for the definition of
a cobordism and that the conditions (2.2.1) and (2.2.3) follow as the
consequences of (2.2.2).  We put all of these conditions as parts of the
definition above to clarify its basic properties.

\vskip.1in

\proclaim{Theorem 2.3} Let $\Delta$ and $\Delta'$ be two simplicial fans in
$N_{\Bbb Q} = N \otimes {\Bbb Q}$ with the same support.  Then there exists a
cobordism $\Sigma$ between $\Delta$ and $\Delta'$.  We may also require $\Sigma$
to be simplicial.
\endproclaim

\demo{Proof}\enddemo First we embed $\Delta$ ``at the
level $-1$" into $N_{\Bbb Q}^+$ so that the embedding $\Delta_-$ maps
isomorphically back onto
$\Delta$ by the projection $\pi$.  Namely, we take the fan $\Delta_-$ in
$N_{\Bbb Q}^+$ consisting of the cones
$\sigma_-$ of the form
$$\sigma_- = \langle(v_1,-1), \cdot\cdot\cdot, (v_k,-1)\rangle,$$
where the corresponding cone $\sigma = \langle\rho_1, \cdot\cdot\cdot,
\rho_k\rangle \in
\Delta$ is generated by the extremal rays $\rho_i$ with the primitive
vectors $v_i = n(\rho_i)$.  Similarly we embed $\Delta'$ ``at the
level $+1$" into $N_{\Bbb Q}^+$ so
that the embedding $\Delta'_+$ maps isomorphically back onto $\Delta'$ by the
projection
$\pi$.

We take $\Gamma$ to be the fan in $N_{\Bbb Q}^+$ consisting of the cones in
$\Delta_-$ and $\Delta'_+$ and the cones $\zeta$ of the form
$$\zeta = \langle(v,-1),(v,+1)\rangle,$$
where the $v$ vary among all the primitive vectors for the extremal rays
$\rho_v$ such that $\rho_v$ is a generator for some $\sigma \in \Delta$ and
some $\sigma' \in \Delta'$ simultaneously.

Now by Sumihiro's equivariant completion theorem [Sumihiro1,2], there exists a
fan $\Sigma^{\circ}$ with $\roman{Supp}(\Sigma^{\circ}) = N_{\Bbb Q}^+$ and
containing
$\Gamma$ as a subfan.

We only have to take $\Sigma$ to be
$$\Sigma = \{\tau \in \Sigma^{\circ}; \roman{Supp}(\tau) \subset S\}$$
where the set $S$ is described as
$$\align
S = \{(x,y) \in N_{\Bbb Q}^+;&x \in \roman{Supp}(\Delta) =
\roman{Supp}(\Delta'), y^x_- \leq y \leq y^x_+\\
&\text{\ with\ }(x,y^x_-) \in \Delta_-, (x,y^x_+) \in \Delta'_+\}.\\
\endalign$$

The cobordism
$\Sigma$ constructed as above may not be simplicial.  We take all
the cones in
$\Sigma$ which are not simplicial, and give them the partial order according to
the inclusion relation.  We take a succession of barycentric star subdivisions
with respect to these cones in the order compatible with the partial order,
starting with the maximal ones.  The resulting fan ${\tilde \Sigma}$ is
simplicial with the property
$$\align
\pi:\partial_-{\tilde \Sigma} = \partial_-\Sigma
&\overset{\sim}\to{\rightarrow}
\Delta\\
\pi:\partial_+{\tilde \Sigma} = \partial_+\Sigma
&\overset{\sim}\to{\rightarrow}
\Delta',\\
\endalign$$
providing a simplicial cobordism between $\Delta$ and $\Delta'$.  (We also
refer the reader to [Oda-Park,Park] for a more systematic treatment.)

\newpage

\S 3. Circuits and Bistellar Operations

\vskip.1in

In this section, we discuss how to read off the information on the
factorization from the circuits of a $\pi$-nonsingular cobordism.

\vskip.1in

$\bold{Definition\ 3.1.}$ Let $\Sigma$ be a simplicial fan in $(N
\oplus Z) \otimes {\Bbb Q} = N_{\Bbb Q}^+$ with the natural projection
$\pi:N_{\Bbb Q}^+ \rightarrow N_{\Bbb Q}$.  Assume that all the cones in
$\Sigma$ are $\pi$-strictly convex.  

A cone $\sigma \in \Sigma$ is $\pi$-independent if $\pi:\sigma \rightarrow
\pi(\sigma)$ is an isomorphism.  Otherwise $\sigma$ is $\pi$-dependent.

A cone $\sigma \in \Sigma$ is called a circuit if it is minimal
among the $\pi$-dependent cones, i.e., if $\sigma$ is $\pi$-dependent and any
proper face of
$\sigma$ is $\pi$-independent.

A cone $\sigma \in \Sigma$ is $\pi$-nonsingular if the projection
$\pi(\tau)$ of each $\pi$-independent face $\tau \subset \sigma$ is nonsingular
as a cone in $N_{\Bbb Q}$ with respect to the lattice $N$.  We say that the fan
$\Sigma$ is
$\pi$-nonsingular if all the cones in $\Sigma$ are
$\pi$-nonsingular.

\vskip.1in

The following theorem describes the transformation, which we call the
bistellar operation, from the lower face
$\partial_-\sigma$ to the upper face $\partial_+\sigma$ of a circuit $\sigma$
of a simplicial and $\pi$-nonsingular cobordism $\Sigma$.  (More generally the
theorem describes the transformation from the lower face
$\partial_-\overline{\roman{Star}(\sigma)}$ to the upper face
$\partial_+\overline{\roman{Star}(\sigma)}$ of the closed star of a circuit
$\sigma$.)  It turns out that the bistellar operation corresponds to a smooth
blowup immediately followed by a smooth blowdown.

\proclaim{Theorem 3.2} Let $\Sigma$ be a simplicial and $\pi$-nonsingular
cobordism in $N_{\Bbb Q}^+$.  Let
$\sigma =
\langle\rho_1, \cdot\cdot\cdot, \rho_k\rangle \in
\Sigma$ be a circuit generated by the extremal rays $\rho_i$.  From each extremal
ray $\rho_i$ we take the vector of the form $(v_i,w_i) \in N_{\Bbb Q}^+ =
N_{\Bbb Q}
\oplus {\Bbb Q}$ where $v_i = n(\pi(\rho_i))$ is the primitive vector of the
projection $\pi(\rho_i)$.

\vskip.1in

$(3.2.1)$ There is a unique linear relation among the $v_i$ (up to renumbering)
of the form
$$\Sigma r_{\alpha}v_{\alpha} = v_1 + \cdot\cdot\cdot + v_l - v_{l+1} -
\cdot\cdot\cdot - v_k = 0
\text{\ for\ some\ }0 \leq l \leq k$$
with
$$\Sigma r_{\alpha}w_{\alpha} = w_1 + \cdot\cdot\cdot + w_l - w_{l+1} -
\cdot\cdot\cdot - w_k > 0.$$

$(3.2.2)$ All the maximal faces $\gamma_i$ (resp. $\gamma_j$) of
$\partial_-\sigma$ (resp.
$\partial_+\sigma$) are of the form
$$\align
\gamma_i &= <\rho_1, \cdot\cdot\cdot, \overset{\vee}\to{\rho_i},
\cdot\cdot\cdot, \rho_l, \rho_{l+1}, \cdot\cdot\cdot, \rho_k> \hskip.1in 1 \leq
i \leq l\\
(\text{resp.\ } \gamma_j &= <\rho_1, \cdot\cdot\cdot, \rho_l, \rho_{l+1},
\cdot\cdot\cdot, \overset{\vee}\to{\rho_j}, \cdot\cdot\cdot, \rho_k> \hskip.1in 
l + 1 \leq j \leq k.\\
\endalign$$

$(3.2.3)$ Let $l_{\sigma}$ be the extremal ray in $N_{\Bbb Q}$ generated by the
vector
$$v_1 + \cdot\cdot\cdot + v_l = v_{l+1} + \cdot\cdot\cdot + v_k.$$
The smooth star subdivision of $\pi(\partial_-\sigma)$ with respect to
$l_{\sigma}$ coincides with the smooth star subdivision of
$\pi(\partial_+\sigma)$ with respect to
$l_{\sigma}$, whose maximal faces are of the form
$$<\pi(\gamma_{ij}),l_{\sigma}> = \langle\pi(\rho_1), \cdot\cdot\cdot,
\overset{\vee}\to{\pi(\rho_i)},
\cdot\cdot\cdot, \pi(\rho_l), \pi(\rho_{l+1}), \cdot\cdot\cdot,
\overset{\vee}\to{\pi(\rho_j)}, \cdot\cdot\cdot,
\pi(\rho_k),l_{\sigma}\rangle.$$

Thus the transformation from$\pi(\partial_-\sigma)$ to $\pi(\partial_+\sigma)$
is a smooth star subdivision followed immediately after by a smooth star
assembling.  We call the transformation a bistellar operation.

Similarly, the transformation from
$\pi(\partial_-\overline{\roman{Star}(\sigma)})$ to
$\pi(\partial_+\overline{\roman{Star}(\sigma)})$ is a smooth star subdivision
followed immediately after by a smooth star assembling. 

\endproclaim

\demo{Proof}\enddemo (3.2.1) Since $\sigma$ is a circuit, it is $\pi$-dependent
and minimal by definition.  Hence we have a linear relation
$$\Sigma r_iv_i = 0 \text{\ with\ } r_i \neq 0 \text{\ for\ all\ }i.$$
Since $\sigma$ is simplicial, the $\rho_i$ are linearly independent in $N_{\Bbb
Q}^+$ and hence $\Sigma r_iw_i \neq 0$.  We choose the signs of the $r_i$ so
that $\Sigma r_iw_i > 0$.  We only have to prove \linebreak
$|r_1| = |r_2| =
\cdot\cdot\cdot = |r_k|$.  Indeed, since $\sigma$ is $\pi$-nonsingular, we have
$$\align
1 &= |\det(\overset{\vee}\to{v_1},v_2, \cdot\cdot\cdot, v_k)| = |\det(v_1,
\cdot\cdot\cdot, \overset{\vee}\to{v_i}, \cdot\cdot\cdot, v_k)|\\ 
&=
|\det(\frac{1}{r_1}(\Sigma_{\alpha \neq 1}r_{\alpha}v_{\alpha}),v_2,
\cdot\cdot\cdot, \overset{\vee}\to{v_i}, \cdot\cdot\cdot, v_k)| =
|\frac{r_i}{r_1}\det(\overset{\vee}\to{v_1}, v_2,
\cdot\cdot\cdot, v_k)|,\\
\endalign$$
which implies $|r_1| = |r_i|$ for all $i$.

\vskip.1in
 
(3.2.2) Note that since $\sigma$ is a circuit, any maximal face $\gamma$ of
$\sigma$ belongs either to $\partial_-\sigma$ or to $\partial_+\sigma$,
exclusively.

Suppose $\gamma_i = \langle\rho_1, \cdot\cdot\cdot, \overset{\vee}\to{\rho_i},
\cdot\cdot\cdot, \rho_k\rangle \in \partial_-\sigma$.  Then since $\sigma$ is a
circuit, for any point
$$p = \Sigma_{\alpha \neq i} c_{\alpha}(v_{\alpha},w_{\alpha}) \in
\roman{RelInt}(\gamma_i)
\text{\ with\ }c_{\alpha} > 0,$$
we have
$$p + (0,\epsilon) \in \sigma \text{\ for\ sufficiently\ small\ }\epsilon > 0.$$
By setting $\epsilon = t_{\epsilon} \cdot \Sigma r_{\alpha}w_{\alpha} \text{\
for\ }t_{\epsilon} > 0$, we obtain
$$p + (0,\epsilon) = \Sigma_{\alpha \neq i} (c_{\alpha} +
t_{\epsilon}r_{\alpha})(v_{\alpha},w_{\alpha}) + t_{\epsilon}r_i(v_i,w_i) \in
\sigma,$$
which implies $c_{\alpha} + t_{\epsilon}r_{\alpha} > 0 \text{\ for\ }\alpha \neq
i$ and $t_{\epsilon}r_i > 0$.  Therefore, we have $r_i > 0$.  Similarly,
if $\gamma_j = \langle\rho_1, \cdot\cdot\cdot, \overset{\vee}\to{\rho_j},
\cdot\cdot\cdot, \rho_k\rangle \in \partial_+\sigma$, then we have $r_j <
0$.  This proves the assertion (3.2.2).

\vskip.1in

The assertion (3.2.3) follows immediately from (3.2.1) and (3.2.2).

\vskip.1in

The assertion about the transformation from
$\pi(\partial_-\overline{\roman{Star}(\sigma)})$ to
$\pi(\partial_+\overline{\roman{Star}(\sigma)})$ is an easy consequence of the
description of the transformation from
$\pi(\partial_-\sigma)$ to $\pi(\partial_+\sigma)$.

\vskip.1in

This completes the proof of Theorem 3.2.

\newpage

\S 4. Collapsibility

\vskip.1in

Let $\Sigma$ be a simplicial cobordism between simplicial fans $\Delta$ and
$\Delta'$.  Noting that
$$\Sigma = \cup_{\sigma} \overline{\roman{Star}(\sigma)} \cup
\partial_-\Sigma$$  where the union is taken over the circuits $\sigma$, we may
try to factorize the transformation from
$\Delta$ to
$\Delta'$ into smooth star subdivisions and smooth star assemblings by replacing
$\partial_-\overline{\roman{Star}(\sigma)}$ with
$\partial_+\overline{\roman{Star}(\sigma)}$, if $\Sigma$ is $\pi$-nonsingular. 
If we think of the cobordism built up out of ``bubbles"
$\overline{\roman{Star}(\sigma)}$, this process might be considered as a
succession of ``collapsing" these bubbles.  The following simple example shows
that this succession of collapsing, which should correspond to the
factorization into smooth star subdivisions and smooth star assemblings, is not
always possible, unless we can arrange the way we break these bubbles in a
certain order.  This possiblility for the certain nice arrangement is what we
call ``collapsibility" in this section.

\vskip.1in

$\bold{Example\ 4.1.}$ We take two sets of vectors in $N_{\Bbb Q} = {\Bbb Z}^2
\otimes {\Bbb Q}$
$$\align
&\{v_1 = (1,0), v_2 = (0,1), v_3 = (-1,0), v_4 = (0,-1)\}\\
&\{v'_1 = (1,1), v'_2 = (-1,1), v'_3 = (-1,-1), v'_4 = (1,-1)\}\\
\endalign$$
and fans $\Delta$ and $\Delta'$ whose maximal cones consist of
$$\align
\Delta &\ni \sigma_{12} = \langle v_1,v_2\rangle, \sigma_{23} =
\langle v_2,v_3\rangle,
\sigma_{34} = \langle v_3,v_4\rangle, \sigma_{41} = \langle v_4,v_1\rangle\\
\Delta' &\ni \sigma'_{12} = \langle v'_1,v'_2\rangle, \sigma'_{23} =
\langle v'_2,v'_3\rangle,
\sigma'_{34} = \langle v'_3,v'_4\rangle, \sigma'_{41} = \langle
v'_4,v'_1\rangle.\\
\endalign$$
If we take the simplicial fan $\Sigma$ (in $N_{\Bbb Q}^+$) whose
maximal cones consist of
$$\align
\sigma_{124'} &= \langle(v_1,0),(v_2,0),(v'_4,1)\rangle\\
\sigma_{231'} &= \langle(v_2,0),(v_3,0),(v'_1,1)\rangle\\
\sigma_{342'} &= \langle(v_3,0),(v_4,0),(v'_2,1)\rangle\\
\sigma_{413'} &= \langle(v_4,0),(v_1,0),(v'_3,1)\rangle\\
\sigma_{4'1'2} &= \langle(v'_4,1),(v'_1,1),(v_2,0)\rangle\\
\sigma_{1'2'3} &= \langle(v'_1,1),(v'_2,1),(v_3,0)\rangle\\
\sigma_{2'3'4} &= \langle(v'_2,1),(v'_3,1),(v_4,0)\rangle\\
\sigma_{3'4'1} &= \langle(v'_3,1),(v'_4,1),(v_1,0)\rangle,\\
\endalign$$
then $\Sigma$ is a simplicial $\pi$-nonsingular cobordism between $\Delta$ and
$\Delta'$.

Observe, however, that we cannot ``collapse" any one of the maximal cones
$\sigma_{ijk'}$ to replace $\partial_-\sigma_{ijk'}$ with
$\partial_+\sigma_{ijk'}$.  In fact, the circuit graph attached to $\Sigma$ as
defined below is a directed cycle consisting of eight vertices 
$$\sigma_{124'} \rightarrow \sigma_{4'1'2} \rightarrow \sigma_{231'} \rightarrow
\sigma_{1'2'3} \rightarrow \sigma_{342'} \rightarrow \sigma_{2'3'4}
\rightarrow \sigma_{413'} \rightarrow \sigma_{3'4'1} \rightarrow
\sigma_{124'}.$$

\vskip.1in

$\bold{Definition\ 4.2.}$ Let $\Sigma$ be a simplicial cobordism in $N_{\Bbb
Q}^+$.  We define a directed graph, which we call the circuit graph of $\Sigma$
as follows: The vertices of the circuit graph consist of the circuits $\sigma$
of
$\Sigma$.  We draw an edge from $\sigma$ to $\sigma'$ if there is a point $p
\in \partial_+\overline{\roman{Star}(\sigma)} \cap
\partial_-\overline{\roman{Star}(\sigma')}$ such that
$$p - (0,\epsilon) \in \overline{\roman{Star}(\sigma)}, p + (0,\epsilon) \in
\overline{\roman{Star}(\sigma')} \text{\ for\ sufficiently\ small\ }\epsilon >
0.$$  We say $\Sigma$ is collapsible if the
circuit graph contains no directed cycle.  When $\Sigma$ is collapsible, the
circuit graph determines a partial order among the circuits: $\sigma \leq
\sigma'$ if there is an edge $\sigma \rightarrow \sigma'$.

\vskip.1in

\proclaim{Theorem 4.3} Let $\Delta$ and $\Delta'$ be two simplicial fans in
$N_{\Bbb Q}$ with the same support.  Then there exists a simplicial and
collapsible cobordism $\Sigma$ in $N_{\Bbb Q}^+$ between $\Delta$ and
$\Delta'$.
\endproclaim

\demo{Proof}\enddemo The proof consists of several steps.  The main idea of
Morelli's is to reduce the collapsibility to the projectivity.

\vskip.1in

Step 1. Show that the projectivity induces the collapsibility.

\vskip.1in

\proclaim{Proposition 4.4} Let $\Sigma$ be a simplicial cobordism in
$N_{\Bbb Q}^+$ and assume that $\Sigma$ is a (part of a) projective fan.  Then
$\Sigma$ is collapsible.
\endproclaim

\demo{Proof}\enddemo Since $\Sigma$ is a part of a projective fan 
(i.e.,a part of a fan $\Sigma'$ whose corresponding toric variety $X_{\Sigma'}$
is projective), there exists a function $h:\roman{Supp}(\Sigma) \rightarrow
{\Bbb Q}$ which is piecewise linear with respect to the fan $\Sigma$ and which
is strictly convex, i.e., we have
$$\frac{1}{2}\{h(v) + h(u)\} \geq h(\frac{1}{2}\{v + u\})$$
whenever the line segment $\overline{vu}$ is in $\roman{Supp}(\Sigma)$ and the
strict inequality holds whenever $v$ and $u$ are in two distinct maximal
cones (cf.[Fulton,Oda2]).

Let $\sigma$ and $\sigma'$ be two circuits with a directed edge, i.e., there
exists a point $p \in \partial_+\overline{\roman{Star}(\sigma)} \cap
\partial_-\overline{\roman{Star}(\sigma')}$ such that
$$p - (0,\epsilon) \in \overline{\roman{Star}(\sigma)}, p + (0,\epsilon) \in
\overline{\roman{Star}(\sigma')} \text{\ for\ sufficiently\ small\ }\epsilon >
0.$$ Take a maximal $\pi$-dependent cone $p \in \eta \supset \sigma$ (resp.
$p \in \eta' \supset \sigma'$) of
$\overline{\roman{Star}(\sigma)}$ (resp. $\overline{\roman{Star}(\sigma')}$)
such that
$p - (0,\epsilon) \in \eta$ (resp. $p + (0,\epsilon) \in \eta'$).

Take also linear functions $h_{\eta}, h_{\eta'}, h_{\sigma}, h_{\sigma'}$
which coincide with $h|_{\eta}, h|_{\eta'}, h|_{\sigma}, h|_{\sigma'}$,
respectively.

Then by the strict convexity of the function $h$, setting the coordinates of
\linebreak $p = (x,y)$ we have $\frac{1}{2}\{h(x,y + \epsilon) + h(x,y -
\epsilon)\} > h(x,y)$ or equivalently \linebreak $h_{\eta'}(0,1) >
h_{\eta}(0,1)$, and hence $h_{\sigma'}(0,1) > h_{\sigma}(0,1)$. (Note that
$(0,1) \in
\roman{span}_{\Bbb Q}(\sigma)$ for any $\pi$-dependent cone
$\sigma$.)

If $\sigma_1, \cdot\cdot\cdot, \sigma_l$ are circuits determining a directed
path in the circuit graph of $\Sigma$, then the above observation
shows $h_{\sigma_1}(0,1) < \cdot\cdot\cdot < h_{\sigma_l}(0,1)$. Thus the path
cannot be a cycle.  Therefore, $\Sigma$ is collapsible.

\newpage

Step 2. Show the toric version of Moishezon's theorem.

\vskip.1in

\proclaim{Theorem 4.5} Let $\Sigma$ be a fan in $N_{\Bbb Q}^+$.  Then there
exists a fan ${\tilde \Sigma}$ obtained from $\Sigma$ by a sequence of star
subdivisions such that ${\tilde \Sigma}$ is a (part of a) projective fan.
\endproclaim

\demo{Proof}\enddemo We may assume that $\roman{Supp}(\Sigma) = N_{\Bbb Q}^+$
and that $\Sigma$ is simplicial and nonsingular by applying some appropriate
sequence of star subdivisions to the original $\Sigma$.

By the toric version of Chow's Lemma (see, e.g., [Oda2], \S 2.3), we have a
projective fan $\Sigma'$ mapping to $\Sigma$, i.e., we have a projective toric
variety $X_{\Sigma'}$ with an equivariant proper birational morphism onto
$X_{\Sigma}$
$$g:X_{\Sigma'} \rightarrow X_{\Sigma}.$$
By the toric version of Hironaka's elimination of points of indeterminacy
(cf.[DeConcini-Procesi].) we can take a fan ${\tilde \Sigma}$ obtained from
$\Sigma$ by a sequence of smooth star subdivisions such that there exists an
equivariant proper birational morphism
$$f:X_{\tilde \Sigma} \rightarrow X_{\Sigma'}.$$
Since $g \circ f$ is projective as it is a sequence of smooth blowups and
since $g$ is separated, $f$ is also projective.  Now since $\Sigma'$ is a
projective fan, so is ${\tilde \Sigma}$.

\vskip.1in

Step 3. Composition of (collapsible) cobordisms.

\vskip.1in

Starting from a simplicial cobordism $\Sigma$ between $\Delta$ and
$\Delta'$ constructed as in Theorem 2.3 and then applying Step 2, we obtain
a simplicial cobordism ${\tilde \Sigma}$ between
$\partial_-{\tilde
\Sigma}$ and $\partial_+{\tilde \Sigma}$, where ${\tilde \Sigma}$ is a (part of
a) projective fan and hence collapsible and where $\partial_-{\tilde \Sigma}$
(resp. $\partial_+{\tilde \Sigma}$) is obtained from $\Delta$ (resp. $\Delta'$)
by a sequence of star subdivisions.  (Or equivalently $\Delta$ (resp.
$\Delta'$) is obtained from $\partial_-{\tilde \Sigma}$
(resp. $\partial_+{\tilde \Sigma}$) by a sequence of star assemblings.)  We
only have to construct a collapsible and simplicial cobordism $\Sigma_{\Delta}$
between
$\Delta$ and
$\partial_-{\tilde \Sigma}$ and another $\Sigma_{\Delta'}$ between
$\partial_+{\tilde \Sigma}$ and $\Delta'$ so that we compose them together
$\Sigma_{\Delta} \circ {\tilde \Sigma} \circ \Sigma_{\Delta'}$ to obtain a
collapsible and simplicial cobordism between $\Delta$ and $\Delta'$.

\proclaim{Proposition-Definition 4.6} Let $\Sigma_1$ and $\Sigma_2$
be cobordisms in $N_{\Bbb Q}^+$ such that

$(4.6.1)$ $\Sigma_1 \cup \Sigma_2$ is again a fan in $N_{\Bbb Q}^+$,

$(4.6.2)$ $\Sigma_1 \cap \Sigma_2 = \partial_+\Sigma_1 \cap \partial_-\Sigma_2$,

$(4.6.3)$ for any cone $\sigma \in \partial_+\Sigma_2$
$$\pi(\sigma) \not\subset \partial\{\pi(\partial_+\Sigma_1 \cup
\partial_+\Sigma_2)\} \text{\ and\ } \pi(\sigma) \subset
\partial(\pi(\partial_+\Sigma_2)) \Longrightarrow
\sigma
\in
\partial_+\Sigma_1,$$ 
and for any cone $\sigma \in \partial_+\Sigma_1$
$$\pi(\sigma) \not\subset \partial\{\pi(\partial_-\Sigma_1 \cup
\partial_-\Sigma_2)\} \text{\ and\ } \pi(\sigma) \subset
\partial(\pi(\partial_-\Sigma_1)) \Longrightarrow \sigma \in
\partial_-\Sigma_2.$$

Then the union $\Sigma_1 \cup \Sigma_2$, which we call the composite of
$\Sigma_1$ with $\Sigma_2$ and denote by $\Sigma_1 \circ \Sigma_2$, is a
cobordism.

Moreover, if both $\Sigma_1$ and $\Sigma_2$ are simplicial and collapsible,
then so is the composite $\Sigma_1 \circ \Sigma_2$.

\endproclaim

\demo{Proof}\enddemo By the condition (4.6.1) the composite $\Sigma_1
\circ \Sigma_2$ is a fan.  The conditions (4.6.2) and (4.6.3) guarantee
$\pi:\partial_-(\Sigma_1 \circ \Sigma_2) \rightarrow N_{\Bbb Q}$ and
$\pi:\partial_+(\Sigma_1 \circ \Sigma_2) \rightarrow N_{\Bbb Q}$ are
isomorphisms of linear complexes onto their images.  Thus $\Sigma_1 \circ
\Sigma_2$ is a cobordism.  The ``Moreover" part of the assertion is also clear.

We note that in case $\partial_+\Sigma_1 = \partial_-\Sigma_2$ the condition
(4.6.3) is automatically satisfied.

\vskip.1in

\proclaim{Proposition 4.7} Let ${\tilde \Delta}$ be a simplicial fan in
$N_{\Bbb Q}$ obtained from another simplicial fan $\Delta$ in $N_{\Bbb Q}$ by
a sequence of star subdivisions and star assemblings.  Suppose $\Delta$ is
embedded in $N_{\Bbb Q}$
$$s:\Delta \hookrightarrow N_{\Bbb Q}^+$$
so that $\pi \circ s$ is the identity of the fan.

Then there exists a simplicial and collapsible cobordism $\Sigma$ between
$\Delta$ and ${\tilde \Delta}$ (resp. between ${\tilde \Delta}$ and
$\Delta$) such that $\partial_-\Sigma = s(\Delta)$ (resp. $\partial_+\Sigma =
s(\Delta)$).
\endproclaim

\demo{Proof}\enddemo We only have to prove the assertion when the sequence
consists of a single star subdivision or a star assembling.

Suppose ${\tilde \Delta}$ is obtained from $\Delta$ by a star subdivision with
respect to a ray $\rho$ passing through the relative interior of a face $\tau
\in \Delta$.  Say that the ray $\rho$ is generated by a primitive vector
$v_{\rho}$.  Then we only have to take by fixing some sufficiently large
$y_{\rho} > 0$
$$\Sigma = s(\Delta) \cup \{\langle s(\zeta),(v_{\rho},y_{\rho})\rangle;\zeta
\in
\Delta,
\zeta \subset \sigma \text{\ for\ some\ }\sigma \in \Delta \text{\ with\
}\sigma \ni
\rho\}.$$

Suppose ${\tilde \Delta}$ is obtained from $\Delta$ by a star assembling,
which is the inverse of a star subdivision with respect to a ray $\rho$
passing through the relative interior of a face $\tau \in {\tilde \Delta}$.
Let $\tau = \langle <\rho_1, \cdot\cdot\cdot, \rho_k\rangle$ be generated by
extremal rays $\rho_i$ with the primitive vectors $v_{\rho_i} = n(\rho_i)$.  We
construct $\Sigma_1, \cdot\cdot\cdot, \Sigma_k$ with $s_i:\Delta
\overset{\sim}\to{\rightarrow} \partial_+\Sigma_i$ and
$\Sigma$ as required inductively. 

Fixing some sufficiently large
$y_{\rho_1} > 0$, we take
$$\Sigma_1 = s(\Delta) \cup \{\langle
s(\zeta),(v_{\rho_1},y_{\rho_1})\rangle;\zeta
\in
\Delta,
\zeta \subset \sigma \text{\ for\ some\ }\sigma \in \Delta \text{\ with\
}\sigma \ni
\rho_1\}.$$
Obviously $\partial_+\Sigma_1$ is isomorphic to $\Delta$ via the projection
$\pi$, and we set the inverse $s_1:\Delta \overset{\sim}\to{\rightarrow}
\partial_+\Sigma_1$.

Suppose we have already constructed $\Sigma_1, \cdot\cdot\cdot, \Sigma_{i-1}$
with a sequence of positive numbers $0 < y_{\rho_1} < \cdot\cdot\cdot <
y_{\rho_{i-1}}$ where each positive number is sufficiently larger than the
previous one, and with the isomorphisms $s_1, \cdot\cdot\cdot, s_{i-1}$ from
$\Delta$ to $\partial_+\Sigma_1, \cdot\cdot\cdot, \partial_+\Sigma_{i-1}$,
respectively.  Then by fixing some positive number $y_{\rho_i}$ which is
sufficiently larger than $y_{\rho_{i-1}}$, we take
$$\Sigma_i = \Sigma_{i-1} \cup \{\langle
s_{i-1}(\zeta),(v_{\rho_i},y_{\rho_i})\rangle;\zeta
\in
\Delta,
\zeta \subset \sigma \text{\ for\ some\ }\sigma \in \Delta \text{\ with\
}\sigma \ni
\rho_i\}.$$
Again clearly $\partial_+\Sigma_i$ is isomorphic to $\Delta$ via the projection
$\pi$, and we set the inverse $s_i:\Delta \overset{\sim}\to{\rightarrow}
\partial_+\Sigma_i$.

Thus we have constructed $\Sigma_1, \cdot\cdot\cdot, \Sigma_k$.  We only have
to set
$${\tilde \Sigma} = \Sigma_k \cup \langle s_k(\rho_1), \cdot\cdot\cdot,
s_k(\rho_k)\rangle \cup \langle s_k(\rho_1), \cdot\cdot\cdot,
s_k(\rho_k), s_k(\rho)\rangle.$$

This completes the proof of Proposition 4.7.

\vskip.1in

Thus we complete Step 3 and hence the proof of Theorem 4.5.

\vskip.1in

In \S 5, starting from a collapsible and simplicial cobordism between two
nonsingular fans $\Delta$ and $\Delta'$ (which we constructed in this
section), we try to construct another cobordism which is not only collapsible
and simplicial but also $\pi$-nonsingular, by further star subdividing the
original cobordism.  It is worthwhile to note that the collapsibility is
preserved under star subdivisions.

\proclaim{Lemma 4.8} Let $\Sigma$ be a simplicial cobordism in
$N_{\Bbb Q}^+$, which is collapsible.  Then any simplicial cobordism ${\tilde
\Sigma}$ obtained from $\Sigma$ by a star subdivision, with respect to a ray
$\rho$, is again collapsible.
\endproclaim

\demo{Proof}\enddemo Note first that if $\Sigma$ consists of the closed star
of a single circuit, then $\rho \cdot \Sigma = \rho \cdot
\overline{\roman{Star}(\sigma)}$ is easily seen to be collapsible.

In general, number the circuits $\sigma_1, \sigma_2, \cdot\cdot\cdot,
\sigma_m$ of $\Sigma$ so that $\sigma_i$ is minimal among $\sigma_i,
\sigma_{i+1}, \cdot\cdot\cdot,
\sigma_m$ according to the order given by the circuit graph.  Then setting
$$\Sigma = \cup_{i = 1}^m \overline{\roman{Star}(\sigma_i)} \cup
\partial_+\Sigma,$$ we have
$$\rho \cdot \Sigma = \cup_{i = 1}^m \rho \cdot
\overline{\roman{Star}(\sigma_i)}
\cup \rho \cdot \partial_+\Sigma$$
and 
$$\rho \cdot \Sigma = \{\rho \cdot
\overline{\roman{Star}(\sigma_1)}\} \circ \cdot\cdot\cdot \circ \{\rho \cdot
\overline{\roman{Star}(\sigma_m)}\} \circ \{\rho \cdot \partial_+\Sigma\}$$
is collapsible by the first observation and by Proposition-Definition 4.6. 

\newpage

\S 5. $\pi$-Desingularization.

\vskip.1in

The purpose of this section, which is technically the most subtle, is to show
the follwoing theorem of ``$\pi$-desingularization".

\proclaim{Theorem 5.1} Let $\Sigma$ be a simplicial cobordism in $N_{\Bbb
Q}^+$.  Then there exists a simplicial cobordism ${\tilde \Sigma}$ obtained from
$\Sigma$ by a sequence of star subdivisions such that ${\tilde \Sigma}$ is
$\pi$-nonsingular.  Moreover, the sequence can be taken so that any
$\pi$-independent and already $\pi$-nonsingular face of $\Sigma$ remains
unaffected during the process. 
\endproclaim

\vskip.1in

Naively, just LIKE the case of the usual desingularization of toric fans, we
would like to subdivide any $\pi$-independent face with
$\pi$-multiplicity bigger than 1 so that its $\pi$-multiplicitiy drops. 
However, UNLIKE the case of the usual desingularization, we might introduce
a new $\pi$-independent face of uncontrollably high $\pi$-multiplicity if
we subdivide blindly, though we may succeed in decreasing the
$\pi$-multiplicity of the $\pi$-independent face that we picked originally. 
This is where the difficulty lies !  We outline Morelli's ingeneous strategy to
subdivide carefully to avoid introducing new
$\pi$-independent faces with high $\pi$-multiplicity and
achieve $\pi$-desingularization.  It consists of the following four steps:

\vskip.1in

Step 1: Introduce the invariant ``$\pi$-multiplicity profile"
$\text{$\pi$-m.p.}(\Sigma)$ of a simplicial cobordism $\Sigma$, which measures
how far $\Sigma$ is from being $\pi$-nonsingular.

\vskip.1in

Step 2: Observe that the star subdivision $\eta' =
\roman{Mid}(\tau,l_q) \cdot \eta$ of a simplex $\eta$ by an interior point of a
face
$\tau$ does not increase the $\pi$-multiplicity profile, i.e.,

$$\text{$\pi$-m.p.}(\eta') \leq \text{$\pi$-m.p.}(\eta)$$

if

\ \ (i) $\tau$ is ``codefinite" with respect to $\eta$, and

\ \ (ii) the interior point corresponds to the midray $\roman{Mid}(\tau,l_q)$,
where the ray $l_q$ is generated by a lattice point $q \in
\roman{par}(\pi(\tau))$.

Moreover, if $\tau$ is contained in a maximal $\pi$-independent face
$\gamma$ of $\eta$ with the maximum $\pi$-multiplicity $h_{\eta}$, i.e., if
$$\tau \subset \gamma \text{\ and\ }\text{$\pi$-mult}(\gamma) =
h_{\eta} = \max\{\text{$\pi$-mult}(\zeta);\zeta \subset \eta\},$$
then the $\pi$-multiplicity profile strictly drops

$$\text{$\pi$-m.p.}(\eta') < \text{$\pi$-m.p.}(\eta).$$

Step 3: Let $\tau$ be a $\pi$-independent face in the closed
star $\overline{\roman{Star}(\sigma)}$ of a circuit $\sigma$ in $\Sigma$. 
Introduce the notion of the star subdivision by the negative or positive center
point of $\sigma$.  We can find
$\Sigma^{\circ}$ such that

\ \ (i) $\Sigma^{\circ}$ is obtained by a succession of appropriate
star subdivisions by negative or positive center points of circuits inside of
$\sigma$,

\ \ (ii) the $\pi$-multiplicity profile does not increase, i.e.,
$$\text{$\pi$-m.p.}(\Sigma^{\circ}) \leq
\text{$\pi$-m.p.}(\Sigma),$$

\ \ (iii) $\tau$ is a face of $\Sigma^{\circ}$ such that $\tau$ is codefinite
with respect to every cone $\eta \in \Sigma^{\circ}$ containing $\tau$.

\vskip.1in

Step 4: Combine Step 2 and Step 3 to find ${\tilde \Sigma}$
obtained from $\Sigma$ by a succession of star subdivisions such that the
$\pi$-multiplicity profile strictly drops
$$\text{$\pi$-m.p.}({\tilde \Sigma}) < \text{$\pi$-m.p.}(\Sigma).$$
As the set of the $\pi$-multiplicity profiles satisfies the descending chain
condition, we reach a $\pi$-nonsingular cobordism after finitely many
star subdivisions as required.  

In fact, by Step 3 we can find a $\pi$-independent face $\tau$
of a maximal cone $\eta' \subset \Sigma^{\circ}$ such that

\ \ (i) $\text{$\pi$-m.p.}(\eta')$ is maximum among the $\pi$-multiplicity
profiles of all the maximal cones of $\Sigma^{\circ}$,

\ \ (ii) $\tau$ is contained in a maximal $\pi$-independent face $\gamma$ of
$\eta'$ with the maximum $\pi$-multiplicity $\text{$\pi$-mult}(\gamma) =
h_{\eta'}$,

\ \ (iii) $\tau$ is codefinite with respect to $\eta'$ and with respect to all
the other maximal cones containing $\tau$,

\ \ (iv) we can find a lattice point $q \in \roman{par}(\pi(\tau))$.

We only have to set ${\tilde \Sigma} = Mid(\tau,l_q) \cdot \Sigma^{\circ}$ to
observe by Step 2 that $\text{$\pi$-m.p.}({\tilde \Sigma}) <
\text{$\pi$-m.p.}(\Sigma)$.

This completes the process of $\pi$-desingularization.

\vskip.1in

Now we discuss the details of each step.

\vskip.1in

$\boxed{\roman{Step}\ 1}$

\vskip.1in

$\bold{Definition\ 5.2.}$ Let $\gamma$ be a simplicial cone in $N_{\Bbb Q}^+$. 
If $\gamma$ is $\pi$-independent, then we define the $\pi$-multiplicity of
$\gamma$ to be
$$\text{$\pi$-mult}(\gamma) = |\det(v_1, \cdot\cdot\cdot, v_k)|,$$
where the $v_i = n(\pi(\rho_i))$ are the primitive vectors of the projections
of the extremal rays $\rho_i$ generating $\gamma = \langle\rho_1,
\cdot\cdot\cdot,
\rho_k\rangle$.  If $\gamma$ is $\pi$-dependent, then we
set $\text{$\pi$-mult}(\gamma) = 0$ by definition.

\vskip.1in

Let $\eta$ be a simplicial and $\pi$-strictly convex cone in
$N_{\Bbb Q}^+$ with
$$\align
h_{\eta} &= \max\{\text{$\pi$-mult}(\gamma);\gamma \text{\ is\ a\
$\pi$-independent\ face\ of\ }\eta\},\\ 
k_{\eta} &= \dim \sigma \text{\ where\ }\sigma
\text{\ is\ the\ unique\ circuit\ contained\ in\ }\eta,\\ 
r_{\eta} &= \text{the\ number\ of\ the\ maximal\ $\pi$-independent\ faces\ of\
$\eta$}\\
&\text{having\ the\ maximum\ $\pi$-multiplicity\ }h_{\eta}.\\
\endalign$$ 
We define the $\pi$-multiplicity profile
$\text{$\pi$-m.p.}(\eta)$ of $\eta$ to be the orderd
quadruple of numbers
$$\text{$\pi$-m.p.}(\eta) = (a_{\eta},b_{\eta},c_{\eta},d_{\eta})$$
where
$$\align
a_{\eta} &= h_{\eta} \\
b_{\eta} &= \left\{\aligned
0 & \text{\ if\ }r_{\eta} \leq 1 \\
1 & \text{\ if\ }r_{\eta} > 1, \\
\endaligned\right.\\
c_{\eta} &= \left\{\aligned
0 & \text{\ if\ }b_{\eta} = 0 \\
k_{\eta} & \text{\ if\ }b_{\eta} =1, \\
\endaligned\right.\\
d_{\eta} &= \left\{\aligned
0 & \text{\ if\ }c_{\eta} = 0 \\
r_{\eta} & \text{\ if\ }c_{\eta} > 0.
\endaligned\right.\\
\endalign$$

We order the set of the $\pi$-multiplicity profiles of all the simplicial
and $\pi$-strictly convex cones in
$N_{\Bbb Q}^+$ lexicographically. 

We define the $\pi$-multiplicity profile $\text{$\pi$-m.p.}(\Sigma)$ of a
simplicial cobordism $\Sigma$ in $N_{\Bbb Q}^+$ to be
$$\text{$\pi$-m.p.}(\Sigma) = [g_{\Sigma};s]$$
where 
$$g_{\Sigma} = \max\{\text{$\pi$-m.p.}(\eta);\eta \text{\ is\ a\ maximal\
simplicial\ cone\ of\ }\Sigma\}$$ 
and where $s$ is the number of the maximal simplicial cones of $\Sigma$ having
the maximum $\pi$-multiplicity profile $g_{\Sigma}$.

When a simplicial cobordism $\Sigma$ consists of only one maximal simplicial
and $\pi$-strictly convex cone
$\eta$ (and its faces), we understand as a convention
$$\text{$\pi$-m.p.}(\Sigma) = [\text{$\pi$-m.p.}(\eta);1] =
\text{$\pi$-m.p.}(\eta).$$ 

\vskip.1in

The definition of the invariant $\pi$-multiplicity profile may look
heuristic at this point.  At the end of the section, we discuss how Morelli
reached this definition after a couple of false trials in [Morelli1,2].  The
behavior of the $\pi$-multiplicity profile under several kinds of
star subdivisions will be the key in Step 3. 

\vskip.2in

$\boxed{\roman{Step}\ 2}$

\vskip.1in

$\bold{Definition\ 5.3.}$ Let $\eta$ be a simplicial, $\pi$-dependent and
$\pi$-strictly convex cone in $N_{\Bbb Q}^+$.  A $\pi$-independent face $\tau$
of $\eta$ is said to be codefinite with respect to $\eta$ if the
set of generators of $\tau$ does not contain both positive and negative extremal
rays
$\rho_i$ of 
$\eta$.  That is to say, if $\Sigma r_iv_i = 0$ is the nontrivial linear
relation for
$\eta$ among the primitive vectors $v_i = n(\pi(\rho_i))$, then the generators
for $\tau$ contain only those extremal rays in the set
$\{\rho_i;r_i < 0\}$ or in the set $\{\rho_i;r_i > 0\}$, exclusively.

\vskip.1in

$\bold{Notation\ 5.4.}$ Let $\tau$ be a cone in a simplicial
cobordism $\Sigma$ in $N_{\Bbb Q}^+$ and
$l$ a ray in $\pi(\tau)$.  Then we define
the ``midray" $\roman{Mid}(\tau,l)$ to be the ray generated by the middle point
of the line segment
$\tau \cap
\pi^{-1}(n(l))$.  (If $\tau \cap
\pi^{-1}(n(l))$ consists of a point, then $\roman{Mid}(\tau,l)$ is the ray
generated by that point.)  

Let $\gamma = \langle\rho_1, \cdot\cdot\cdot, \rho_k\rangle$ be a
$\pi$-independent cone in
$N_{\Bbb Q}^+$ generated by the extremal rays $\rho_i$ with the corresponding
primitive generators
$v_i = n(\pi(\rho_i)) \in N$.  Then
we define the set
$$\roman{par}(\pi(\gamma)) = \{m \in N; m = \Sigma_i a_iv_i, 0 < a_i < 1\}.$$

\proclaim{Proposition 5.5} Let $\tau$ be a $\pi$-independent face of
a simplicial, $\pi$-dependent and $\pi$-strictly convex cone $\eta$ in
$N_{\Bbb Q}^+$.  Assume
$\tau$ is codefinite
with respect to $\eta$.  Let $\eta' = \roman{Mid}(\tau,l_q) \cdot \eta$ be the
star subdivision of $\eta$ by the midray $\roman{Mid}(\tau, l_q)$ where the ray
$l_q$ is generated by a lattice point $q \in \roman{par}(\pi(\tau))$.  Then the
$\pi$-multiplicity profile does not increase under the star subdivision, i.e.,
$$\text{$\pi$-\rm{m.p.}}(\eta') \leq \text{$\pi$-\rm{m.p.}}(\eta).$$
Moreover, if $\tau$ is contained in a maximal codimension one face $\gamma$ of
$\eta$ with
$$\text{$\pi$-\rm{mult}}(\gamma) = h_{\eta} =
\max\{\text{$\pi$-\rm{mult}}(\zeta);\zeta \subset \eta\},$$ 
then the $\pi$-multiplicity strictly decreases, i.e.,
$$\text{$\pi$-\rm{m.p.}}(\eta') < \text{$\pi$-\rm{m.p.}}(\eta).$$
\endproclaim

\demo{Proof}\enddemo We claim first that all the NEW maximal $\pi$-independent
faces $\gamma'$ of $\eta'$ have $\pi$-multiplicities strictly smaller than
$h_{\eta}$, i.e.,
$$\text{$\pi$-mult}(\gamma') < h_{\eta}.$$

Let $\tau = \langle \rho_1, \cdot\cdot\cdot, \rho_n\rangle$ be generated by
the extremal rays $\rho_i$ with the corresponding primitive vectors $v_i =
n(\pi(\rho_i)), i = 1, \cdot\cdot\cdot, n$.  We can write $0 \neq q = \Sigma_i
a_iv_i$ with
$0 < a_i < 1$ for all $i$, as $q \in \roman{par}(\pi(\tau))$.  

Any new maximal $\pi$-independent face $\gamma'$ in $\eta'$ has the form
$$\gamma' = \rho' + \tau' + \nu$$
where $\rho' = \roman{Mid}(\tau,l_q)$, $\tau'$ is a proper face of $\tau$
with $\rho'
\not\in
\tau'$ and where $\nu \in \roman{link}_{\eta}(\tau)$.

Observe that in general a maximal $\pi$-independent face of a simplicial cone
in $N_{\Bbb Q}^+$ has codimension at most one and hence we may assume that in
the above expression $\tau'$ has codimension at most two in $\tau$.

\vskip.1in

Case: $\tau'$ has codimension one in $\tau$.

\vskip.1in

The face $\tau'$ omits, say, $\rho_j$ among the extremal rays of $\tau$.  Then

$$\text{$\pi$-mult}(\rho' + \tau' + \nu') \leq a_j \cdot \text{$\pi$-mult}(\tau
+
\nu) \leq a_j \cdot h_{\eta} < h_{\eta}.$$

\vskip.1in

Case: $\tau'$ has codimension two in $\tau$.

\vskip.1in

The face $\tau'$ omits, say, $\rho_j$ and $\rho_k$ among the extremal rays of
$\tau$.  Observe that in this case $\tau + \nu$ is necessarily
$\pi$-dependent.  Indeed, if $\tau + \nu$ is $\pi$-independent, then there
exists a codimension one face $\tau'' (\supset \tau')$ of $\tau$ such that
we have
$\pi$-independent faces
$$\tau + \nu \supset \rho' + \tau'' + \nu \underset{\neq}\to{\supset} \rho' +
\tau' + \nu,$$ contradicting the maximality of $\rho' + \tau' + \nu$.  

Let $\nu = \langle\rho_{n+1}, \cdot\cdot\cdot, \rho_m\rangle$
be generated by the extremal rays $\rho_i$ with the corresponding primitive
vectors $v_i = n(\pi(\rho_i)), i = n+1, \cdot\cdot\cdot, m$ as before.  Then,
since $\tau + \nu$ is $\pi$-dependent, we have a nontrivial linear dependence
relation $\Sigma_{i = 1}^m r_iv_i = 0$.  In order to compute the
$\pi$-multiplicity, choose a basis of $\{\roman{span}_{\Bbb Q}\pi(\tau + \nu)\}
\cap N$.  Then
$$\align
(\diamond) \hskip.1in &\text{$\pi$-mult}(\rho' + \tau' + \nu) \\
&\leq |\det(q, v_1,
\cdot\cdot\cdot,
\overset{\vee}\to{v_j}, \cdot\cdot\cdot, \overset{\vee}\to{v_k},
\cdot\cdot\cdot, v_m)| \\
&= |\Sigma_i a_i \cdot \det(v_i, v_1, \cdot\cdot\cdot, \overset{\vee}\to{v_j},
\cdot\cdot\cdot, \overset{\vee}\to{v_k}, \cdot\cdot\cdot, v_m)| \\
&= |a_j \cdot \det(v_j, v_1, \cdot\cdot\cdot, \overset{\vee}\to{v_j},
\cdot\cdot\cdot, \overset{\vee}\to{v_k}, \cdot\cdot\cdot, v_m) + a_k \cdot
\det(v_k, v_1, \cdot\cdot\cdot, \overset{\vee}\to{v_j},
\cdot\cdot\cdot, \overset{\vee}\to{v_k}, \cdot\cdot\cdot, v_m)|.\\
\endalign$$ 
On the other hand, we have
$$\align
0 &= |\Sigma_i r_i \cdot \det(v_i, v_1, \cdot\cdot\cdot, \overset{\vee}\to{v_j},
\cdot\cdot\cdot, \overset{\vee}\to{v_k}, \cdot\cdot\cdot, v_m)| \\
&= |r_j \cdot \det(v_j, v_1, \cdot\cdot\cdot, \overset{\vee}\to{v_j},
\cdot\cdot\cdot, \overset{\vee}\to{v_k}, \cdot\cdot\cdot, v_m) + r_k \cdot
\det(v_k, v_1, \cdot\cdot\cdot, \overset{\vee}\to{v_j},
\cdot\cdot\cdot, \overset{\vee}\to{v_k}, \cdot\cdot\cdot, v_m)|.\\
\endalign$$ 
Since $\tau$ is codefinite with respect to $\eta$, either $r_j$ and $r_k$
have the same sign or one of them is 0.  (If $r_j = r_k = 0$, then $\tau' +
\nu$ would be $\pi$-dependent since $\Sigma_{i \neq j, k}r_iv_i =
\Sigma_ir_iv_i = 0$.  But $\rho' + \tau' + \nu$, containing $\tau' + \nu$,
is
$\pi$-independent, a contradiction!)  In the former case, $\det(v_j, v_1,
\cdot\cdot\cdot, \overset{\vee}\to{v_j},
\cdot\cdot\cdot, \overset{\vee}\to{v_k}, \cdot\cdot\cdot, v_m)$ and $\det(v_k,
v_1, \cdot\cdot\cdot, \overset{\vee}\to{v_j},
\cdot\cdot\cdot, \overset{\vee}\to{v_k}, \cdot\cdot\cdot, v_m)$ have opposite
signs and hence continuing the formula $(\diamond)$ we have
$$\align
&\leq \text{max}\{a_j \cdot |\det(v_j, v_1, \cdot\cdot\cdot,
\overset{\vee}\to{v_j},
\cdot\cdot\cdot, \overset{\vee}\to{v_k}, \cdot\cdot\cdot, v_m)|, a_k \cdot
|\det(v_k, v_1, \cdot\cdot\cdot, \overset{\vee}\to{v_j},
\cdot\cdot\cdot, \overset{\vee}\to{v_k}, \cdot\cdot\cdot, v_m)|\} \\
&\leq \text{max}\{a_j \cdot h_{\eta}, a_k \cdot h_{\eta}\} < h_{\eta}\\
\endalign$$
In the latter case (say,
$r_j = 0$ while $r_k \neq 0$), we have 
$$\det(v_k, v_1, \cdot\cdot\cdot,
\overset{\vee}\to{v_j},
\cdot\cdot\cdot, \overset{\vee}\to{v_k}, \cdot\cdot\cdot, v_m) = 0$$ 
and hence
continuing the formula $(\diamond)$ we obtain
$$\align
&= |a_j \cdot \det(v_j, v_1, \cdot\cdot\cdot, \overset{\vee}\to{v_j},
\cdot\cdot\cdot, \overset{\vee}\to{v_k}, \cdot\cdot\cdot, v_m)| \leq a_j \cdot
h_{\eta} < h_{\eta}. \\
\endalign$$

This completes the proof of the claim.

\vskip.1in

Obesrve that $\eta = \langle\rho_1, \cdot\cdot\cdot, \rho_n, \rho_{n+1},
\cdot\cdot\cdot, \rho_m\rangle$ and that a maximal cone
$\zeta'$ of
$\eta'$ has the form
$$\zeta' = \langle\rho', \rho_1, \cdot\cdot\cdot, \overset{\vee}\to{\rho_j},
\cdot\cdot\cdot, \rho_m\rangle \text{\ for\ some\ }j = 1, \cdot\cdot\cdot, m.$$
The only possible and old maximal $\pi$-independent face of $\zeta'$ is
\linebreak
$\langle\rho_1, \cdot\cdot\cdot, \overset{\vee}\to{\rho_j},
\cdot\cdot\cdot, \rho_m\rangle$ and hence the above claim
implies $\text{$\pi$-m.p.}(\zeta') \leq (h_{\eta},0,0,0)$.  Note
that $\text{$\pi$-m.p.}(\eta) \geq (h_{\eta},0,0,0)$ and that if the equality
holds then there is only one maximal $\pi$-independent face $\gamma \subset
\eta$ with $\text{$\pi$-mult}(\gamma) = h_{\eta}$ and hence we have possibly
only one maximal cone
$\zeta'$ of
$\eta'$, namely the one containing $\gamma$, having the $\pi$-multiplicity
profile equal to $(h_{\eta},0,0,0)$.  Therefore, we have either
$$\text{$\pi$-m.p.}(\eta) = (h_{\eta},1,*,*) = [(h_{\eta},1,*,*);1] >
[(h_{\eta},0,0,0),s] \geq \text{$\pi$-m.p.}(\eta')$$
or
$$\text{$\pi$-m.p.}(\eta) = (h_{\eta},0,0,0) = [(h_{\eta},0,0,0);1] \geq
\text{$\pi$-m.p.}(\eta').$$
If $\tau$ is contained in a maximal codimension one face $\gamma$ of
$\eta$ with $\text{$\pi$-mult}(\gamma) = h_{\eta}$, then in the latter
case we have the strict inequality.

\vskip.1in

This completes the proof of Proposition 5.5. 

\vskip.1in

As shown above, the star subdivision by a $\pi$-independent face behaves well
(choosing an appropriate division point in the interior) if it is codefinite
with respect to a $\pi$-dependent cone containing it, i.e., if it is codefinite
with respect to a circuit in its closed star.  In the following, we study how to
make a given
$\pi$-independent face codefinite with
respect to all the circuits in its closed star, after some specific
star subdivisions.

\vskip.1in

$\boxed{\roman{Step}\ 3}$

\vskip.1in

Let $\sigma = \langle\rho_1, \cdot\cdot\cdot, \rho_k\rangle$ be a simplicial
and $\pi$-strictly convex cone which is a circuit of dimension $k$ in
$N_{\Bbb Q}^+$, where the extremal rays $\rho_i$ of $\sigma$ are generated by
\linebreak
$(v_i,w_i) \in N_{\Bbb Q}^+ = N_{\Bbb Q} \oplus {\Bbb Q}$ with $v_i =
n(\pi(\rho_i))$, $i = 1, \cdot\cdot\cdot, k$ being the primitive
vectors in $N$.  Let $\tau$ be a codimension one face of
$\sigma$ with the maximum
$\pi$-multiplicity $h_{\sigma}$ among all the $\pi$-independent faces of
$\sigma$.  Say,
$$\tau = \tau_{\alpha} = \langle\rho_1, \cdot\cdot\cdot,
\overset{\vee}\to{\rho_{\alpha}}, \cdot\cdot\cdot, \rho_k\rangle.$$ 
We have the unique linear dependence
relation
$$(\natural) \hskip.1in \Sigma_{i = 1}^kr_iv_i = 0 \text{\ with\ the\
conditions\ }|r_{\alpha}| = 1 \text{\ and\ } r_1w_1 + \cdot\cdot\cdot + r_kw_k
> 0.$$  
We note that $0 < |r_i| \leq 1 \text{\ for\ } i = 1, \cdot\cdot\cdot, k$ where
$|r_i| = 1$ if and only if $\text{$\pi$-mult}(\tau_i) = h_{\sigma}$
for $\tau_i = \langle\rho_1, \cdot\cdot\cdot, \overset{\vee}\to{\rho_i},
\cdot\cdot\cdot, \rho_k\rangle$.

The first inequality $0 < |r_i|$ comes from the fact that $\sigma$ is a circuit
and the second inequality and the assertion about the equality comes from the
easy observation
$$\align
\text{$\pi$-mult}(\tau_i) &= \text{$\pi$-mult}(\langle\rho_1, \cdot\cdot\cdot,
\overset{\vee}\to{\rho_i},
\cdot\cdot\cdot, \rho_k\rangle) \\ 
&= |\det(v_1,
\cdot\cdot\cdot, \overset{\vee}\to{v_i}, \cdot\cdot\cdot, -
r_{\alpha}v_{\alpha} =
\Sigma_{j \neq \alpha}r_jv_j, \cdot\cdot\cdot, v_k)| \\ 
&= |r_i| \cdot |\det(v_1, \cdot\cdot\cdot, v_i, \cdot\cdot\cdot,
\overset{\vee}\to{v_{\alpha}},
\cdot\cdot\cdot, v_k)|
\\ 
&= |r_i| \cdot \text{$\pi$-mult}(\tau) \leq \text{$\pi$-mult}(\tau). \\
\endalign$$

Thus we conclude that the relation $(\natural)$ is independent of the choice
of a codimension one $\pi$-independent face $\tau$ of $\sigma$ as long as $\tau$
has the maximum $\pi$-multiplicity $h_{\sigma}$.

\vskip.1in

$\bold{Definition\ 5.6.}$ Let $\sigma$ be a circuit of $N_{\Bbb Q}^+$ as above. 
Then the negative (resp. positive) center point $\roman{Ctr}_-(\sigma)$ (resp.
$\roman{Ctr}_+(\sigma)$) of
$\sigma$ is defined to be
$$\roman{Ctr}_-(\sigma) = \Sigma_{r_i < 0}v_i \hskip.1in (\text{resp.\ }
\roman{Ctr}_+(\sigma) =
\Sigma_{r_i > 0}v_i).$$

\proclaim{Lemma 5.7} Let $\sigma$ be a circuit of $N_{\Bbb Q}^+$.  Then
$$\roman{Ctr}_-(\sigma), \roman{Ctr}_+(\sigma) \in
\roman{RelInt}(\pi(\sigma)).$$
\endproclaim

\demo{Proof}\enddemo We observe
$$\align
Ctr_-(\sigma) &= \Sigma_{r_i < 0}v_i = \Sigma_{r_i < 0}v_i + \Sigma_{i =
1}^kr_iv_i = \Sigma_{r_i > 0}r_iv_i + \Sigma_{r_i < 0}(1 +
r_i)v_i \\ 
&= (1 - \epsilon)\{\Sigma_{r_i < 0}v_i\} +
\epsilon\{\Sigma_{r_i > 0}r_iv_i + \Sigma_{r_i < 0}(1 + r_i)v_i\} \text{\ for\
} 0 < \epsilon < 1 \\ 
&= \Sigma_{i = 1}^k c_iv_i,\\
\endalign$$
where
$$c_i =
\left\{\aligned 
&= \epsilon r_i \text{\ when\ }r_i > 0 \\
&= 1 - \epsilon + \epsilon(1 + r_i) \text{\ when\ }r_i < 0
\endaligned 
\right.$$

Since $r_i \neq 0$ and $r_i \geq -1$ for all $i$, we
see
$$c_i > 0 \text{\ for\ all\ } i = 1, \cdot\cdot\cdot, k.$$
Thus we conclude
$$\roman{Ctr}_-(\sigma) \in \roman{RelInt}(\pi(\sigma)).$$
The argument for the statement $\roman{Ctr}_+(\sigma) \in
\roman{RelInt}(\pi(\sigma))$ is identical.

\vskip.1in

\proclaim{Lemma 5.8} Let $\sigma$ be a circuit in $N_{\Bbb Q}^+$ as
above with the negative center point $\roman{Ctr}_-(\sigma)$ (resp. the positive
center point $\roman{Ctr}_+(\sigma)$).  Let
$l_-$ (resp. $l_+$) be the ray generated by $\roman{Ctr}_-(\sigma)$ (resp.
$\roman{Ctr}_+(\sigma)$) and $\sigma' =
\roman{Mid}(\sigma,l_-)
\cdot
\sigma$ (resp. $\sigma' = \roman{Mid}(\sigma,l_+) \cdot
\sigma$) be the subdivision of $\sigma$ by the midray $\roman{Mid}(\sigma,l_-)$
(resp. $\roman{Mid}(\sigma,l_+)$).  Then every codimension one face $\zeta$ of
$\sigma$ with the maximum
$\pi$-multiplicity
$h_{\sigma}$ (which stays unchanged through the subdivision and hence can be
considered a face
$\zeta \in
\sigma'$) is codefinite with respect to the (unique) maximal cone in the closed
star of $\zeta$ in $\sigma'$.
\endproclaim

\demo{Proof}\enddemo We use the same notation as above.  We only prove the
statement for the negative center as the proof is identical for the positive
center. 

\vskip.1in

Observe fisrt that we have
$$\roman{Mid}(\sigma,l_-) \in \roman{RelInt}(\sigma),$$ 
since $\roman{Ctr}_-(\sigma) \in
\roman{RelInt}(\pi(\sigma))$ by Lemma 5.7.  Therefore, the
star subdivision with respect to $\roman{Mid}(\sigma,l_-)$ does not affect
$\zeta$, i.e., $\zeta \in \sigma'$.

\vskip.1in

Observe secondly (say, $\zeta = \zeta_j = \langle\rho_1, \cdot\cdot\cdot,
\overset{\vee}\to{\rho_j},
\cdot\cdot\cdot, \rho_k\rangle$) that
$$\zeta_j \text{\ has\ maximal\ $\pi$-multiplicity\ }\Longleftrightarrow |r_j|
= 1.$$

\newpage

Case: $r_j = 1$.

\vskip.1in

In this case, since $Ctr_-(\sigma) = \Sigma_{r_i < 0}v_i$ and since the maximal
cone $\sigma_j$ containing $\zeta_j$ in $\sigma'$ is of the form $\sigma_j =
\langle\roman{Mid}(\sigma,l_-), \rho_1, \cdot\cdot\cdot,
\overset{\vee}\to{\rho_j}, \cdot\cdot\cdot, \rho_k\rangle$, 
the linear relation for $\sigma_j$ is given by
$$\roman{Ctr}_-(\sigma) - \Sigma_{r_i < 0}v_i = 0.$$
As $\zeta_j$ contains only the extremal rays corresponding to the $v_i$, which
have the same sign (or 0) in the linear relation, $\zeta_j$ is codefinite with
respect to the (unique) maximal cone $\sigma_j$ in the closed star of $\zeta_j$
in
$\sigma'$.

\vskip.1in

Case: $r_j = -1$.

\vskip.1in

In this case, since $\roman{Ctr}_-(\sigma) = \Sigma_{r_i > 0}r_iv_i + \Sigma_{-1
< r_i < 0}(1 + r_i)v_i$ and since the maximal cone $\sigma_j$ containing
$\zeta_j$ is of the form $\sigma_j = \langle\roman{Mid}(\sigma,l_-), \rho_1,
\cdot\cdot\cdot,
\overset{\vee}\to{\rho_j}, \cdot\cdot\cdot, \rho_k\rangle$, 
the linear relation for $\sigma_j$ is given by
$$\roman{Ctr}_-(\sigma) - \Sigma_{r_i > 0}r_iv_i - \Sigma_{-1 < r_i < 0}(1 +
r_i)v_i = 0.$$ 
As $\zeta_j$ contains only the extremal rays corresponding to the $v_i$, which
have the same sign (or 0) in the linear relation, $\zeta_j$ is codefinite with
respect to the (unique) maximal cone in the closed star $\sigma_j$ of $\zeta_j$
in
$\sigma'$.

\vskip.1in

This completes the proof of Lemma 5.8.

\vskip.1in

This lemma suggests that we should use the star subdivision by the negative or
positive center of a circuit to achieve codefinireness of a face $\tau$ in
order to bring the situation of Proposition 5.4 in Step 2.  But the lemma
only achieves the codefiniteness for a face $\tau$ which is contained in a
maximal $\pi$-independent face with the maximum $\pi$-multiplicity but does
not analyze the behavior of the $\pi$-multiplicity profile.  In our process of
$\pi$-desingularization, we need to achieve codefiniteness for a face $\tau$
which is not contanied in a maximal $\pi$-independent face with the maximum
$\pi$-multiplicity and the analysis of the $\pi$-multiplicity profile is
crucial.  Both of these needs are fulfilled by the following proposition,
which is at the technical heart of this section.

\proclaim{Proposition 5.9} Let $\sigma$ be a circuit of $\dim \sigma > 2$ in
$N_{\Bbb Q}^+$.  Then by choosing $\sigma'$ to be either the star subdivision of
$\sigma$ corresponding to the negative center point or the one by the positive
center point, i.e.,
$$\sigma' = \roman{Mid}(\sigma,l_-) \cdot \sigma \text{\ or\ }
\roman{Mid}(\sigma,l_+)
\cdot \sigma$$
where $l_-$ (resp. $l_+$) is the ray generated by the negative (resp.
positive) center point
$\roman{Ctr}_-(\sigma)$ (resp. $\roman{Ctr}_+(\sigma)$), we see $\sigma'$
satisfies one of the following:

\vskip.1in

A: Every maximal cone $\delta'$ of $\sigma'$ has the $\pi$-multiplicity profile
strictly smaller than that of $\sigma$, i.e.,
$$\text{$\pi$-\rm{m.p.}}(\delta') < \text{$\pi$-\rm{m.p.}}(\sigma).$$
In particular, we have
$$\text{$\pi$-\rm{m.p.}}(\sigma') < \text{$\pi$-\rm{m.p.}}(\sigma).$$

\vskip.1in

B: Every maximal cone $\delta'$ of $\sigma'$, except for one maximal cone
$\kappa'$, has the $\pi$-multiplicity profile strictly smaller than that of
$\sigma$, i.e., 
$$\text{$\pi$-\rm{m.p.}}(\delta') < \text{$\pi$-\rm{m.p.}}(\sigma)$$
and the exceptional maximal cone $\kappa'$ has the same $\pi$-multiplicity
profile as that of $\sigma$, i.e.,
$$\text{$\pi$-\rm{m.p.}}(\kappa') = \text{$\pi$-\rm{m.p.}}(\sigma).$$
In particular, we have
$$\text{$\pi$-\rm{m.p.}}(\sigma') = \text{$\pi$-\rm{m.p.}}(\sigma).$$
Moreover, there exists a maximal $\pi$-independent face $\gamma'$ of $\kappa'$
such that

\ \ (\rm{B-o}) $\gamma'$ is also a face of $\sigma$ (i.e., $\gamma'$ remains
untouched under th subdivision),

\ \ (\rm{B-ii}) $\gamma'$ has the maximum $\pi$-multiplicity,
i.e., $\text{$\pi$-mult}(\gamma') = h_{\sigma'} = h_{\sigma}$,

\ \ (\rm{B-iii}) $\gamma'$ is codefinite with respect to $\kappa'$.

\endproclaim

\demo{Proof}\enddemo Let $\sigma = \langle\rho_1, \cdot\cdot\cdot,
\rho_k\rangle$, where the extremal rays $\rho_i$ are generated by $(v_i,w_i) \in
N_{\Bbb Q}^+$ with $v_i = n(\pi(\rho_i))$, $i = 1, \cdot\cdot\cdot,
k,$ being the primitive vectors for the projections.  

Let $\Sigma r_iv_i = 0$ be the nontrivial linear relation so that $\Sigma
r_iw_i > 0$ and 
$$|r_i| = 1 \Longleftrightarrow \text{$\pi$-mult}(\tau_i) =
h_{\sigma} \text{\ for\ }\tau_i = \langle\rho_1, \cdot\cdot\cdot,
\overset{\vee}\to{\rho_i},
\cdot\cdot\cdot, \rho_k\rangle.$$

Note that the maximal cones $\sigma'_i$ of $\sigma'$ are of the form
$$\sigma'_i = \langle\rho_0, \rho_1, \cdot\cdot\cdot, \overset{\vee}\to{\rho_i},
\cdot\cdot, \rho_k\rangle$$
where $\rho_0$ is the midray $\roman{Mid}(\sigma,l_-)$ or
$\roman{Mid}(\sigma,l_+)$ depending on the choice of the negative or positive
center point.

We compute the $\pi$-multiplicity of the maximal faces $\tau'_{ij}$ of
$\sigma'_i$
$$\tau'_{ij} = \langle\rho_0, \rho_1, \cdot\cdot\cdot,
\overset{\vee}\to{\rho_i},
\cdot\cdot\cdot, \overset{\vee}\to{\rho_j}, \cdot\cdot\cdot, \rho_k\rangle$$
as follows:

\vskip.1in

Case of the negative center point: $\rho_0 = \roman{Mid}(\sigma,l_-)$.

\vskip.1in

We let $e_- \in {\Bbb N}$ be the integer such that
$\Sigma_{r_{\alpha} < 0}v_{\alpha} = e_- \cdot n(\pi(\rho_0))$ with
$n(\pi(\rho_0))$ being the primitive vector.

\newpage

\ \ Subcase $r_i > 0$:

$$\align
\text{$\pi$-mult}(\tau'_{i0}) &= \text{$\pi$-mult}(\tau_i) \\
\text{$\pi$-mult}(\tau'_{ij}) &= \frac{1}{e_-}|\det(\Sigma_{r_{\alpha} <
0}v_{\alpha},v_1, \cdot\cdot\cdot, \overset{\vee}\to{v_i}, \cdot\cdot\cdot,
\overset{\vee}\to{v_j}, \cdot\cdot\cdot, v_k)|\\
&= \left\{\aligned
0 \hskip.65in &\text{\ when\ }r_j > 0 \\
\frac{1}{e_-}\text{$\pi$-mult}(\tau_i) &\text{\ when\ }r_j < 0. \\
\endaligned
\right.\\
\endalign$$

\ \ Subcase $r_i < 0$:

$$\align
\text{$\pi$-mult}(\tau'_{i0}) &= \text{$\pi$-mult}(\tau_i) \\
\text{$\pi$-mult}(\tau'_{ij}) &= \frac{1}{e_-}|\det(\Sigma_{r_{\alpha} <
0}v_{\alpha},v_1, \cdot\cdot\cdot, \overset{\vee}\to{v_i}, \cdot\cdot\cdot,
\overset{\vee}\to{v_j}, \cdot\cdot\cdot, v_k)|\\
&= \left\{\aligned
\frac{1}{e_-}\text{$\pi$-mult}(\tau_j) \hskip.87in &\text{\ when\ }r_j > 0 \\
\frac{1}{e_-}|\text{$\pi$-mult}(\tau_j) - \text{$\pi$-mult}(\tau_i)| &\text{\
when\ }r_j < 0. \\
\endaligned
\right.\\
\endalign$$

Note that $\det(v_i,v_1, \cdot\cdot\cdot, \overset{\vee}\to{v_i},
\cdot\cdot\cdot, \overset{\vee}\to{v_j}, \cdot\cdot\cdot, v_k)$ and
$\det(v_j,v_1, \cdot\cdot\cdot, \overset{\vee}\to{v_i},
\cdot\cdot\cdot, \overset{\vee}\to{v_j}, \cdot\cdot\cdot, v_k)$ have opposite
signs, since 
$$\align
0 &= \det(\Sigma r_{\alpha}v_{\alpha}, v_1, \cdot\cdot\cdot,
\overset{\vee}\to{v_i},
\cdot\cdot\cdot, \overset{\vee}\to{v_j}, \cdot\cdot\cdot, v_k) \\
&= r_i(v_i, v_1, \cdot\cdot\cdot, \overset{\vee}\to{v_i},
\cdot\cdot\cdot, \overset{\vee}\to{v_j}, \cdot\cdot\cdot, v_k) + r_j(v_j,v_1, \cdot\cdot\cdot, \overset{\vee}\to{v_i},
\cdot\cdot\cdot, \overset{\vee}\to{v_j}, \cdot\cdot\cdot, v_k) = 0.\\
\endalign$$
 
\vskip.1in

Symmetrically we compute the other case.

\vskip.1in

Case of the positive center point: $\rho_0 = \roman{Mid}(\sigma,l_+)$.

\vskip.1in

We let $e_+ \in {\Bbb N}$ be the integer such that
$\Sigma_{r_{\alpha} > 0}v_{\alpha} = e_+ \cdot n(\pi(\rho_0))$ with
$n(\pi(\rho_0))$ being the primitive vector.

\vskip.1in

\ \ Subcase $r_i < 0$:

$$\align
\text{$\pi$-mult}(\tau'_{i0}) &= \text{$\pi$-mult}(\tau_i) \\
\text{$\pi$-mult}(\tau'_{ij}) &= \frac{1}{e_+}|\det(\Sigma_{r_{\alpha} >
0}v_{\alpha},v_1, \cdot\cdot\cdot, \overset{\vee}\to{v_i}, \cdot\cdot\cdot,
\overset{\vee}\to{v_j}, \cdot\cdot\cdot, v_k)|\\
&= \left\{\aligned
0 \hskip.65in &\text{\ when\ }r_j < 0 \\
\frac{1}{e_+}\text{$\pi$-mult}(\tau_i) &\text{\ when\ }r_j > 0. \\
\endaligned
\right.\\
\endalign$$

\newpage

\ \ Subcase $r_i > 0$:

$$\align
\text{$\pi$-mult}(\tau'_{i0}) &= \text{$\pi$-mult}(\tau_i) \\
\text{$\pi$-mult}(\tau'_{ij}) &= \frac{1}{e_+}|\det(\Sigma_{r_{\alpha} >
0}v_{\alpha},v_1, \cdot\cdot\cdot, \overset{\vee}\to{v_i}, \cdot\cdot\cdot,
\overset{\vee}\to{v_j}, \cdot\cdot\cdot, v_k)|\\
&= \left\{\aligned
\frac{1}{e_+}\text{$\pi$-mult}(\tau_j) \hskip.87in &\text{\ when\ }r_j < 0 \\
\frac{1}{e_+}|\text{$\pi$-mult}(\tau_j) - \text{$\pi$-mult}(\tau_i)| &\text{\
when\ }r_j > 0. \\
\endaligned
\right.\\
\endalign$$

Using this computation, we can now easily derive the conclusion of the
proposition dividing it into the cases according to the cardinalities of the
following sets:
$$I^+ = \{i;r_i > 0\}, I^- = \{i;r_i < 0\}, I^+_1 = \{i;r_i = 1\}, I^-_1 =
\{i;r_i = -1\}.$$

\vskip.2in

\noindent Case:\ $2 \leq \#I^-_1 \leq \#I^+_1$

\vskip.1in

In this case, we choose the negative center point and let $\rho_0 =
\roman{Mid}(\sigma,l_-)$.

When $0 < r_i < 1$, we have $a_{\sigma'_i} = h_{\sigma'_i} < h_{\sigma} =
a_{\sigma}$.

When $0 < r_i = 1$, we have $a_{\sigma'_i} = h_{\sigma'_i} = h_{\sigma}
= a_{\sigma}$.  If $e_- > 1$, then $r_{\sigma_i'} = 1$ and hence $b_{\sigma_i'}
= 0 < 1 = b_{\sigma}$.  If $e_- = 1$, then $r_{\sigma'_i}
\geq
\# I^-_1
\geq 2$ and hence
$b_{\sigma'_i} = 1 = b_{\sigma}$.  But $c_{\sigma'_i} = k_{\sigma'_i} <
k_{\sigma} = c_{\sigma}$, since
$\text{$\pi$-mult}(\tau'_{ij}) = 0$ for $j \in I^+ \supset I^+_1$.

When $-1 < r_i < 0$ and $e_- > 1$, we have $a_{\sigma_i'} = h_{\sigma_i'} <
h_{\sigma} = a_{\sigma}$.

When $- 1 = r_i < 0$ and $e_- > 1$, we have $a_{\sigma_i'} = h_{\sigma_i'} =
h_{\sigma} = a_{\sigma}$.  But $r_{\sigma_i'} = 1$ and hence $b_{\sigma_i'} = 0
< 1 = b_{\sigma}$.

When $- 1 \leq r_i < 0$ and $e_- = 1$, we have $a_{\sigma'_i} = h_{\sigma'_i} =
h_{\sigma} = a_{\sigma}$, \linebreak
$r_{\sigma'_i} \geq \# I^+_1 \geq 2$ and
hence
$b_{\sigma'_i} = 1 = b_{\sigma}$.  We also have $c_{\sigma'_i} = k_{\sigma'_i}
\leq k_{\sigma} = c_{\sigma}$, since
$\sigma$ is a circuit.  But $d_{\sigma'_i} = r_{\sigma'_i} = r_{\sigma} - \#
I^-_1 + 1 < r_{\sigma} = d_{\sigma}$, since
$\tau'_{ij} = |\text{$\pi$-mult}(\tau_j) - \text{$\pi$-mult}(\tau_i)| <
h_{\sigma}$ for $j
\in I^-_1, \hskip.1in j \neq i$.

Thus we have 
$$\text{$\pi$-m.p.}(\sigma_i') < \text{$\pi$-m.p.}(\sigma)$$ 
for all the maximal cones $\sigma_i'$ of $\sigma'$.

Therefore, in this case with the choice of the negative center we conclude we
are in Case A and
$$\text{$\pi$-m.p.}(\sigma') < \text{$\pi$-m.p.}(\sigma).$$

\vskip.1in

\noindent Case:\ $1 = \#I^-_1 < \#I^+_1$

\vskip.1in

In this case, we choose the negative center point and let $\rho_0 =
\roman{Mid}(\sigma,l_-)$.

When $0 < r_i < 1$, we have $a_{\sigma'_i} = h_{\sigma'_i} < h_{\sigma} =
a_{\sigma}$ and hence \linebreak
$\text{$\pi$-m.p.}(\sigma'_i) <
\text{$\pi$-m.p.}(\sigma)$.

When $0 < r_i = 1$, we have $a_{\sigma'_i} = h_{\sigma'_i} = h_{\sigma}
= a_{\sigma}$.  If $\# I^-  = \# I^-_1 = 1$ or $e_- > 1$, then $r_{\sigma'_i} =
1$ and hence
$b_{\sigma'_i} = 0 < 1 = b_{\sigma}$.  If $\# I^- > 1$ and $e_- = 1$, then
$r_{\sigma'_i} > 1$ and hence $b_{\sigma'_i} = 1 = b_{\sigma}$.  But we
have
$c_{\sigma'_i} = k_{\sigma'_i} < k_{\sigma} = c_{\sigma}$, since
$\text{$\pi$-mult}(\tau'_{ij}) = 0$ for $j \in I^+ \supset I^+_1$.  Thus we
have $\text{$\pi$-m.p.}(\sigma'_i) < \text{$\pi$-m.p.}(\sigma)$.

When $-1 < r_i < 0$ and $e_- > 1$, we have $a_{\sigma'_i} = h_{\sigma'_i} <
h_{\sigma} = a_{\sigma}$ and hence $\text{$\pi$-m.p.}(\sigma'_i) <
\text{$\pi$-m.p.}(\sigma)$.

When $-1 = r_i < 0$ and $e_- > 1$, we have $a_{\sigma'_i} = h_{\sigma'_i} = h_{\sigma}
= a_{\sigma}$.  But $r_{\sigma'} = 1$ and hence $b_{\sigma'} = 0
< 1 = b_{\sigma}$.  Thus we have $\text{$\pi$-m.p.}(\sigma'_i) <
\text{$\pi$-m.p.}(\sigma)$.

When $- 1 \leq r_i < 0$ and $e_- = 1$, we have $a_{\sigma'_i} = h_{\sigma'_i} =
h_{\sigma} = a_{\sigma}$, $r_{\sigma'_i} \geq \# I^+_1 \geq 2$ and hence
$b_{\sigma'_i} = 1 = b_{\sigma}$.  We also have $c_{\sigma'_i} = k_{\sigma'_i}
\leq k_{\sigma} = c_{\sigma}$, since
$\sigma$ is a circuit.  But $d_{\sigma'_i} = r_{\sigma'_i} = r_{\sigma} - \#
I^-_1 + 1 \leq r_{\sigma} = d_{\sigma}$, since
$\tau'_{ij} = |\text{$\pi$-mult}(\tau_j) - \text{$\pi$-mult}(\tau_i)| <
h_{\sigma}$ for $j
\in I^-_1, \hskip.1in j \neq i$.  Thus we have $\text{$\pi$-m.p.}(\sigma'_i)
\leq \text{$\pi$-m.p.}(\sigma)$. The equality holds only when $r_i = -1$ with
$i$ being the sole member of
$I^-_1$, in which case the face $\tau'_{i0} = \tau_i$  has the maximum
$\pi$-multiplicity $h_{\sigma}$ and it is codefinite with respect to
$\sigma'_i$ by Lemma 5.8.

Therefore, in this case with the choice of the negative center we conclude
that we are in Case A and 
$$\text{$\pi$-m.p.}(\sigma') < \text{$\pi$-m.p.}(\sigma) \text{\ if\ }e_- >
1$$ 
and that we are in Case B and
$$\text{$\pi$-m.p.}(\sigma') = \text{$\pi$-m.p.}(\sigma) \text{\ if\ }e_- =
1.$$

\vskip.1in

\noindent Case:\ $1 = \#I^-_1 = \#I^+_1 < \# I^+$

\vskip.1in

In this case, we choose the negative center point and let $\rho_0 =
\roman{Mid}(\sigma,l_-)$.

When $0 < r_i < 1$, we have $a_{\sigma'_i} = h_{\sigma'_i} < h_{\sigma} =
a_{\sigma}$ and hence \linebreak
$\text{$\pi$-m.p.}(\sigma'_i) <
\text{$\pi$-m.p.}(\sigma)$.

When $0 < r_i = 1$ and $e_- > 1$, we have $a_{\sigma'_i} = h_{\sigma'_i} =
h_{\sigma} = a_{\sigma}$.  But $r_{\sigma'} = 1$ and hence $b_{\sigma_i'} = 0
< 1 = b_{\sigma}$.  Thus we have $\text{$\pi$-m.p.}(\sigma'_i) <
\text{$\pi$-m.p.}(\sigma)$.

When $0 < r_i = 1$ and $e_- = 1$, we have $a_{\sigma'_i} = h_{\sigma'_i} =
h_{\sigma} = a_{\sigma}$.  If $\# I^-  = \# I^-_1 = 1$, then $r_{\sigma'_i} =
1$ and hence
$b_{\sigma'_i} = 0 < 1 = b_{\sigma}$.  If $\# I^- > 1$, then $r_{\sigma'_i} >
1$ and hence $b_{\sigma'_i} = 1 = b_{\sigma}$.  But we have
$c_{\sigma'_i} = k_{\sigma'_i} < k_{\sigma} = c_{\sigma}$, since
$\text{$\pi$-mult}(\tau'_{ij}) = 0$ for $j \in I^+ \supset I^+_1, \hskip.1in
j \neq i$.  Thus we have $\text{$\pi$-m.p.}(\sigma'_i) <
\text{$\pi$-m.p.}(\sigma)$.

When $-1 < r_i < 0$ and $e_- > 1$, we have $a_{\sigma'_i} = h_{\sigma'_i} <
h_{\sigma} = a_{\sigma}$ and hence $\text{$\pi$-m.p.}(\sigma'_i) <
\text{$\pi$-m.p.}(\sigma)$.

When $- 1 < r_i < 0$ and $e_- = 1$, we have $a_{\sigma'_i} = h_{\sigma'_i} =
h_{\sigma} = a_{\sigma}$.  But $r_{\sigma'_i} = 1$ and hence $b_{\sigma'_i} = 0
< 1 = b_{\sigma}$.  Thus we have $\text{$\pi$-m.p.}(\sigma'_i) <
\text{$\pi$-m.p.}(\sigma)$.

When $r_i = -1$ and $e_- > 1$, we have $a_{\sigma'_i} = h_{\sigma'_i} =
h_{\sigma} = a_{\sigma}$.  But $r_{\sigma'_i} = 1$ and hence $b_{\sigma'_i} = 0
< 1 = b_{\sigma}$.  Thus we have $\text{$\pi$-m.p.}(\sigma'_i) <
\text{$\pi$-m.p.}(\sigma)$.

When $r_i = -1$ and $e_- = 1$, $i$ is the sole member of $I^-_1$ and we have
$a_{\sigma'_i} = h_{\sigma'_i} = h_{\sigma} = a_{\sigma}$, $r_{\sigma'_i} = 2$
and hence
$b_{\sigma'_i} = 1 = b_{\sigma}$.  Moreover, we have $c_{\sigma'_i} =
k_{\sigma'_i} = k_{\sigma} = c_{\sigma}$ and  $d_{\sigma'_i} = r_{\sigma'_i} = 2
= r_{\sigma} = d_{\sigma_i}$.  Thus we have $\text{$\pi$-m.p.}(\sigma'_i) =
\text{$\pi$-m.p.}(\sigma)$.  The face
$\tau'_{i0} = \tau_i$  has the maximum
$\pi$-multiplicity $h_{\sigma}$ and it is codefinite with respect to
$\sigma'_i$ by Lemma 5.8.

Therefore, in this case with the choice of the negative center we conclude
that we are in Case A and
$$\text{$\pi$-m.p.}(\sigma') < \text{$\pi$-m.p.}(\sigma) \text{\ if\ }e_- >
1$$
and that we are in Case B and
$$\text{$\pi$-m.p.}(\sigma') = \text{$\pi$-m.p.}(\sigma) \text{\ if\ }e_- =
1.$$

\noindent Case:\ $0 = \#I^-_1 < 2 \leq \#I^+_1$

\vskip.1in

In this case, we choose the positive center point and let $\rho_0 =
\roman{Mid}(\sigma,l_+)$.

When $-1 < r_i < 0$, we have $a_{\sigma'_i} = h_{\sigma'_i} < h_{\sigma} =
a_{\sigma}$ and hence \linebreak
$\text{$\pi$-m.p.}(\sigma'_i) <
\text{$\pi$-m.p.}(\sigma)$.

When $0 < r_i < 1$, we have $a_{\sigma'_i} = h_{\sigma'_i} < h_{\sigma}
= a_{\sigma}$ and hence \linebreak
$\text{$\pi$-m.p.}(\sigma'_i) <
\text{$\pi$-m.p.}(\sigma)$.

When $r_i = 1$, we have $a_{\sigma'_i} = h_{\sigma'_i} = h_{\sigma}
= a_{\sigma}$.  But $r_{\sigma'_i} = 1$ and hence \linebreak
$b_{\sigma'_i} - 0 <
1 = b_{\sigma}$.  Thus we have $\text{$\pi$-m.p.}(\sigma'_i) <
\text{$\pi$-m.p.}(\sigma)$.

Therefore, in this case with the choice of the positive center we conclude
that we are in Case A and
$$\text{$\pi$-m.p.}(\sigma') < \text{$\pi$-m.p.}(\sigma).$$

\vskip.1in

\noindent Case:\ $0 = \#I^-_1 < 1 = \# I^+_1$

\vskip.1in

In this case, we choose the positive center point and let $\rho_0 =
\roman{Mid}(\sigma,l_+)$.

When $-1 < r_i < 0$, we have $a_{\sigma'_i} = h_{\sigma'_i} < h_{\sigma} =
a_{\sigma}$ and hence \linebreak
$\text{$\pi$-m.p.}(\sigma'_i) <
\text{$\pi$-m.p.}(\sigma)$.

When $0 < r_i < 1$, we have $a_{\sigma'_i} = h_{\sigma'_i} < h_{\sigma}
= a_{\sigma}$ and hence \linebreak
$\text{$\pi$-m.p.}(\sigma'_i) <
\text{$\pi$-m.p.}(\sigma)$.

When $r_i = 1$, i.e., $i$ is the sole member of $I^+_1$, we have $a_{\sigma'_i}
= h_{\sigma'_i} = h_{\sigma} = a_{\sigma}$, $r_{\sigma'_i} = 1$ and hence
$b_{\sigma'_i} = 0 = b_{\sigma}$.  Moreover, we have $c_{\sigma'_i} =
k_{\sigma'_i} = k_{\sigma} = c_{\sigma}$ and  $d_{\sigma'_i} = r_{\sigma'_i} = 1
= r_{\sigma} = d_{\sigma}$.  Thus we have $\text{$\pi$-m.p.}(\sigma'_i) =
\text{$\pi$-m.p.}(\sigma)$.  The face $\tau'_{i0} = \tau_i$ has the maximum
$\pi$-multiplicity $h_{\sigma}$ and it is codefinite with respect to
$\sigma'_i$ by Lemma 5.8.

Therefore, in this case with the choice of the positive center we conclude we
are in Case B and
$$\text{$\pi$-m.p.}(\sigma') = \text{$\pi$-m.p.}(\sigma).$$
  
\vskip.2in

Symmetrically, we also conclude:

\vskip.1in

\noindent Case:\ $2 \leq \# I^+_1 \leq \# I^-_1$

\vskip.1in

With the choice of the positive center point, we are in Case A.

\vskip.1in

\noindent Case:\ $1 = \#I^+_1 < \#I^-_1$

\vskip.1in

With the choice of the positive center point, we are in Case A if $e_+ > 1$
and \linebreak
in Case B if $e_+ = 1$.

\vskip.1in

\noindent Case:\ $1 = \#I^+_1 = \#I^-_1 < \# I^-$

\vskip.1in

With the choice of the positive center point, we are in Case A if $e_+ > 1$
and \linebreak
in Case B if $e_+ = 1$.

\newpage

\noindent Case:\ $0 = \#I^+_1 < 2 \leq \#I^-_1$

\vskip.1in

With the choice of the negative center point, we are in Case A.

\vskip.1in

\noindent Case:\ $0 = \#I^+_1 < 1 = \# I^-_1$

\vskip.1in

With the choice of the negative center point, we are in Case B.

\vskip.2in

Since the above cases exhaust all the possibilities, we complete the proof for
Proposition 5.9.

\vskip.1in

The next lemma shows that the $\pi$-multiplicity of a cone can be computed
easily from that of the unique circuit contained in it. 

\proclaim{Lemma 5.10} Let $\sigma$ be a circuit in a simplicial
cobordism $\Sigma$ in $N_{\Bbb Q}^+$ and
$\eta$ be a maximal cone in $\overline{\roman{Star}(\sigma)}$.  Then any maximal
$\pi$-independent face $\gamma$ of $\eta$ is of the form
$$\gamma = \tau + \nu$$
where $\tau = \gamma \cap \sigma$ is a maximal $\pi$-independent face of
$\sigma$ and
$\nu$ is the unique maximal cone of $\roman{link}_{\eta}(\sigma)$.  

Moreover, there exists $e \in {\Bbb N}$ such
that for any $\gamma$ as above (once $\eta$ is fixed) we have the formula
$$\text{$\pi$-\rm{mult}}(\gamma) = \text{$\pi$-\rm{mult}}(\tau) \cdot e.$$
In particular, we have
$$\text{$\pi$-\rm{m.p.}}(\eta) = (a_{\eta},b_{\eta},c_{\eta},d_{\eta}) =
(e \cdot a_{\sigma},b_{\sigma},c_{\sigma},d_{\sigma}).$$
\endproclaim

\demo{Proof}\enddemo Let $\sigma = \langle\rho_1, \cdot\cdot\cdot,
\rho_k\rangle$ and $\eta = \langle\rho_1, \cdot\cdot\cdot, \rho_k, \rho_{k+1},
\cdot\cdot\cdot,
\rho_l\rangle$ be generated by the extremal rays $\rho_i$ with the coresponding
primitive vectors of the projections \linebreak
$v_i = n(\pi(\rho_i)) \in N$. 
Then a maximal
$\pi$-independent face
$\gamma$ of $\eta$ is of the form
$$\gamma = \langle\rho_1, \cdot\cdot\cdot, \overset{\vee}\to{\rho_j},
\cdot\cdot\cdot,
\rho_k,
\rho_{k+1},
\cdot\cdot\cdot,
\rho_l\rangle = \tau + \nu$$
where $\tau = \langle\rho_1, \cdot\cdot\cdot, \overset{\vee}\to{\rho_j},
\cdot\cdot\cdot,
\rho_k\rangle = \gamma \cap \sigma$ and $\nu = \langle\rho_{k+1},
\cdot\cdot\cdot,
\rho_l\rangle$ is the unique maximal cone of $\roman{link}_{\eta}(\sigma)$. 
This proves the first assertion.

\vskip.1in

For ``Moreover" part, we have the exact sequence
$$0 \rightarrow L \rightarrow N_{\eta} \rightarrow Q \rightarrow 0$$
where 
$L = \roman{span}_{\Bbb Q}(\pi(\sigma)) \cap N$, $N_{\eta} = \roman{span}_{\Bbb
Q}(\pi(\eta))
\cap N$ and $Q$ is the cokernel, which is
torsion free and hence a free ${\Bbb Z}$-module.
Take a ${\Bbb Z}$-basis $\{u_1, \cdot\cdot\cdot , u_{k-1}, u_{k+1},
\cdot\cdot\cdot, u_l\}$ of $N_{\eta}$ so that $\{u_1, \cdot\cdot,
u_{k-1}\}$ is a ${\Bbb Z}$-basis of $L$ and $\{u_{k+1}, \cdot\cdot\cdot,
u_l\}$ maps to a ${\Bbb Z}$-basis of $Q$.  With respect to this basis of
$N_{\eta}$, the
$\pi$-multiplicity of
$\gamma$ can be computed
$$\text{$\pi$-mult}(\gamma) = \det \left(\matrix
A & B \\
0 & E \\
\endmatrix
\right) = \det A \cdot \det E = \text{$\pi$-mult}(\tau) \cdot e,$$
where
$$(v_1, \cdot\cdot\cdot, \overset{\vee}\to{v_j}, \cdot\cdot\cdot, v_k) =
\left(\matrix A \\
0 \\
\endmatrix
\right)
\hskip.1in \text{and} \hskip.1in
(v_{k+1}, \cdot\cdot\cdot, v_l) = \left(\matrix
B \\
E \\
\endmatrix
\right).$$
This completes the proof of the lemma.

\vskip.1in

Now it is easy to see the following main consequence of Step 3.

\proclaim{Corollary 5.11} Let $\tau$ be a $\pi$-independent face contained in
the closed star $\overline{\roman{Star}(\sigma)}$ of a circuit $\sigma$.  Then
there exists $\{\overline{\roman{Star}(\sigma)}\}^{\circ}$ obtained by a
succession of star subdivisions by the negative or positive center points of
the circuits (of the intermediate subdivisions) inside of $\sigma$ such that

$(5.11.1)$ the $\pi$-multiplicity profile does not increase, i.e.,
$$\text{$\pi$-\rm{m.p.}}(\{\overline{\roman{Star}(\sigma)}\}^{\circ}) \leq
\text{$\pi$-\rm{m.p.}}(\overline{\roman{Star}(\sigma)}),$$

$(5.11.2)$ $\tau$ is a face of
$\{\overline{\roman{Star}(\sigma)}\}^{\circ}$ and
$\tau$ is codefinite with respect to every cone $\nu \in
\{\overline{\roman{Star}(\sigma)}\}^{\circ}$ containing $\tau$.

\endproclaim

\demo{Proof}\enddemo We prove the assertion by induction on the
$\pi$-multiplicity profile of $\sigma$.

\vskip.1in

If $\text{$\pi$-m.p.}(\sigma) = (1,*,*,*)$, then by taking
$\{\overline{\roman{Star}(\sigma)}\}^{\circ}$ to be the star subdivision
corresponding to the negative or positive center point of $\sigma$, we see
easily that the condition (5.11.1) is satisfied, while the condition (5.11.2) is
a consequence of Lemma 5.8.

\vskip.1in

We assume that the assertion holds for the case with the $\pi$-multiplicity
profile smaller than $\text{$\pi$-m.p.}(\sigma)$.  If $\dim \sigma = 2$,
then $\tau$ is already codefinite with respect to $\sigma$ and there is
nothing more to prove.  So we may assume $\dim \sigma > 2$.  We take the
star subdivision by the negative or positive center of $\sigma$, according to
Proposition 5.9, so that either the case A or the case B holds and hence the
$\pi$-multiplicity does not increase.

\vskip.1in

If the case A holds, noting that the circuits of all the maximal cones of
the star subdivision are contained in $\sigma$ we see the assertion holds
immediately by the induction hypothesis, since all the maximal cones have
the $\pi$-multiplicity profile strictly smaller than
$\text{$\pi$-m.p.}(\sigma)$.   

\vskip.1in

Suppose the case B holds.  If $\tau \cap \sigma$ is contained in
$\kappa'$, then $\tau \cap \sigma$ is necessarily contained in $\gamma'$ and
hence codefinite with respect to $\kappa'$.  The other maximal cones have the
$\pi$-multiplicity profile strictly smaller than $\text{$\pi$-m.p.}(\sigma)$
and the assertion again holds by the induction hypothesis.

\vskip.1in

This completes the proof of Corollary 5.11 and Step 3.

\vskip.2in

Now we discuss Step 4.

\vskip.1in

$\boxed{\roman{Step}\ 4}$

\vskip.1in

We start from a simplicial cobordism $\Sigma$.

\vskip.1in

If $\Sigma$ is $\pi$-nonsingular, then we are done.  

\vskip.1in

So we may assume $\Sigma$ is not $\pi$-nonsingular and
hence $\text{$\pi$-m.p.}(\Sigma) = (g_{\Sigma};s)$ \linebreak with $g_{\Sigma} >
(1,*,*,*)$.  We only have to construct ${\tilde \Sigma}$ obtained from $\Sigma$
by a succession of star subdivisions such that
$\text{$\pi$-m.p.}({\tilde
\Sigma}) <
\text{$\pi$-m.p.}(\Sigma)$.

Let $\eta$ be a maximal cone of $\Sigma$ such that $\text{$\pi$-m.p.}(\eta) =
g_{\Sigma}$ with $\sigma$ being the unique circuit contained in $\eta$.  

If $\dim \sigma \leq 2$, then we let $\gamma$ be a maximal
$\pi$-independent face of $\eta$ with \linebreak
$\text{$\pi$-mult}(\gamma) =
h_{\eta}$.  We let $\tau$ be a minimal $\pi$-singular (i.e. not
$\pi$-nonsingular) face of $\gamma$ so that we can pick a point $q \in
\roman{par}(\pi(\tau))$.

If $\dim \sigma > 2$, then we take the
star subdivision $\Sigma'$ of $\Sigma$ with respect to the negative or
positive center point of $\sigma$ so that either the case A or the case B
occurs according to Proposition 5.9. 

If the case A occurs, then $\text{$\pi$-m.p.}(\Sigma') <
\text{$\pi$-m.p.}(\Sigma)$ and we simply have to set $\Sigma^{\circ} =
\Sigma'$.

If the case B occurs, then we take the exceptional cone $\kappa'$ of $\sigma'$
with \linebreak
$\text{$\pi$-m.p.}(\kappa') = \text{$\pi$-m.p.}(\sigma)$ as
described in Proposition 5.9 and take the maximal $\pi$-independent face
$\gamma$ of
$\eta$ such that $\gamma \cap \sigma' = \gamma'$, where $\gamma'$ is a face of
$\kappa'$ satisfying the conditions (B-o), (B-i) and (B-ii) in Proposition
5.9.  Observe that by Lemma 5.10 there is a maximal cone $\eta'$ of $\Sigma'$
such that $\eta' \cap \sigma' = \kappa'$, $\text{$\pi$-m.p.}(\eta') =
g_{\Sigma'} = g_{\Sigma}$, $\gamma$ is a face of $\eta'$ as well as that of
$\eta$, $\text{$\pi$-mult}(\gamma) = h_{\eta'} = h_{\eta}$ and that
$\gamma$ is codefinite with respect to $\eta'$.

We also take $\tau$ to be a minimal $\pi$-singular (i.e. not
$\pi$-nonsingular) face of $\gamma$ so that we can pick a point $q \in
\roman{par}(\pi(\tau))$. 

\vskip.1in

Now we consider the situation where $\dim \sigma \leq 2$ and the situation where
$\dim
\sigma > 2$ with the case B together. 

\vskip.1in

Take all the circuits $\theta'$ (except for the one contained in $\kappa'$) of
$\Sigma'$ such that $\tau \subset \overline{\roman{Star}(\theta')}$.  By
Corollary 5.11 of Step 3 for each $\theta'$ we can find
$\{\overline{\roman{Star}(\theta')}\}^{\circ}$ obtained by a succession of
star subdivisions by the negative or positive center points of the circuits (of
the intermediate subdivisions) inside of
$\theta'$ such that the $\pi$-multiplicity profile does not increase, i.e.,
$$\text{$\pi$-m.p.}(\{\overline{\roman{Star}(\theta')}\}^{\circ}) \leq
\text{$\pi$-m.p.}(\overline{\roman{Star}(\theta')}),$$
and that $\tau$ is a face of $\{\overline{\roman{Star}(\theta')}\}^{\circ}$ and
$\tau$ is codefinite with respect to every cone $\nu \in
\{\overline{\roman{Star}(\theta)}\}^{\circ}$ containing $\tau$.

Note that these star subdivisions can be carried out simultaneously without
affecting each other and that hence we obtain a simplicial cobordism
$\Sigma^{\circ}$ obtained from $\Sigma$ by a successive star subdivisions such
that

(o) the $\pi$-multiplicity profile does not increase, i.e.,

$$\text{$\pi$-m.p.}(\Sigma^{\circ}) \leq \text{$\pi$-m.p.}(\Sigma),$$

(i) $\eta'$ ($\eta' = \eta$ in the case $\dim \sigma = 2$) is a maximal
cone in
$\Sigma^{\circ}$ with
$$\text{$\pi$-m.p.}(\eta') = g_{\Sigma^{\circ}} = g_{\Sigma'} = g_{\Sigma},$$

(ii) $\tau$ is contained in a maximal $\pi$-independent face $\gamma$ of
$\eta'$ with the maximum $\pi$-multiplicity $\text{$\pi$-mult}(\gamma) =
h_{\eta'}$,

(iii) $\tau$ is codefinite with respect to $\eta'$ and with respect to
all the other maximal cones containing $\tau$,

(iv) we can find a lattice point $q \in \roman{par}(\pi(\tau))$.

\vskip.1in

We only have to set ${\tilde \Sigma} = \roman{Mid}(\tau,l_q) \cdot
\Sigma^{\circ}$ to observe by Proposition 5.5 in Step 2
that 
$$\text{$\pi$-m.p.}({\tilde \Sigma}) < \text{$\pi$-m.p.}(\Sigma).$$

By the descending chain condition of the set of the $\pi$-multiplicity
profiles, this completes the process of $\pi$-desingularization.

Remark that by construction the process leaves any $\pi$-independent and
already $\pi$-nonsingular cone of $\Sigma$ unaffected.

\vskip.2in

\proclaim{Remark 5.12}\endproclaim

We discuss the comparison of our arguments with the original papers
[Morelli1,2].

\vskip.1in

\noindent $(5.12.1)$ (Definition of the negative or positive center point.)

\vskip.1in

The definition of the negative or positive center point $\roman{Ctr}_-(\sigma),
\roman{Ctr}_+(\sigma)$ as presented here and in [Morelli2] is different from the
original definition of the center point $\roman{Ctr}(\sigma,\tau)$ in
[Morelli1].  In spite of the assertions in [Morelli1],
$\roman{Ctr}(\sigma,\tau)$ is not always in
$\roman{RelInt}(\pi(\tau))$, as one can see in
some easy examples.  This causes a problem in the original argument in
[Morelli1], as the subdivision corresponding to the center point may affect
not only the cones in the closed star $\overline{\roman{Star}(\sigma)}$ but also
possibly some other cones, which we do not have any control over.  This is
the first problematic point in the argument of [Morelli1] noticed by [King2].

\vskip.1in

\noindent $(5.12.2)$ (Definition of the $\pi$-multiplicity profile.)

\vskip.1in

In [Morelli1], the $\pi$-multiplicity pofile $\text{$\pi$-m.p.}(\eta)$ of a
simplicial cone $\eta$ was defined to be
$$\text{$\pi$-m.p.}(\eta) = (\text{$\pi$-mult}(\gamma_1), \cdot\cdot\cdot,
\text{$\pi$-mult}(\gamma_l))$$
where $\gamma_1, \cdot\cdot\cdot, \gamma_l$ are the maximal $\pi$-independent
faces of $\eta$ with  
$$\text{$\pi$-mult}(\gamma_1) \geq \text{$\pi$-mult}(\gamma_2)
\geq \cdot\cdot\cdot \geq \text{$\pi$-mult}(\gamma_l)).$$

Proposition 5.5 holds with this definition, while Proposition 5.9 fails to
hold, as [King2] noticed.

(With the slightly coarser definition of the $\pi$-multiplicity profile
$$\text{$\pi$-m.p.}(\eta) = (h_{\eta},r_{\eta}),$$
Proposition 5.5 holds, while Proposition 5.9 fails to
hold in a similar way.)

\vskip.1in

In [Morelli2], the $\pi$-multiplicity pofile $\text{$\pi$-m.p.}(\eta)$ of a
simplicial cone $\eta$ was changed and defined to be
$$\text{$\pi$-m.p.}(\eta) = (h_{\eta},k_{\eta},r_{\eta}).$$

Proposition 5.9 holds with this definition, while now in turn Proposition 5.5
fails to hold.

\vskip.1in

The current and correct definition of the $\pi$-multiplicity profile, as
presented here, was suggested to us by Morelli after we discussed the
dilemma as above through e-mail.

\vskip.1in

\noindent $(5.12.3)$ (How to choose $\tau$ with $q \in \roman{par}(\pi(\tau))$
and make it codefinite.)

\vskip.1in

[Morelli1] could be read (by a naive reader like us) in such a way that it
suggests that for a maximal $\pi$-independent face $\gamma$ with the maximum
$\pi$-multiplicity
$\text{$\pi$-mult}(\gamma) = h > 1$ we could take $q \in
\roman{par}(\pi(\gamma))$, which is clearly false in the case $\dim N_{\Bbb Q}
\geq 3$.  The subdivision with respect to $q
\in \roman{par}(\pi(\gamma))$ would only affect the cones in the star
$\roman{Star}({\gamma})$ and we would only have to analyze those circuits
$\sigma$ such that $\gamma \subset \overline{\roman{Star}(\sigma)}$.  Then the
face
$\zeta =
\gamma
\cap \sigma$ has the maximum $\pi$-multiplicity $h_{\sigma}$ and only Lemma
5.8 would suffice to achieve codefiniteness after the subdivision by the
negative or positive center point. 

But in general it is only a subface $\tau \subset \gamma$ which contains a
point $q \in \roman{par}(\pi(\tau))$.  Now we have to analyze those circuits
$\sigma$ such that $\tau \subset \overline{\roman{Star}(\sigma)}$ but maybe
\linebreak
$\gamma
\not\subset \overline{\roman{Star}(\sigma)}$.  Lemma 5.8 is not sufficient any
more to achieve the codefiniteness for $\tau$.  This is another problematic
point in the argument of [Morelli1] noticed by [King2].

\vskip.1in

[Morelli2] tries to fix this problem via the use of Proposition 5.9 and what
Morelli calls the trivial subdivision of a circuit $\sigma$.

\vskip.1in

Our argument here to achieve Corollary 5.11 solves the problem by induction on
$\pi$-multiplicity profile based upon Proposition 5.9 and does not use the
trivial subdivision.

\newpage

\S 6. The Weak Factorization Theorem

\vskip.1in

In this section, we harvest the fruit ``Weak Factorization Theorem" grown upon
the tree of the results of the previous sections.

\proclaim{Proposition 6.1} We have the weak factorization of a proper
equivariant birational map between two nonsingular toric varieties $X_{\Delta}$
and
$X_{\Delta'}$ if and only if there exists a simplicial, collapsible and
$\pi$-nonsingular cobordism $\Sigma$ between the fans $\Delta$ and $\Delta'$.
\endproclaim

\demo{Proof}\enddemo Suppose we have the weak factorization of a proper
equivariant birational map between two nonsingular toric varieties $X_{\Delta}$
and
$X_{\Delta'}$.  Then the fan $\Delta'$ is obtained from $\Delta$ by a sequence
of smooth star subdivisions and smooth star assemblings (in arbitrary order). 
By Proposition 4.7 there exists a simlicial and collapsible cobordism $\Sigma$
between $\Delta$ and $\Delta'$, which is also $\pi$-nonsingular by
construction (cf. the proof of Proposition 4.7).

Conversely, suppose there exists a simplicial, collapsible and
$\pi$-nonsingular cobordism $\Sigma$ between the fans $\Delta$ and $\Delta'$. 
Write
$$\Sigma = \cup_{\sigma}\overline{\roman{Star}(\sigma)}
\cup \partial_-\Sigma$$
where the union is taken over all the circuits $\sigma$.  By the collapsibility
of
$\Sigma$, we can order the circuits $\sigma_1, \cdot\cdot\cdot, \sigma_m$ so
that each $\sigma_i$ is minimal among the circuits $\sigma_i, \sigma_{i+1},
\cdot\cdot\cdot, \sigma_m$ with respect to the partial order given by the
circuit graph of $\Sigma$.  Accordingly, we have a sequence of fans
$$\align
\Delta = \Delta_0 &= \pi(\partial_-\Sigma) =
\pi(\partial_-\{\cup_{i = 1}^k\overline{\roman{Star}(\sigma_i)}
\cup\partial_+\Sigma\})\\
\Delta_1 &= \pi(\partial_-\{\cup_{i = 2}^k\overline{\roman{Star}(\sigma_i)}
\cup \partial_+\Sigma\})\\
&\cdot\cdot\cdot \\
\Delta_j &= \pi(\partial_-\{\cup_{i = j+1}^k\overline{\roman{Star}(\sigma_i)}
\cup \partial_+\Sigma\})\\
&\cdot\cdot\cdot \\
\Delta_k &= \partial_+\Sigma = \Delta'.\\
\endalign$$

Note that the fan $\Delta_{j+1}$ is obtained from $\Delta_j$ by replacing
$\partial_-\overline{\roman{Star}(\sigma_j)}$ with
$\partial_+\overline{\roman{Star}(\sigma_j)}$, which is the bistellar operation
analyzed in \S 3 and corresponds to a smooth star subdivision followed by a
smooth star assembling.  Therefore, we conclude $X_{\Delta'}$ is obtained
from $X_{\Delta}$ by a sequence of equivariant smooth blowups and smooth
blowdowns. 

\vskip.1in

\proclaim{Theorem 6.2 (The Weak Factorization Theorem)} We have the weak
factorization for every proper and equivariant birational map between two nonsingular toric varieties $X_{\Delta}$
and
$X_{\Delta'}$, i.e., Conjecture 1.1 holds in the weak from.
\endproclaim

\demo{Proof}\enddemo Let $\Delta$ and $\Delta'$ be the corresponding
nonsingular fans in $N_{\Bbb Q}$ with the same support.  Then by Theorem 4.3
there exists a simplicial and collapsible cobordsim $\Sigma$ in $N_{\Bbb Q}^+$
between $\Delta$ and $\Delta'$.  Theorem 5.1 implies there is a simplicial fan
${\tilde \Sigma}$ obtained from $\Sigma$ by a sequence of star subdivisions
such that ${\tilde \Sigma}$ is $\pi$-nonsingular and that the process leaves
all the $\pi$-independent and $\pi$-nonsingular cones of $\Sigma$ unaffected. 
By Lemma 4.8 we see that ${\tilde \Sigma}$ is also collapsible as well as
simplicial and
$\pi$-nonsingular and that the lower face and upper face of ${\tilde \Sigma}$
are unaffected and hence isomorphic to $\Delta$ and $\Delta'$, respectively. 
Thus ${\tilde \Sigma}$ is a simplicial, collapsible and $\pi$-nonsingular
cobordism between $\Delta$ and $\Delta'$.  By Proposition 6.1, we have the
weak factorization between $X_{\Delta}$ and $X_{\Delta'}$.  This completes the
proof of Theorem 6.2.

\newpage

\S 7. The Strong Factorization Theorem.

\vskip.1in

The purpose of this section is to show the strong factorization theorem, i.e.,
a proper and equivariant birational map $X_{\Delta} \dashrightarrow
X_{\Delta'}$ between smooth toric varieties can be factored into a sequence
of smooth equivariant blowups $X_{\Delta} \leftarrow X_{\Delta''}$ followed
immediately by smooth equivariant blowdowns $X_{\Delta''}
\rightarrow X_{\Delta'}$, based upon the weak factorization theorem (of \S 6 or
[W{\l}odarczyk1]).  The main difference between the weak and strong
factorization theorems is that the former allows the sequence to consist of
blowups and blowdowns in any order for the factorization, while the latter
allows the sequence to consist only of blowups first and immediately followed by
blowdowns.  We should emphasize that this section uses only the statement of
the weak factorization theorem and hence is independent of the methods of the
previous sections and that the reader, if he wishes, can use [W{\l}odarczyk1]'s
result as the starting point for this section (though we continue to phrase the
statements in Morelli's terminology that we have been using up to
\S 6).  

\vskip.1in

Our strategy goes as follows.  We start with a simplicial, collapsible and
$\pi$-nonsingular cobordism $\Sigma$ between $\Delta$ and $\Delta'$, whose
existence is guaranteed by Theorem 6.2.  We construct a new cobordism ${\tilde
\Sigma}$ from $\Sigma$ applying an appropriate sequence of star subdivisions
such that $\partial_-\Sigma = \partial_-{\tilde \Sigma}$ is
unaffected through the process of the star subdivisions and that the
cobordism ${\tilde \Sigma}$ represents, via the bistellar operations (cf.
Theorem 3.2), a sequence consisting only of smooth star subdivisions starting
from
$\Delta = \pi(\partial_-\Sigma) = \pi(\partial_-{\tilde
\Sigma})$ and ending with $\pi(\partial_+{\tilde \Sigma})$.  Observing that
$\pi(\partial_+{\tilde
\Sigma})$ is obtained from $\pi(\partial_+\Sigma) = \Delta'$ by a sequence
consisting only of smooth star subdivisions, or equivalently
$\Delta' = \pi(\partial_+\Sigma)$ is obtained from $\pi(\partial_+{\tilde
\Sigma})$ by a sequence consisting only of smooth star assemblings, we achieve
the strong factorization
$$\Delta = \pi(\partial_-\Sigma) = \pi(\partial_-{\tilde \Sigma}) \leftarrow
\pi(\partial_+{\tilde \Sigma}) \rightarrow \pi(\partial_+\Sigma) = \Delta'.$$

\vskip.1in

First we identify the condition for the bistellar operation to consist of a
single smooth star subdivision.

\vskip.1in

$\bold{Definition\ 7.1.}$ A $\pi$-nonsingular simplicial circuit 
$$\sigma = \langle(v_1,w_1), (v_2,w_2), \cdot\cdot\cdot, (v_k,w_k)\rangle
\subset N_{\Bbb Q}
\oplus {\Bbb Q} = N_{\Bbb Q}^+$$
is called pointing up (resp. pointing down) if it has
exactly one positive (resp. negative) extremal ray, i.e., we have the linear
relation among the primitive vectors
$v_i = n(\pi(\rho_i))$ of the projections of the extremal rays $\rho_i$ for
$\sigma$ (after re-numbering)
$$\align
v_1 - v_2 - \cdot\cdot\cdot - v_k &= 0 \text{\ with\ }w_1 - w_2 -
\cdot\cdot\cdot - w_k > 0 \\
(\text{resp.\ }- v_1 + v_2 + \cdot\cdot\cdot + v_k &= 0 \text{\ with\ } - w_1 +
w_2 +
\cdot\cdot\cdot + w_k > 0).\\
\endalign$$ 

\proclaim{Lemma 7.2} Let $\Sigma$ be a simplicial and
$\pi$-nonsingular cobordism in
$N_{\Bbb Q}^+$ and $\sigma \in \Sigma$ a circuit which is
pointing up.  Let  
$$\sigma = \langle(v_1,w_1), (v_2,w_2), \cdot\cdot\cdot, (v_k,w_k)\rangle
\subset N_{\Bbb Q}
\oplus {\Bbb Q} = N_{\Bbb Q}^+$$
with the linear relation among the primitive vectors $v_i = n(\pi(\rho_i))$
of the projections of the extremal rays $\rho_i$ for $\sigma$
$$\align
v_1 - v_2 - \cdot\cdot\cdot - v_k &= 0 \text{\ with\ }w_1 - w_2 -
\cdot\cdot\cdot - w_k > 0. \\
\endalign$$ 
Then the bistellar operation going from
$\pi(\partial_-\overline{\roman{Star}(\sigma)})$ to
$\pi(\partial_+\overline{\roman{Star}(\sigma)})$ is a smooth star subdivision
with respect to the ray generated by 
$$v_1 = v_2 + \cdot\cdot\cdot + v_k.$$
If $\sigma$ is pointing down with the linear relation
$$- v_1 + v_2 + \cdot\cdot\cdot + v_k = 0 \text{\ with\ }- w_1 + w_2 +
\cdot\cdot\cdot + w_k > 0,$$
then the the bistellar operation going from
$\pi(\partial_-\overline{\roman{Star}(\sigma)})$ to
$\pi(\partial_+\overline{\roman{Star}(\sigma)})$ is a smooth star assembling,
the inverse of a smooth star subdivision going from
$\pi(\partial_+\overline{\roman{Star}(\sigma)})$ to
$\pi(\partial_-\overline{\roman{Star}(\sigma)})$ with respect to the ray
generated by 
$$v_1 = v_2 + \cdot\cdot\cdot + v_k.$$
\endproclaim 

The proof is immediate from Theorem 3.2.

\vskip.1in

\proclaim{Lemma 7.3} Let $\Sigma$ be a simplicial and $\pi$-nonsingular
cobordism.  Let 
$$\tau = \langle(v_1,w_1), \cdot\cdot\cdot, (v_l,w_l)\rangle$$
be a $\pi$-independent cone of $\Sigma$ with the $v_i = n(\pi(\rho_i))$ being
the primitive vectors of the projections of the extremal rays $\rho_i$ for
$\tau$.  Let
$\rho_{\tau}$ be the midray $\roman{Mid}(\tau,l_{r(\tau)})$, where
$r(\tau) \in N$ is the vector $r(\tau) = v_1 + \cdot\cdot\cdot + v_l$, called
the ``$\pi$-barycenter" of $\tau$.  If
$\tau$ is codefinite with respect to all the circuits $\sigma
\in \Sigma$ with $\tau \in \overline{\roman{Star}(\sigma)}$, then $\rho_{\tau}
\cdot \Sigma$ stays
$\pi$-nonsingular.
\endproclaim

\demo{Proof}\enddemo Note that though in the statement of Proposition 5.5 the
point $q$ was assumed to be taken from $\roman{par}(\pi(\tau))$, we only need
the description
$$q = \Sigma a_iv_i \text{\ with\ }0 \leq a_i \leq 1$$
(allowing the equality $a_i = 1$) to conclude that the maximum of the
$\pi$-multiplicities of the $\pi$-independent cones does not increase.  Thus we
can apply the argument in the proof of Proposition 5.5 with
$$q = r(\tau) = v_1 + \cdot\cdot\cdot + v_l$$
to conclude that the maximum of the $\pi$-multiplicities of the
$\pi$-independent cones does not increase and in particular
$\rho_{\tau} \cdot \Sigma = \roman{Mid}(\tau,l_{\tau})
\cdot \Sigma$ stays $\pi$-nonsingular.

\vskip.1in

$\bold{Definition\ 7.4.}$ Let $I$ be a subset, consiting only of
$\pi$-independent cones, of a simplicial cobordism
$\Sigma$.  Assume
$I$ is join closed, i.e., 
$$\tau, \tau' \in I \Longrightarrow \tau +
\tau'
\in I \hskip.1in (\text{provided\ } \tau + \tau' \in \Sigma).$$  
We denote
$$I \cdot \Sigma = \rho_{\tau_n} \cdot\cdot\cdot \rho_{\tau_1} \cdot \Sigma$$
where $\rho_{\tau_i}$ is the midray $\roman{Mid}(\tau_i,l_{r(\tau_i)})$ with
$r(\tau_i)$ being the $\pi$-barycenter of $\tau_i$, as
described in Lemma 7.3, and where the
$\tau_i$ are cones in
$I$ so ordered that \linebreak
$\dim \tau_i \geq \dim \tau_{i+1} \hskip.1in
\text{for\ all\ } i$.  (Observe that, as $I$ is join closed, $I \cdot \Sigma$ is
independent of the choice of the order and is well-defined.)

\vskip.1in

The following simple observation of Morelli is the basis of our method in this
section.

\proclaim{Lemma 7.5} Let $\sigma$ be a circuit in a simplicial
and $\pi$-nonsingular cobordism
$\Sigma$.  Let
$$\sigma = \langle(v_1,w_1), \cdot\cdot\cdot, (v_m,w_m),(v_{m+1},w_{m+1}),
\cdot\cdot\cdot, (v_k,w_k)\rangle,$$
where $v_1, \cdot\cdot\cdot, v_m, v_{m+1}, \cdot\cdot\cdot, v_k$ are
the primitive vectors in $N$ of the projections of the extremal rays for
$\sigma$, having the unique linear relation
$$v_1 + \cdot\cdot\cdot + v_m - v_{m+1} - \cdot\cdot\cdot - v_k = 0
\text{\ with\ \ }w_1 +
\cdot\cdot\cdot + w_m - w_{m+1} - \cdot\cdot\cdot - w_k > 0.$$
Let
$$\sigma_+ = \langle(v_1,w_1), \cdot\cdot\cdot, (v_m,w_m)\rangle
\text{\ and\ }\sigma_- =
\langle(v_{m+1},w_{m+1}), \cdot\cdot\cdot, (v_k,w_k)\rangle.$$

$(7.5.1)$ The fan $\rho_{\sigma_+} \cdot \overline{\roman{Star}(\sigma)}$, where
$\rho_{\sigma_+}$ is the midray $\roman{Mid}(\sigma_+,l_{r(\sigma_+)})$ with
$r(\sigma_+)$ being the
$\pi$-barycenter of $\sigma_+$, is
$\pi$-nonsingular and the closed star of a
$\pi$-nonsingular pointing up circuit $\sigma'$.

$(7.5.2)$ If $\sigma$ is pointing up and $I$ is a join closed subset of
$\sigma_-$, then $I \cdot
\overline{\roman{Star}(\sigma)}$ is
$\pi$-nonsingular and the closed star of a
$\pi$-nonsingular pointing up circuit.

\endproclaim

\demo{Proof}\enddemo (7.5.1) First note that, since $\sigma$ is $\pi$-strongly
convex and hence does not contain a nonzero vector $0 \neq (0,w) \in N_{\Bbb
Q}^+ = N_{\Bbb Q} \oplus {\Bbb Q}$, it is impossible to have all the
coefficients in the linear relation to be
$+1$ or all to be
$-1$. 

Let $\eta \in \roman{Star}(\sigma)$ be a simplicial cone of the form
$$\eta = \langle(u_1,w_1'), \cdot\cdot\cdot, (u_l,w_l'), (v_1,w_1),
\cdot\cdot\cdot, (v_k,w_k)\rangle.$$
Then the maximal cones of $\rho_{\sigma_+} \cdot \eta$ are of the form
$$\align
\langle(u_1,w_1'), \cdot\cdot\cdot, (u_l,w_l'), (v_1,w_1), \cdot\cdot\cdot,
& \overset{\vee}\to{(v_i,w_i)}, \cdot\cdot\cdot, (v_m,w_m), \\
 \hskip.1in & 1 \leq i \leq m\\
&(v_{m+1},w_{m+1}),
\cdot\cdot\cdot, (v_k,w_k), (r(\sigma_+),\Sigma_{i = 1}^m w_i)\rangle \\
\endalign$$
omitting one of $(v_i,w_i)$, $1 \leq i \leq m$, from the generators of
$\sigma_+$.  Therefore, 
$$\sigma' = \langle (r(\sigma_+),\Sigma_{i = 1}^m w_i), (v_{m+1},w_{m+1}),
\cdot\cdot\cdot, (v_k,w_k)\rangle$$
is the unique circuit in 
$\rho_{\sigma_+} \cdot
\overline{\roman{Star}(\sigma)}$ and
$$\rho_{\sigma_+} \cdot \overline{\roman{Star}(\sigma)} =
\overline{\roman{Star}(\sigma')}.$$ As $\rho_{\sigma_+}$ is generated by the
vector
$(r(\sigma_+),\Sigma_{i = 1}^mw_i) = (\Sigma_{i = 1}^mv_i,\Sigma_{i = 1}^mw_i)$,
the unique linear relation for $\sigma'$ is 
$$\align
&n(\pi(\rho_{\sigma_+})) - v_{m+1} - \cdot\cdot\cdot - v_k = 0 \text{\ where\
}n(\pi(\rho_{\sigma_+})) = \Sigma_{i = 1}^mv_i\\
&\hskip2in \text{\ with\ }(\Sigma_{i = 1}^mw_i) - w_{m+1} - \cdot\cdot\cdot -
w_k > 0.\\
\endalign$$
Therefore, the circuit $\sigma'$ is pointing up.  We note
that $\pi$-nonsingularity is preserved as $\sigma_+$ is obviously codefinite
with respect to the circuit
$\sigma$.

\vskip.1in

(7.5.2) We use the same notation as in (7.5.1) with $\sigma_+ =
\langle(v_1,w_1)\rangle$ being the only positive extremal ray of the pointing up
circuit $\sigma$.  Let
$\zeta$ be the maximal cone in $I$.  Then the maximal cones $\eta'$ of
$\rho_{\zeta} \cdot
\eta$, where $\rho_{\zeta}$ is the midray
$\roman{Mid}(\zeta,l_{r(\zeta)})$, are of the form
$$\align
\langle(u_1,w_1'), \cdot\cdot\cdot, (u_l,w_l'), &(v_1,w_1), \\
(v_{1+1},w_{1+1}), \cdot\cdot\cdot, &\overset{\vee}\to{(v_j,w_j)},
\cdot\cdot\cdot, (v_k,w_k), (r(\zeta),\Sigma_{(v_i,w_i) \in \zeta}w_i)\rangle\\ 
&(v_j,w_j) \in \zeta. \\
\endalign$$
Therefore,
$$\sigma_{\zeta} = \langle(v_1,w_1), \text{all\ the\ } (v_i,w_i) \hskip.1in
\not\in
\zeta, (r(\zeta),\Sigma_{(v_i,w_i) \in \zeta}w_i)\rangle$$ 
is the unique circuit in $\rho_{\zeta} \cdot \overline{\roman{Star}(\sigma)}$,
which is pointing up with the unique linear relation
$$\align
&v_1 - \Sigma_{(v_i,w_i) \not\in \zeta}v_i - n(\pi(\rho_{\zeta})) = 0 \text{\
where\ }n(\pi(\rho_{\zeta})) = \Sigma_{(v_i,w_i) \in \zeta}v_i\\
&\hskip2in \text{\ with\ }w_1 - \Sigma_{(v_i,w_i) \not\in \zeta}w_i -
(\Sigma_{(v_i,w_i)
\in \zeta}w_i) > 0.\\
\endalign$$
With $\eta \in \roman{Star}(\sigma)$ being arbitrary, we also have
$$\rho_{\zeta} \cdot \overline{\roman{Star}(\sigma)} =
\overline{\roman{Star}(\sigma_{\zeta})}.$$ 
Moreover, every cone in the
complement
$I'$ of $\zeta$ in $I$ (i.e., $I'$ consits of the proper subfaces of $\zeta$) is
disjoint from
$\sigma_{\zeta}$.  Therefore,
$\sigma_{\zeta}$ is still the unique circuit, which is pointing up, in
$$I \cdot \overline{\roman{Star}(\sigma)} = I' \cdot \rho_{\zeta} \cdot
\overline{\roman{Star}(\sigma)}$$ 
and
$$I \cdot \overline{\roman{Star}(\sigma)} =
\overline{\roman{Star}(\sigma_{\zeta})}.$$

This completes the proof of Lemma 7.5.

\vskip.1in

The following is an easy consequence of Lemma 7.5.

\proclaim{Lemma 7.6} Let $\Sigma$ be a simplicial, collapsible and
$\pi$-nonsingular cobordism whose circuits are all pointing up and let $I
\subset
\partial_-\Sigma$ be a join closed subset.  Assume the condition $(\star)$:

$$(\star) \hskip.1in I \cap \overline{\roman{Star}(\sigma)} \subset \{\tau \in
\Sigma;\tau \subset \sigma_-\} = \partial_-\sigma
\hskip.1in
\text{\ for\ any\ circuit\ } \sigma \in \Sigma.$$  

Then $\Sigma' = I \cdot \Sigma$ is again a simplicial, collapsible and
$\pi$-nonsingular cobordism containing only pointing up circuits. 

\endproclaim

\demo{Proof}\enddemo By Lemma 4.8 and Lemma 7.3 the cobordism $\Sigma'$ is
again simplicial, collapsible and $\pi$-nonsingular.  We only have to check that
$I \cdot
\overline{\roman{Star}(\sigma)} = (I \cap \overline{\roman{Star}(\sigma)}) \cdot
\overline{\roman{Star}(\sigma)}$ contains only pointing up
circuits for any circuit
$\sigma \in \Sigma$, which follows immediately from the condition
$(\star)$ and (7.5.2) in Lemma 7.5.

\vskip.1in

\proclaim{Remark 7.7}\endproclaim

Lemma 7.6 is a modification of Lemma 9.7 in [Morelli1] (together with the notion
of ``neatly founded"), which unfortunately has a counter-example as below. 
We observe that the notion of ``neatly founded" is used only in the form of the
condition $(\star)$ in the argument of [Morelli1] and we carry out our argument
here all through with the condition $(\star)$ instead of the notion of ``neatly
founded".

Below we recall the definition of ``neatly founded" and Lemma 9.7 in
[Morelli1] and then present a counter-example.

\vskip.1in 

[Morelli1] defines that $\Sigma$ is ``neatly founded" if for each down
definite face $\tau \in \Sigma$ (A face $\tau \in \Sigma$ is down definite
if $\tau \in \partial_-\Sigma$ but $\tau \not\in \partial_+\Sigma$.), there is a
circuit
$\sigma
\in
\Sigma$ such that
$\tau =
\sigma_-$.

\proclaim{Lemma 9.7 in [Morelli1]} Let $\Sigma$ be a neatly founded,
simplicail, collapsible and $\pi$-nonsingular cobordism whose circuits are all
pointing up, and let
$I
\subset
\partial_-\Sigma$ be join closed.  Then $\Sigma' = I \cdot \Sigma$ is
again a simplicial, collapsible and $\pi$-nonsingular cobordism containing only
pointing up circuits.
\endproclaim

A counter-example to Lemma 9.7 in [Morelli1]:

\vskip.1in

We take
$$\align
\rho_1 &= (v_1, 0) \\
\rho_2 &= (v_2, 0) \\
\rho_3 &= (v_3, 0) \\
\rho_4 &= (v_1 + v_2 + v_3, 1) \\
\rho_5 &= (v_1 + v_2 + 2v_3, 2) \\
\endalign$$ 
in $N_{\Bbb Q}^+ = (N \oplus {\Bbb Z}) \otimes {\Bbb Q} = N_{\Bbb Q} \oplus
{\Bbb Q}$ with
$\dim N_{\Bbb Q} = 3$ where $v_1, v_2, v_3$ form a ${\Bbb Z}$-basis for $N$. 
We set
$\Sigma$ to be
$$\Sigma = \left\{\aligned
&\langle\rho_1, \rho_2, \rho_3, \rho_4\rangle \text{\ and\ its\ faces},\\
&\langle\rho_2, \rho_3, \rho_4, \rho_5\rangle \text{\ and\ its\ faces},\\
&\langle\rho_1, \rho_3, \rho_4, \rho_5\rangle \text{\ and\ its\ faces}\\
\endaligned
\right\}.$$

The fan $\Sigma$ is by construction a simplicial, collapsible and
$\pi$-nonsingular cobordism between $\Delta = \partial_-\Sigma$ and $\Delta' =
\partial_+\Sigma$.

The cobordism $\Sigma$ is neatly founded as $\langle\rho_1, \rho_2,
\rho_3\rangle$ is the only down definite face and there is a circuit
$\langle\rho_1,
\rho_2, \rho_3,
\rho_4\rangle$ such that
$$\langle\rho_1, \rho_2, \rho_3\rangle = \langle\rho_1, \rho_2, \rho_3,
\rho_4\rangle_-.$$

All circuits $\langle\rho_1, \rho_2, \rho_3, \rho_4\rangle$ and $\langle\rho_3,
\rho_4,
\rho_5\rangle$ are pointing up.

Take
$$I = \left\{\langle\rho_2, \rho_3\rangle \text{\ and\ its\ faces} \right\}.$$

Now $\Sigma$ and $I$ satisfy all the conditions of Lemma 9.7.  On the other
hand, $\Sigma' = I \cdot \Sigma$ contains a circuit
$$\langle\rho_2, M, \rho_4, \rho_5\rangle \text{\ where\ } M = (v_2 + v_3,
0),$$ 
which is NOT pointing up!

\vskip.1in

We resume our proof of the implication the ``weak" factorization
$\Longrightarrow$ the ``strong" factorization. 

\vskip.1in

\proclaim{Proposition 7.8} Let $\Sigma$ be a simplicial, collapsible and
$\pi$-nonsingular cobordism containing only pointing up circuits.  Then there is
a simplicial, collapsible and $\pi$-nonsingular cobordism $\Sigma'$ such that

$(7.8.1)$ $\Sigma'$ contains only pointing up circuits,

$(7.8.2)$ $\Sigma'$ satisfies the condition $(\star)$ for any join closed subset
$I
\subset
\partial_-\Sigma'$,

$(7.8.3)$ $\Sigma'$ is obtained from $\Sigma$ by a sequence of star
subdivisions, none of which involve $\partial_-\Sigma$, of the
$\pi$-independent faces which are codefinite with respect to all the circuits.

\endproclaim

\demo{Proof}\enddemo Express the collapsible $\Sigma$ as
$$\Sigma = \overline{\roman{Star}(\sigma_m)} \circ
\overline{\roman{Star}(\sigma_{m-1})}
\circ
\cdot\cdot\cdot \overline{\roman{Star}(\sigma_1)} \circ \partial_+\Sigma$$
for the circuits $\sigma_m, \sigma_{m-1}, \cdot\cdot\cdot, \sigma_1 \in \Sigma$
so that $\sigma_i$ is minimal among $\sigma_i, \sigma_{i-1}, \cdot\cdot\cdot,
\sigma_1$ according to the partial order given by the circuit graph.  We prove
the lemma by induction on
$m$.

\vskip.1in

Case $m = 1$: This case is the building block of the construction in the
induction step and we state it in the form of a lemma as below.

\proclaim{Lemma 7.9} Let $\Sigma$ be a simplicial, collapsible and
$\pi$-nonsingular cobordism containing only pointing up circuits.  Let
$\overline{\roman{Star}(\sigma)}$ be the closed star of a circuit
$\sigma \in \Sigma$.  Let
$$J = \{\langle\sigma_+,\nu\rangle; \nu \in \roman{link}_{\Sigma}(\sigma)\}.$$
Then

$(7.9.1)$ $J \cdot \overline{\roman{Star}(\sigma)}$ contains only pointing up
circuits,

$(7.9.2)$ $J \cdot \overline{\roman{Star}(\sigma)}$ satisfies the condition
$(\star)$ for any join closed subset \linebreak
$I \subset \partial_-\{J \cdot
\overline{\roman{Star}(\sigma)}\}$, and

$(7.9.3)$ $J \cdot \overline{\roman{Star}(\sigma)}$ is obtained from
$\overline{\roman{Star}(\sigma)}$ by a sequence of star subdivisions, none of
which involve
$\partial_-\overline{\roman{Star}(\sigma)}$, of the $\pi$-independent faces
which are codefinite with respect to all the circuits.

\endproclaim

\demo{Proof}\enddemo Let
$$\sigma = \langle(v_1,w_1), (v_2,w_2), \cdot\cdot\cdot, (v_k,w_k)\rangle$$
where $v_1, v_2, \cdot\cdot\cdot, v_k$ are primitive vectors in $N$
satisfying the unique linear relation
$$v_1 - v_2 - \cdot\cdot\cdot - v_k = 0 \text{\ with\ }w_1 - w_2 -
\cdot\cdot\cdot - w_k > 0.$$
Let $\eta \in \roman{Star}(\sigma)$ be a simplicial cone of the form
$$\eta = \langle(u_1,w_1'), \cdot\cdot\cdot, (u_l,w_l'), (v_1,w_1),
\cdot\cdot\cdot, (v_k,w_k)\rangle.$$
Then the circuits of $J \cdot \eta = \{J \cap \eta\} \cdot \eta$ are the cones
of the form
$$\sigma_{\nu} = \langle(r(\langle\sigma_+,\nu\rangle),w_1 +
\Sigma_{(u_j,w_j') \in \nu}w_j'), (v_2,w_2),
\cdot\cdot\cdot, (v_k,w_k), \text{all\ the\ }(u_j,w_j') \in \nu\rangle
\text{\ for\ }\nu
\in \roman{link}_{\eta}(\sigma)$$
(including $\sigma = \sigma_{\emptyset} = \langle(r(\sigma_+),w_1) = (v_1,w_1),
(v_2,w_2),
\cdot\cdot\cdot, (v_k,w_k)\rangle$ for
$\nu =
\emptyset$) satisfying the unique linear relation
$$\align
(v_1 + \Sigma_{(u_j,w_j') \in \nu}u_j) - v_2 - \cdot\cdot\cdot - v_k -
\Sigma_{(u_j,w_j') \in \nu}u_j &= 0 \\
\text{\ with\ }(w_1 + \Sigma_{(u_j,w_j') \in \nu}w_j') - w_2 - \cdot\cdot\cdot -
w_k -
\Sigma_{(u_j,w_j') \in \nu}w_j' &> 0. \\
\endalign$$
Thus $J \cdot \eta$ contains only pointing up circuits.  Since $\eta \in
\roman{Star}(\sigma)$ is arbitrary, we conclude $J \cdot
\overline{\roman{Star}(\sigma)}$ contains only pointing up circuits, proving
(7.9.1).

We also observe that the
maximal cones of $\overline{\roman{Star}(\sigma_{\nu})}$ are of the form
$$\align
\langle\sigma_{\nu}, &\roman{Mid}(\langle v_1,\nu,(u_{p(1)},w_{p(1)}'),
\cdot\cdot\cdot, (u_{p(s)},w_{p(s)}')\rangle,l_{r(\roman{Mid}(\langle
v_1,\nu,(u_{p(1)},w_{p(1)}'),
\cdot\cdot\cdot, (u_{p(s)},w_{p(s)}')\rangle)}),\\
&\hskip2in s = 1,
\cdot\cdot\cdot, l' = l - \#\{(u_j,w_j') \in \nu\}\rangle\\
\endalign$$
where
$$(u_{p(1)},w_{p(1)}'), (u_{p(2)},w_{p(2)}'), \cdot\cdot\cdot,
(u_{p(l')},w_{p(l')}')$$ are the $(u_j,w_j')$'s NOT belonging to $\nu$, ordered
in the specified way by a permutation $p$.  Therefore, any cone in
the lower face $\partial_-\overline{\roman{Star}(\sigma_{\nu})}$, if not
included in
$\sigma_{\nu}$, is also in the upper face but not in the lower face of the
closed star of some other circuit of $J \cdot
\overline{\roman{Star}(\sigma)}$.  Therefore, we conclude that for
any join closed subset \linebreak
$I \subset \partial_-\{J \cdot
\overline{\roman{Star}(\sigma)}\}$ we have
$$\align
I \cap \overline{\roman{Star}(\sigma_{\nu})} &= \partial_-\{J \cdot
\overline{\roman{Star}(\sigma)}\} \cap \{\tau \in J
\cdot
\overline{\roman{Star}(\sigma)};\tau \subset \sigma_{\nu}\} \\
&\subset
\hskip.1in \{\tau \in J \cdot
\overline{\roman{Star}(\sigma)};\tau \subset (\sigma_{\nu})_-\} =
\partial_-\sigma_{\nu}.\\
\endalign$$    
Since $\eta \in \roman{Star}(\sigma)$ is
arbitrary, this proves (7.9.2).

The condition (7.9.3) is obvious from the construction.

This completes the proof of Lemma 7.9.

\vskip.1in

We go back to the proof of Proposition 7.8 resuming the induction.

\vskip.1in

Suppose $m > 1$.  Set
$$\Sigma_{m-1} = \overline{\roman{Star}(\sigma_{m-1})} \circ \cdot\cdot\cdot
\circ
\overline{\roman{Star}(\sigma_1)} \circ \partial_+\Sigma$$
and apply the induction hypothesis to $\Sigma_{m-1}$ to obtain $\Sigma_{m-1}'$
satisfying the conditions (7.8.1), (7.8.2) and (7.8.3).  Then
$\overline{\roman{Star}(\sigma_m)}
\circ
\Sigma_{m-1}'$ is the result of a sequence of star subdivisions, none of which
involve $\partial_-\Sigma$, of the $\pi$-indepenednt faces which are
codefinite with respect to all the circuits.  Let
$$J = \{\langle({\sigma_m})_+,\nu\rangle;\nu \in
\roman{link}_{\Sigma}(\sigma_m)\}.$$ We show that $J \cdot
(\overline{\roman{Star}(\sigma_m)} \circ \Sigma_{m-1}')$ satisfies the
conditions (7.8.1), (7.8.2) and (7.8.3).

Since $\Sigma_{m-1}'$ satisfies the condition (7.8.1) and $J \subset
\partial_-\Sigma_{m-1}'$ is join closed, the condition $(\star)$ for $J$ with
Lemma 7.6 implies that $J
\cdot \Sigma_{m-1}'$ is a simplicial, collapsible and
$\pi$-nonsingular cobordism containing only pointing up circuits.  Lemma 7.9
implies that $J \cdot
\overline{\roman{Star}(\sigma_m)}$ is also a
simplicial, collapsible and $\pi$-nonsingular cobordism containing only pointing
up circuits.  Therefore,
$$\Sigma' = J \cdot (\overline{\roman{Star}(\sigma_m)} \circ \Sigma_{m-1}') = (J
\cdot
\overline{\roman{Star}(\sigma_m)}) \circ (J \cdot \Sigma_{m-1}')$$
is a simplicial, collapsible and $\pi$-nonsingular cobordism satisfying
the condition (7.8.1).  

Observe that
$$\partial_-\Sigma' = \partial_-\overline{\roman{Star}(\sigma_m)} \cup
(\partial_-\Sigma'_{m-1} - \roman{RelInt}(J)).$$
Thus by construction we have the condition (7.8.2).

Let $I$ be any join closed subset of
$\partial_-\Sigma'$.  Let $\sigma' \in
\Sigma'$ be a circuit.  If $\sigma' \in J \cdot
\overline{\roman{Star}(\sigma_m)}$, then by Lemma 7.9 we have
$$I \cap \overline{\roman{Star}(\sigma')} = (I \cap J \cdot
\overline{\roman{Star}(\sigma_m)}) \cap \overline{\roman{Star}(\sigma')} \subset
\partial_-\sigma'.$$   
If $\sigma' \in J \cdot \Sigma_{m-1}'$ and $\sigma' \not\in \Sigma_{m-1}'$, then
there exists a circuit
$\sigma = \langle(v_1,w_1), \cdot\cdot\cdot, (v_k,w_k)> \in \Sigma_{m-1}'$ such
that
$$\sigma' = \sigma_{\zeta} = <(v_1,w_1), \text{all\ the\ }(v_i,w_i) \hskip.1in i
\not\in \zeta, (r(\zeta),\Sigma_{(v_i,w_i) \in \zeta}w_i)\rangle$$
where $\zeta$ is the maximal cone in $J \cap \{\tau \in
\Sigma'_{m-1};\tau \subset \sigma\}$, using the same notation as in Lemma 7.5. 
Observe that for any maximal cone $\eta'' \in
\overline{\roman{Star}(\sigma')}$ if a face $\tau \subset \eta''$ contains a
new ray used for the subdividing operation ``$J \cdot$" as one of the generators
then
$\tau
\not\in I$.  Therefore, by looking at the description of $\eta'$ in Lemma 7.5
and $\eta''$ obtained from $\eta'$ by the star subdivision of some faces of
$\zeta$, we conclude 
$$I \cap \{\tau \i \Sigma';\tau \subset \eta''\} = (I \cap \{\tau \in
\Sigma';\tau \subset \sigma'\})
\cap
\partial_-\sigma'
\subset
\partial_-\sigma'.$$
If $\sigma' \in J \cdot \Sigma_{m-1}'$ and also $\sigma' \in \Sigma_{m-1}'$,
then the condition $(\star)$ for $\Sigma_{m-1}'$ implies
$$I \cap \overline{\roman{Star}(\sigma')} \subset \partial_-\sigma'.$$

Thus we have the condition $(\star)$ for $\Sigma'$ proving the condition
(7.8.2).

This completes the proof of Proposition 7.8.

\vskip.1in

\proclaim{Theorem 7.10} Any simplicial, collapsible and
$\pi$-nonsingular cobordism
$\Sigma$ between
$\Delta$ and
$\Delta'$ can be made into a simplicial, collapsible and
$\pi$-nonsingular cobordism
$\Sigma'$ between
$\Delta$ and
$\Delta''$ by a sequence of star subdivisions such that $\Sigma'$ contains
only pointing up circuits and that $\Delta''$ is obtained from $\Delta'$ by a
sequence of smooth star subdivisions.
\endproclaim

\demo{Proof}\enddemo Express
$$\Sigma = \overline{\roman{Star}(\sigma_m)} \circ
\overline{\roman{Star}(\sigma_{m-1})}
\circ
\cdot\cdot\cdot \circ \overline{\roman{Star}(\sigma_1)} \circ
\partial_+\Sigma$$ 
for the circuits $\sigma_m, \sigma_{m-1}, \cdot\cdot\cdot, \sigma_1 \in
\Sigma$ so that $\sigma_i$ is minimal among $\sigma_i, \sigma_{i-1},
\cdot\cdot\cdot, \sigma_1$ according to the partial order given by the
circuit graph. 

Define a sequence of cobordisms ${\tilde \Sigma}_k, {\tilde
\Sigma}_k'$ inductively as follows: Let
$$\align
{\tilde \Sigma}_0 &= {\tilde \Sigma}_0' = \partial_+\Sigma \\
{\tilde \Sigma}_k &= \rho_{\sigma_k^+} \cdot (\overline{\roman{Star}(\sigma_k)}
\circ {\tilde \Sigma}_{k-1}')\\
\endalign$$ 
where ${\tilde \Sigma}_{k-1}'$ for $k \geq 2$ is obtained from ${\tilde
\Sigma}_{k-1}$ by the procedure described in Proposition 7.8 to satisfy the
conditions (7.8.1), (7.8.2) and (7.8.3).  We remark that 
$$\partial_-{\tilde \Sigma}_k = \partial_-{\tilde \Sigma}_k' =
\partial_-(\overline{\roman{Star}(\sigma_k)} \circ
\overline{\roman{Star}(\sigma_{k-1})}
\circ \cdot\cdot\cdot \overline{\roman{Star}(\sigma_1)} \circ
\partial_+\Sigma.$$   Note then that inductively  by Lemma 7.5, Lemma 7.6 and
Proposition 7.8 ${\tilde \Sigma}_k$ is a simplicial, collapsible and
$\pi$-nonsingular cobordism containing only pointing up circuits.  Finally
${\tilde \Sigma} = {\tilde
\Sigma}_m$ is a
simplicial, collapsible and $\pi$-nonsingular cobordism containing only pointing
up circuits between
$\Delta = \pi(\partial_-\Sigma) =
\pi(\partial_-{\tilde \Sigma})$ and
$\Delta'' = \pi(\partial_+{\tilde \Sigma})$, which is obtained from $\Delta'$ by
a sequence of smooth star subdivisions.

This completes the proof of Theorem 7.10.

\vskip.1in  

\proclaim{Corollary 7.11 (The Strong Factorization Theorem)} We have the
strong factorization for every proper and equivariant birational map between
two nonsingular toric varieties $X_{\Delta}$ and $X_{\Delta'}$, i.e.,
Conjecture 1.1 holds in the strong form.  In particular, if both $X_{\Delta}$
and $X_{\Delta'}$ are projective, then the factorization can be chosen so that
all the intermediate toric varieties are also projective.
\endproclaim

\demo{Proof}\enddemo Let $\Delta$ and $\Delta'$ be the corresponding two
nonsingular fans in $N_{\Bbb Q}$ with the same support.  Then by Proposition
6.1 and Theorem 6.2 there exists a simplicial, collapsible and
$\pi$-nonsingular cobordism $\Sigma$ between $\Delta$ and $\Delta'$.  By
Theorem 7.10 we can make $\Sigma$ into a simplicial, collapsible and
$\pi$-nonsingular cobordism ${\tilde \Sigma}$ with only pointing up cuircuits
between
$\Delta$ and
a fan $\Delta''$ such that $\Delta''$ is obtained from $\Delta'$ by a sequence
of smooth star subdivisions.  By Lemma 7.2 $\Delta'' =
\pi(\partial_+{\tilde \Sigma})$ is also obtained from
$\pi(\partial_-{\tilde \Sigma}) = \Delta$ by a sequence of smooth star
subdivisions.  Thus we have the factorization
$$\Delta = \pi(\partial_-\Sigma) = \pi(\partial_-{\tilde \Sigma})
\leftarrow
\pi(\partial_+{\tilde \Sigma}) \rightarrow \pi(\partial_+\Sigma)
= \Delta',$$ which corresponds to the strong factorization
$$X_{\Delta} \leftarrow X_{\Delta''} \rightarrow X_{\Delta'}.$$

\newpage

\S 8. The Toroidal Case

\vskip.1in

The purpose of this section is to generalize the main theorem of the previous
sections, namely the strong factorization of a proper and equivariant
birational map between two nonsingular toric varieties, to the one in the
toroidal case.

First we recall several definitions about the toroidal embeddings
(cf.[Kempf-Knudsen-Mumford-SaintDonat]) and the notion of a ``toroidal" morphism
as in [Abramovich-Karu].

\vskip.1in

$\bold{Definition 8.1\ (Toroidal\ Embeddings).}$ Given a normal variety $X$ and
an open subset
$U_X
\subset X$, the embedding $U_X \subset X$ is called toroidal if for every
closed point $x \in X$ there exist an affine toric variety $X_{\sigma}$, a
closed point $s \in X_{\sigma}$ and an isomorphism of complete local algebras
$${\hat {\Cal O}}_{X,x} \cong {\hat {\Cal O}}_{X_{{\sigma},s}}$$
so that the ideal in ${\hat {\Cal O}}_{X,x}$ generated by the ideal of $X -
U_X$ corresponds under this isomorphism to the ideal in ${\hat {\Cal
O}}_{X_{{\sigma},s}}$ generated by the ideal of $X_{\sigma} - T$, where $T$ is
the torus.  The affine toric variety$X_{\sigma}$ is called a local model of
$X$ at $x$.

We will always assume that the irreducible components of $\bigcup_{i \in I}E_i =
X - U_X$ are normal, i.e., $U_X \subset X$ is a toroidal embedding without
self-intersection.  (In fact, in most of the cases $X$ is nonsingular and
$\bigcup_{i \in I}E_i \subset X$ is a divisor with normal crossings whose
irreducible components are all nonsingular.)

The irreducible components of $\bigcap_{i \in J}E_i$ for $J \subset I$,
together with $U_X$, define a stratification of $X$.  (These components and
$X$ are the closures of the strata.  The closures of the strata formally
correspond to the closures of the orbits in local models.)

Let $S$ be a stratum in $X$, which is by definition an open set in an
irreducible component of $\bigcap_{i \in J}E_i$ for some $J \subset I$.  The
star
$\roman{Star}(S)$ is the union of those strata containing $S$ in their closure
(each of them corresponds to some $K \subset J \subset I$).  To the stratum $S$
one associates the following data:

$M^S$: the group of Cartier divisors in $\roman{Star}(S)$ supported in
$\roman{Star}(S) - U_X$

$N^S := \roman{Hom}(M^S,{\Bbb Z})$

$M^S_+ \subset M^S$: effective Cartier divisors

$\sigma^S \subset N_{\Bbb R}^S$: the dual of $M_S^+$.

\vskip.1in

If $(X_{\sigma},s)$ is a local model at $x \in X$ in the stratum $S$, then
$$M^S \cong M_{\sigma}/\sigma^{\perp}, N^S \cong N_{\sigma} \cap
\roman{span}(\sigma) \text{\ and\ } \sigma^S \cong \sigma.$$

The cones glue together to form a conical complex
$$\Delta_X = (|\Delta_X|,\{\sigma^S\},\{N^S\}),$$
where $|\Delta_X| = \bigcup_S\sigma^S$ is the support of $\Delta_X$ and the
lattices $N^S$ form an integral structure on $\Delta_X$ with $\sigma^S
\hookrightarrow N^S_{\Bbb R}$.

\newpage

$\bold{Definition 8.2\ (Toroidal\ Morphisms).}$ A dominant morphism
$$f:(U_X \subset X) \rightarrow (U_Y \subset Y)$$
of toroidal embeddings is called toroidal if for every closed point $x \in X$
there exist local models $(X_{\sigma},s)$ at $x$ and $X_{{\tau},t}$ at $y =
f(x)$ and a toric morphism $g:X_{\sigma} \rightarrow X_{\tau}$ such that the
following diagram commutes
$$\CD
{\hat {\Cal O}_{X,x}} \hskip.1in @. \cong @. \hskip.1in {\hat {\Cal
O}_{X_{\sigma},s}}\\ @AA{{\hat f}^*}A @. @AA{{\hat g}^*}A \\
{\hat {\Cal O}_{Y,y}} \hskip.1in @. \cong @. \hskip.1in {\hat {\Cal
O}_{X_{\tau},t}}.\\
\endCD$$

Now we can state our main result of this section.

\proclaim{Theorem 8.3} Let
$$f:(U_X \subset X) \rightarrow (U_Y \subset Y)$$
be a proper birational and toroidal morphism between toroidal embeddings where
$X$ and $Y$ are nonsingular and $\bigcup_{i \in I}E_i = X - U_X$ and
$\bigcup_{j
\in J}F_j$ are divisors with normal crossings whose irreducible components are
all nonsingular.
Then there exist a toroidal embedding $(U_V,V)$ and sequences of blowups, with
centers being smooth closed strata, which factor $f$
$$(U_X,X) \leftarrow (U_V,V) \rightarrow (U_Y,Y).$$
\endproclaim

\proclaim{Lemma 8.4} Let
$$f:(U_X \subset X) \rightarrow (U_Y \subset Y)$$
be a toroidal morphism between two toroidal embeddings.

$(8.4.1)$ $f$ induces a morphism $f_{\Delta}:\Delta_X \rightarrow \Delta_Y$ of
complexes such that each $\sigma^S \in \Delta_X$ maps to some $\sigma^{S'}
\in \Delta_Y$ linearly $f_{\Delta}:\sigma^S \hookrightarrow \sigma^{S'}$ with
the map of lattices of the integral structures $N_{\sigma^S} \rightarrow
N_{\sigma^{S'}}$.

$(8.4.2)$ If $f$ is proper and birational, then each $\sigma^S \in \Delta_X$
maps injectively into some $\sigma^{S'}
\in \Delta_Y$ linearly $f_{\Delta}:\sigma^S \hookrightarrow \sigma^{S'}$ and
the lattice $N_{\sigma^S}$ is a saturated sublattice of
$N_{\sigma^{S'}}$.  In short, $\Delta_X$ is a refinement of $\Delta_Y$ with
$|\Delta_X| = |\Delta_Y|$ preserving the integral structure.  Moreover, once
we fix the toroidal embedding $(U_Y \subset Y)$, there is a
one-to-one corespondence between the set of refinements $f_{\Delta}:\Delta_X
\rightarrow \Delta_Y$ preserving the integral structures and the set of toroidal
embeddings mapping proper birationally onto $(U_Y \subset Y)$ by toroidal
morphisms $f:(U_X \subset X) \rightarrow (U_Y \subset Y)$.
\endproclaim

\demo{Proof}\enddemo For a proof, we refer the reader to
[Kempf-Knudsen-Mumford-SaintDonat] and [Abramovich-Karu].  We only note that a
proper birational toroidal morphism between toroidal embeddings without
self-intersection is always allowable in the sense of
[Kempf-Knudsen-Mumford-SaintDonat]. 

\vskip.1in

We can reformulate via the lemma our main theorem of this section in terms of
the conical complexes (which are always assumed to be finite in this section).

\proclaim{Theorem 8.5} Let $f_{\Delta}:\Delta' \rightarrow \Delta$ be a map
between two nonsingular conical complexes, which represents a
refinement preserving the integral structure.  Then there exist a nonsingular
conical complex $\Delta''$ obtained both from $\Delta'$ and from
$\Delta$ by some sequences of smooth star subdivisions which factor
$f_{\Delta}$
$$\Delta \leftarrow \Delta'' \rightarrow \Delta.$$
\endproclaim

Given a conical comlex $\Delta$, we consider the space $N^S \oplus
{\Bbb Z}$, for each $N^S = N_{\sigma^S}$ associated to the cone ${\sigma^S}
\in \Delta$, which can be glued together naturally via the glueing of $N^S$
to form the integral structure.  We denote this space $N_{\Delta} \oplus
{\Bbb Z}$.  By considering the spaces $(N^S \oplus
{\Bbb Z}) \otimes {\Bbb Q}$ and glueing them together, we obtain the space
$$(N_{\Delta})_{\Bbb Q}^+ = (N_{\Delta} \oplus {\Bbb Z}) \otimes {\Bbb Q} =
(N_{\Delta})_{\Bbb Q} \oplus {\Bbb Q}$$
with the lattices $N^S \oplus
{\Bbb Z}$ also glued together to form the integral structure $N_{\Delta}
\oplus {\Bbb Z}$. 

If $f_{\Delta}:\Delta' \rightarrow \Delta$ is a refinement of $\Delta$
preserving the integral structures, then we can identify $(N_{\Delta'})_{\Bbb
Q}^+$ with $(N_{\Delta})_{\Bbb Q}^+$ having the same integral structure
$N_{\Delta'} \oplus {\Bbb Z} = N_{\Delta} \oplus {\Bbb Z}$.

Observe that as in the case of toric fans we can define a cobordism $\Sigma$
in the space $(N_{\Delta})_{\Bbb Q}^+$ between $\Delta'$ and $\Delta$ as
well as the notions of collapsibility, $\pi$-nonsingularity, pointing up, etc.

Once this is understood, we can carry out the same strategy as
the one presented in \S 1 $\sim$ \S 7 by Morelli to factor a proper birational
toroidal morphism and we only have to prove:

\proclaim{Theorem 8.6} Let $f_{\Delta}:\Delta' \rightarrow \Delta$ be a map
between two nonsingular conical complexes, which represents a
refinement preserving the integral structure.  Then there exists a simplicial,
collapsible and $\pi$-nonsingular cobodism $\Sigma$ in
$(N_{\Delta})_{\Bbb Q}^+$ betwwen
conical complexes $\Delta''$ and $\Delta$ such that $\Delta''$ is obtained from
$\Delta'$ by a sequence of smooth star subdivisions and that $\Sigma$ consists
only of pointing up circuits and hence $\Delta''$ is also obtained from $\Delta$
by a sequence of smooth star subdivisions.
\endproclaim

\demo{Proof}\enddemo We follow exactly the line of argument developed in the
previous sections.

First we claim that there exists a simplicial and collapsible cobordism
$\Sigma$ between $\Delta$ and $\Delta'$.  Recall that in order to construct a
cobordism and make it collapsible in the argument for the toric case we have
utilized such global theorems as Sumihiro's and Moishezon's, which are no
longer applicable in the toroidal case.  This calamity can be avoided by
using the following simple lemma.

\proclaim{Lemma 8.7} Let $\Delta$ be a simplicial conical complex.  Then we
can embed the barycentric star subdivision $\Delta_B$ (cf. Definition 2.1)
into a toric fan
$\Delta_B^T$ in some vector space $N_{\Bbb Q}$, i.e., there is a bijective map
$i:|\Delta_B| \rightarrow |\Delta_B^T|$ such that it restricts to a linear
isomorphism to each cone $i:\sigma \rightarrow \sigma^T$.  (Note that we do
NOT require $i$ to preserve the integral structure.)
\endproclaim 

\demo{Proof}\enddemo We prove by induction on the dimension $d$ of $\Delta$
and the number of the cones of the maximal dimension $d$.

When $d = 1$, i.e., $\Delta$ is a finite number of lines, the assertion is
obvious.

Suppose the assertion is proved already for a simplicial conical complex of
either dimension $< d$ or dimension $d$ with $k -1$ number of the cones of
the maximal dimension $d$.  Take a simplicial conical complex $\Delta$ of
dimension
$d$ with $k$ number of the cones of the maximal dimesion $d$.  Choose one cone
$\sigma$ of dimension $d$ and let $\Delta_{\sigma} = \Delta - \{\sigma\}$.  By
the induction hypothesis, we can embed the barycentric star subdivision
$(\Delta_{\sigma})_B$ into a toric fan $(\Delta_{\sigma})_B^T$ in some vector
space $N'_{\Bbb Q}$
$$i':|(\Delta_{\sigma})_B| \overset{\sim}\to{\rightarrow}
|(\Delta_{\sigma})_B^T|.$$ 
We take $N_{\Bbb Q} =
N'_{\Bbb Q} \oplus {\Bbb Q}$ and regard $N'_{\Bbb Q} = N'_{\Bbb Q} \oplus
\{0\} \subset N_{\Bbb Q}$.  We only have to take the embedding $i:\Delta_B
\rightarrow \Delta_B^T$ to be the one such that
$$i|_{(\Delta_{\sigma})_B} = i':|(\Delta_{\sigma})_B| \rightarrow
|(\Delta_{\sigma})_B^T|
\subset N'_{\Bbb Q}
\subset N_{\Bbb Q} \text{\ and\ }i(r(\sigma)) = (0,1) \in N'_{\Bbb Q} \oplus
{\Bbb Q},$$ where $r(\sigma)$ is the barycenter of $\sigma$ (in the sense of
Definition 2.1 and hence corresponding to the sum of the primitive vectors of
the extremal rays for
$\sigma$) and the map
$i$ on the cones in $\Delta_B$ containing $r(\sigma)$ is defined in the obvious
way.

\vskip.1in

We resume the proof of Theorem 8.6.

\vskip.1in

Take the barycentric star subdivisions $\Delta'_B$ and $\Delta_B$ of the
conical complexes $\Delta'$ and $\Delta$, respectively, and let ${\tilde
\Delta}_B$ be a simplicial common refinement of $\Delta'_B$ and $\Delta_B$. 
By Lemma 8.7 we can embed $\Delta'_B$ into a toric fan ${\Delta'}_B^T$ in some
vector space $N'_{\Bbb Q}$.  As ${\tilde \Delta}_B$ is a refinement of
$\Delta'_B$, it can also be embedded as a toric fan ${\tilde \Delta}_B^T$ in
the same space $N'_{\Bbb Q}$ by the extension of the same map.  We can take
${\Delta^{\circ}}_B^T$, obtained by a sequence of star subdivisions from
${\Delta'}_B^T$ such that it is a refinement of ${\tilde \Delta}_B^T$ (cf.
[DeConcini-Procesi]).  By replacing the original ${\tilde \Delta}_B$ with
the pull-back of ${\Delta^{\circ}}_B^T$, we may assume that ${\tilde
\Delta}_B$ is a refinement of $\Delta_B$ and $\Delta'_B$ and that ${\tilde
\Delta}_B$ is obtained from $\Delta'_B$ by a sequence of star subdivisions.

By Lemma 8.7 we can embed $\Delta_B$ into a toric fan $\Delta_B^T$ in some
vector space $N_{\Bbb Q}$.  As ${\tilde \Delta}_B$ is a refinement of
$\Delta_B$, it can also be embedded as a toric fan ${\tilde \Delta}_B^T$ in
the same space by the extension of the same map.  Now we can apply the arguments
in \S 3 and \S 4 to conclude there is a simplicial and collapsible cobordism in
$N_{\Bbb Q}^+$ between ${\widehat {(\Delta_B^T)}}$ and ${\widehat {({\tilde
\Delta}_B^T)}}$, where ${\widehat {(\Delta_B^T)}}$ is obtained from $\Delta_B^T$
by a sequence of star subdivisions and ${\widehat {({\tilde \Delta}_B^T)}}$ is
obtained from
${\tilde \Delta}_B^T$ by another sequence of star subdivisions.  We can pull
back this cobordism to obtain a simplicial and collapsible cobordism
${\tilde \Sigma}$ in $(N_{\Delta})_{\Bbb
Q}^+$ between
${\widehat {(\Delta_B)}}$ and
${\widehat {({\tilde
\Delta}_B)}}$, where ${\widehat {(\Delta_B)}}$ is obtained from $\Delta$ by a
sequence of star subdivisions (via the barycentric star subdivision $\Delta_B$)
and
${\widehat {({\tilde
\Delta}_B)}}$ is obtained from $\Delta'$ by a sequence of star subdivisions (via
the barycentric satr subdivision $\Delta'_B$ and ${\tilde \Delta}_B$).  Now we
apply Proposition 4.8, which is also valid in the toroidal case, to the lower
face $\partial_-{\tilde \Sigma}$ and to the upper face $\partial_+{\tilde
\Sigma}$ to extend it to a simplicial and collapsible cobordism $\Sigma$
between $\Delta$ and $\Delta'$.

Now apply the process of $\pi$-desingularization described in \S 5, which is
word for word valid also in the toroidal case to make $\Sigma$ a simplicial,
collapsible and $\pi$-nonsingular cobordism between $\Delta$ and $\Delta'$.

Finally apply the process described in \S 7, which is again word for word valid
in the toroidal case, to the cobordism above to obtain the desired simplicial,
collapsible and
$\pi$-nonsingular cobordism $\Sigma$ between $\Delta''$ and $\Delta$ such that
$\Delta''$ is obtained from $\Delta'$ by a sequnce of smooth star subdivisions
and that
$\Sigma$ consists only of pointing up circuits and hence $\Delta''$ is also
obtained from $\Delta$ by a sequence of smooth star subdivisions.

\vskip.1in

This completes the proof of Theorem 8.6 and the verification of the strong
factorization theorem for proper birational toroidal morphisms.

\newpage

$$\bold{REFERENCES}$$

\ref \by [Abramovich-Karu] D. Abramovich and K. Karu
\paper Weak semistable reduction in characteristic 0
\jour preprint
\yr 1997
\pages 24 pp
\endref

\ref \by [Abramovich-Karu-Matsuki-W{\l}odarczyk] D. Abramovich, K. Karu, K.
Matsuki and J. W{\l}odarczyk
\paper Torification and factorization of birational maps
\jour preprint
\yr 1999
\pages 29 pp
\endref

\ref \by [Christensen] C. Christensen
\paper Strong domination / weak factorization of three dimensional regular
local rings
\jour J. Indian Math. Soc.
\vol 45
\yr 1981
\pages 21-47
\endref

\ref \by [Corti] A. Corti
\paper Factorizing birational maps of 3-folds after Sarkisov
\jour J. Alg. Geom. 
\vol 4
\yr 1995
\pages 23-254
\endref

\ref \by [Cutkosky1] S.D. Cutkosky
\paper Local factorization of birational maps
\jour Adv. in Math. 
\vol 132
\yr 1997
\pages 167-315
\endref

\ref \by [Cutkosky2] S.D. Cutkosky
\paper Local factorization and monomialization of morphisms
\jour preprint
\yr 1997
\pages 133 pp
\endref

\ref \by [Cutkosky3] S.D. Cutkosky
\paper Local factorization and monomialization of morphisms
\jour math.AG/9803078
\yr March 17, 1998
\pages 141 pp
\endref

\ref \by [Danilov1] V.I. Danilov
\paper The birational geometry of toric varieties
\jour Russian Math. Surveys
\vol 33
\yr 1978
\pages 97-154
\endref

\ref \by [Danilov2] V.I. Danilov 
\paper The birational geometry of toric 3-folds
\yr 1983  
\vol 21 
\jour Math.
USSR-Izv. 
\pages 269-280
\endref

\ref \by [DeConcini-Procesi] C. De Concini and C. Procesi 
\paper Complete Symmetric Varieties II \yr 1985   \jour in Algebraic Groups
and Related Topics (R. Hotta, ed.) Adv. Studies in Pure Math. \vol 6 \pages
481-513 \publ Kinokuniya, Tokyo and North Holland, Amsterdam, New York, Oxford
\endref

\ref \by [DeJong] A.J. de Jong
\paper Smoothness, semistability, and alterations
\jour Publ. Math. I.H.E.S. \vol 83
\yr 1996
\pages 51-93
\endref

\ref \by [Ewald] G. Ewald
\paper Blowups of smooth toric 3-varieties
\jour Abh. Math. Sem. Univ. Hamburg
\vol 57
\yr 1987
\pages 193-201
\endref

\ref \by [Fulton] W. Fulton
\paper Introduction to Toric varieties
\jour Ann. of Math. Stud. \vol 131
\publ Princeton University Press
\yr 1993
\endref

\ref \by [Iitaka] S. Iitaka
\paper Algebraic Geometry (An Introduction to Birational Geometry of Algebraic
Varieties)
\jour Graduate Texts in Math. \vol 76
\yr 1982
\publ Springer-Verlag
\endref

\ref \by [Kawamata1] Y. Kawamata
\paper On the finiteness of generators of a pluricanonical ring for a 3-fold
of general type
\jour Amer. J. Math.
\vol 106
\yr 1984
\pages 1503-1512
\endref

\ref \by [Kawamata2] Y. Kawamata
\paper The cone of curves of algebraic varieties
\jour Ann. of Math.
\vol 119
\yr 1984
\pages 603-633
\endref

\ref \by [Kawamata3] Y. Kawamata
\paper Crepant blowing-ups of three dimensional canonical singularities and
its application to degenerations of surfaces
\jour Ann. of Math.
\vol 127
\yr 1988
\pages 93-163
\endref

\ref \by [Kempf-Knudsen-Mumford-SaintDonat] G. Kempf, F. Knudsen, D. Mumford
and B. Saint-Donat 
\paper Toroidal Embeddings I
\yr 1973 
\vol 339
\jour Lecture Notes in Math. 
\publ Springer-Verlag 
\endref

\ref \by [King1] H. King
\paper Resolving Singularities of Maps
\jour Real algebraic geometry and topology (East Lansing, MI, 1993) 
\yr 1995
\publ Contemp. Math., Amer. Math. Soc.
\endref

\ref \by [King2] H. King
\paper A private e-mail to Morelli
\yr 1996
\endref

\ref \by [Koll\'ar] J. Koll\'ar
\paper The cone theorem.  Note to a paper: ``The cone of curves of algebraic
varieties" by Kawamata
\jour Ann. of Math.
\vol 120
\yr 1984
\pages 1-5
\endref

\ref \by [Matsuki] K. Matsuki
\paper Introduction to the Mori Program
\jour the manuscript of a textbook to
appear from Springer-Verlag
\yr 1999 
\endref

\ref \by [Morelli1] R. Morelli
\paper The birational geometry of toric varieties
\jour J. Alg. Geom. \vol 5
\yr 1996
\pages 751-782
\endref

\ref \by [Morelli2] R. Morelli
\paper Correction to ``The birational geometry of toric varieties"
\jour homepage at the Univ. of Utah
\yr 1997
\pages 767-770
\endref

\ref \by [Mori1] S. Mori
\paper Threefolds whose canonical bundles are not numerically effective
\jour Ann. of Math.
\vol 116
\yr 1982
\pages 133-176
\endref

\ref \by [Mori2] S. Mori
\paper Flip theorem and the existence of minimal models for 3-folds
\jour J. of Amer. Math. Soc.
\vol 1
\yr 1988
\pages 117-253
\endref

\ref \by [Oda1] T. Oda 
\paper Lectures on Torus Embeddings and Applications, Based on joint work with
Katsuya Miyake
\yr 1966 
\jour Tata Inst. of Fund. Research
\publ Springer-Verlag 
\vol 58
\endref

\ref \by [Oda2] T. Oda
\paper Convex Bodies and Algebraic Geometry (An Introduction to the Theory of
Toric Varieties)
\jour Ergebnisse der Mathematik und ihrer Grenzgebiete 3 Folge Band
\publ Springer-Verlag
\vol 15 
\yr 1988
\endref

\ref \by [Oda-Park] T. Oda and H. Park
\paper Linear Gale transforms and Gelfand-Kapranov-Zelvinsky decompositions
\jour Tohoku Math. J.
\vol 43
\yr 1991
\pages 375-399
\endref

\ref \by [Park] H. Park
\paper The Chow rings and GKZ decompositions for ${\Bbb Q}$-factorial toric
varieties
\jour Tohoku Math. J.
\vol 45
\yr 1993
\pages 109-145
\endref

\ref \by [Reid1] M. Reid
\paper Canonical threefolds
\jour G\'eom\'etrie Alg\'ebrique Angers
\yr 1980
\publ A. Beauville ed. Sijthoff and Nordhoff
\pages 273-310
\endref

\ref \by [Reid2] M. Reid
\paper Minimal models of canonical threefolds
\jour Adv. Stud. in Pure Math.
\vol 1
\yr 1983
\pages 131-180
\endref

\ref \by [Reid3] M. Reid
\paper Decomposition of toric morphisms
\jour Arithmetic and Geometry, papers dedicated to I. R. Shafarevich on the
occasion of his 60th birthday, vol. II, Progress in Math. (M. Artin and J.
Tate, eds.) 
\vol 36
\yr 1983
\pages 395-418
\endref

\ref \by [Reid4] M. Reid
\paper Birational geometry of 3-folds according to Sarkisov
\jour preprint
\yr 1991
\endref

\ref \by [Sarkisov] V.G. Sarkisov
\paper Birational maps of standard ${\Bbb Q}$-Fano fiberings
\jour I. V. Kurchatov Institute for Atomic Energy preprint
\yr 1989
\endref

\ref \by [Shokurov] V.V. Shokurov
\paper A non-vanishing theorem
\jour Izv. Akad. Nauk SSSR Ser. Mat.
\vol 49
\yr 1985
\pages 635-651
\endref

\ref \by [Sumihiro1] H. Sumihiro
\paper Equivariant Completion I
\jour J. Math. Kyoto Univ.
\vol 14
\yr 1974
\pages 1-28
\endref

\ref \by [Sumihiro2] H. Sumihiro
\paper Equivariant Completion II
\jour J. Math. Kyoto Univ.
\vol 15
\yr 1975
\pages 573-605
\endref

\ref \by [W{\l}odarczyk1] J. W{\l}odarczyk
\paper Decomposition of birational toric maps in blow-ups and blow-downs
\yr 1997 \jour Trans. Amer. Math. Soc. \vol 349 \page 373-411
\endref

\ref \by [W{\l}odarczyk2] J. W{\l}odarczyk
\paper Birational cobordism and factorization of birational maps
\jour math.AG/9904074
\yr 1999
\pages 23 pp
\endref

\ref \by [W{\l}odarczyk3] J. W{\l}odarczyk
\paper Combinatorial structures on toroidal varieties and a proof of the Weak
Factorization Theorem
\jour math.AG/9904076
\yr 1999
\pages 32 pp
\endref

\vskip.3in

Dan Abramovich

Department of Mathematics

Boston University

Boston, MA 02215-2411

e-mail:abrmovic$\@$math.bu.edu

\vskip.1in

Kenji Matsuki

Department of Mathematics

Purdue University

West Lafayette, IN 47907-1395

e-mail:kmatsuki$\@$math.purdue.edu

\vskip.1in

Suliman Rashid

Department of Mathematics

Purdue University

West Lafayette, IN 47907-1395

e-mail:rashid$\@$math.purdue.edu

\enddocument